\documentclass{article}
\usepackage{latexsym,amssymb,amsmath}
\usepackage{graphicx}
\usepackage{color}
\usepackage[emtex,matrix,arrow]{xy}
\newtheorem{thm}{Theorem}[section]
\newtheorem{cor}[thm]{Corollary}
\newtheorem{exs}[thm]{Examples}
\newtheorem{nota}[thm]{Notation}
\newtheorem{ex}[thm]{Example}
\newtheorem{lem}[thm]{Lemma}
\newtheorem{prop}[thm]{Proposition}
\newtheorem{defn}[thm]{Definition}
\newtheorem{rem}[thm]{Remark}
\numberwithin{equation}{section}
\newcommand{\EE}{\mathbb E}
\newcommand{\PP}{\mathbb P}
\newcommand{\R}{\mathbb R}
\newcommand{\G}{\mathbb G}

\newcommand{\F}{\mathbb{F}}

\newcommand{\OO}{\mathbb{O}}
\newcommand{\I}{\mathbb{I}}
\newcommand{\Integer}{\mathbb Z}
\newcommand{\Natural}{\mathbb N}

\newcommand{\Complex}{\mathbb C}

\newcommand{\To}{\rightarrow}

\newcommand{\Cc}{\mathcal{C}}

\newcommand{\Oo}{\mathcal{O}}

\newcommand{\Ff}{\mathcal{F}}
\newcommand{\V}{\mathcal{V}}
\newcommand{\Rr}{\mathcal{R}}
\newcommand{\Ss}{\mathcal{S}}

\newcommand{\Bgg}{\mathfrak{B}}
\newcommand{\Agg}{\mathfrak{A}}
\newcommand{\Cgg}{\mathfrak{C}}

\newcommand{\Dgg}{\mathfrak{D}}
\newcommand{\mgg}{\mathfrak{m}}
\newcommand{\ngg}{\mathfrak{n}}
\newcommand{\id}{\rm{id}}
\newcommand{\m}{\bf{Meas}}
\newcommand{\mg}{\mathfrak{Meas}}
\newcommand{\repg}{\mathfrak{Rep}}
\newcommand{\dv}{\mathbf{2Vect}}

\newcommand{\qed}{\hspace*{\fill}$\Box$  \ifmmode \else
    \par\addvspace\topsep\fi}
\newenvironment {proof}{\par\addvspace\topsep\noindent{\it Proof.}
    \ignorespaces }{\qed}

\begin{document}

\title{Representation theory of 2-groups on Kapranov and Voevodsky's 2-vector spaces }
\author{Josep Elgueta \\ Dept. Matem\`atica Aplicada II \\ Universitat
  Polit\`ecnica de Catalunya \\ email: Josep.Elgueta@upc.edu}


\maketitle

\begin{abstract}
In this paper the 2-category $\repg_{{\bf 2Mat}_{\Complex}}(\G)$ of (weak) representations of an arbitrary (weak) 2-group $\G$ on (some
version of) Kapranov and Voevodsky's 2-category of (complex) 2-vector spaces is studied. In particular, the set of equivalence classes of representations is computed in terms of the invariants $\pi_0(\G)$, $\pi_1(\G)$ and $[\alpha]\in H^3(\pi_0(\G),\pi_1(\G))$ classifying $\G$. Also the categories of morphisms (up to equivalence) and the composition functors are determined explicitly. As a consequence, we obtain that the {\it monoidal} category of linear representations (more generally, the category of $[z]$-projective representations, for any given cohomology class $[z]\in H^2(\pi_0(\G),\Complex^*)$) of the first homotopy group $\pi_0(\G)$ as well as its category of representations on finite sets both live in $\repg_{{\bf 2Mat}_{\Complex}}(\G)$, the first as the monoidal category of endomorphisms of the trivial representation (more generally, as the category of morphisms between suitable 1-dimensional representations) and the second as a subcategory of the homotopy category of $\repg_{{\bf 2Mat}_{\Complex}}(\G)$.
\end{abstract}




\section{Introduction}

\label{introduccio} 

Representation theory of groups pervaded many areas in the last century
mathematics and 
theoretical physics. At the same time, among some authors, there is the
feeling that the so called {\it higher
dimensional algebra} (i.e., algebra in the category or, more generally,
$n$-category setting) can play a significant role in the mathematics
and theoretical physics of the new century. Proof of this is the two weeks
workshop held at the IMA (Minneapolis, USA) in June 2004. Therefore, it
seems natural to develop a representation theory for {\it  
  higher dimensional groups} (i.e., group structures in the $n$-category
setting, for $n\geq 1$) and to explore its possible applications.

In this paper, we consider the representation theory of 2-groups on (a
suitable   
version of) the 2-category of (complex) finite dimensional 2-vector
spaces as defined by Kapranov and Voevodsky \cite{KV94}.  
 
By a {\sl 2-group}, also called {\it gr-category} in \cite{lB92},
\cite{aR03} or {\it categorical group} in \cite{JS93}, \cite{CC96}, we
mean the suitably weakened version of the notion of a group object
(Maclane \cite{sM98}) in the category {\bf Cat} of (small) categories
and functors. Thus, roughly, it is a category $\G$ equipped with a
product functor $m:\G\times\G\To\G$, a distinguished unit 
object $I\in|\G|$ and a functor $\iota:\G\To\G$ giving inverses,
all these data satisfying the usual axioms of a group up to suitable
isomorphisms. A well known result establishes that these objects are
classified, up to the suitable notion of equivalence, by its {\it
  homotopy groups} $\pi_0(\G)$ and $\pi_1(\G)$ (the second one is in
fact a $\pi_0(\G)$-module) together with a cohomology class
$[\alpha]\in H^3(\pi_0(\G),\pi_1(\G))$ (see Section~\ref{2grups}).  

Generally, elements of a group are represented as automorphisms of some object
in some category, mostly the category ${\bf Vect}_K$ of (finite
dimensional) vector spaces 
over a field $K$. Similarly, objects in a 2-group are to be represented as
1-automorphisms (more generally, 1-autoequivalences) of some object in a given 2-category. Therefore, the first thing to be decided when studying representations of
2-groups is on which 2-category we wish to represent our 2-group.

A natural choice seems to be Kapranov and Voevodsky's 2-category
${\bf 2Vect}_K$ of finite dimensional 2-vector spaces over $K$ \cite{KV94}, as it looks
like the analog of ${\bf Vect}_K$ in the 2-category setting. {\it Strict} representation theory of {\it strict} 2-groups on this 2-category has been explored by
Barrett and Mackaay \cite{BM04}, who noticed that the resulting
representation theory is poor except when the set $\G_0$ of objects,
which has a group structure when $\G$ is strict, is a profinite group,
excluding interesting cases of Lie groups. More concretely, they
observed that, except for such a choice of 
2-groups, the representations can not be collectivelly faithful with respect
to $\G_0$ (as regards this, see Remark~\ref{remarca_fidelitat_representacions} below) \footnote{In the finally available version of their work \cite{BM04}, this result is not mentioned. But it is mentioned by Crane and Yetter in the introduction to \cite{CY03}. }. 
 
This led Yetter to introduce a more involved candidate of 2-category
on which to represent 2-groups, the 2-category of {\it
measurable categories}, denoted ${\mg}$ \cite{dY03}. This 2-category is a
two-fold generalization of ${\bf 2Vect}_K$. 
On the one hand, while objects in ${\bf 2Vect}_K$ are essentially
products of the category ${\bf Vect}_K$ of {\it finite dimensional}
vector spaces, the basic building block to get objects in ${\mg}$
is the category ${\bf Hilb}_s$ of {\it possibly infinite
dimensional} separable Hilbert spaces. On the other hand, and most
important, together with the {\it finite-dimensional} 2-vector
spaces of ${\bf 2Vect}_K$, which are products of a {\it finite} number
of copies of ${\bf Vect}_K$, ${\mg}$ further includes as objects
suitably defined ``continuous products'' of ${\bf Hilb}_s$. Thus,
an arbitrary object in ${\mg}$ is a category ${\m}(X)$, for $X$
some measurable space, whose objects are (measurable) fields of
Hilbert spaces over $X$. One recovers the (Hilbert space version of
the) categories ${\bf Vect}^n$ (the objects of ${\dv}_K$) by taking
$X=\{1,\ldots,n\}$ with the 
discrete measurable structure. First attempts to determine how the
corresponding representation theory looks like have been carried
out by Crane and Yetter \cite{CY03} and Crane and Sheppeard
\cite{CS03}. But, althought the theory seems to be interesting
and richer than the representation theory on ${\bf 2Vect}_K$, it is
also more cumbersome, involving some non trivial measure theoretic notions
such as the direct integral of measurable fields of Hilbert spaces on a
measure space and the theory of disintegration of measures.

Another alternative, not explored to our knowledge, is the 2-category
of 2-vector spaces as defined by Baez and Crans \cite{BC03}. The
notion of 2-vector space introduced by 
these authors is different from that of Kapranov and Voevodsky and it
turns out to be equivalent to a 2-term chain complex of vector spaces,
so that the corresponding representation 
theory seems to be more accessible than the previous one.

In this paper, however, we wish to consider the representation theory on
Kapranov and Voevodsky's 2-vector spaces because the resulting theory
also has its interest and it is not so poor, in our opinion. For
instance, we shall see that the {\it monoidal} category of linear
representations of $\pi_0(\G)$ arises as the category of {\it
  endomorphisms} of the trivial representation of $\G$. Moreover,
current investigations \cite{jE5} seem to suggest that the whole
2-group can be reconstructed (up to equivalence, of course) from its
2-category of representations on Kapranov and Voevodsky's 2-vector
spaces, at least when the 2-group $\G$ is finite (i.e., such that
$\pi_0(\G)$ and $\pi_1(\G)$ are finite), even when the classifying
cohomology class $[\alpha]$ is non trivial. 

In this work we do not undertake a systematic development of the
representation theory of 2-groups on Kapranov and Voevodsky's 2-vector
spaces. Thus, we do not consider analogs of the notions of direct sum,
irreducibility and character of an arbitrary representation and do 
not try to state and prove results similar to those established in classical
linear representation theory of groups. We do not consider either the
important question of the monoidal structure on the obtained 2-category of
representations (some work in this direction appears in the above
mentioned work by Barrett and Mackaay \cite{BM04}). Instead, this
paper is intended as a first step towards this representation
theory. The main results we get can be stated as follows.

\vspace{0.4 truecm}
{\bf Theorem.} {\it There is a 1-1 correspondence between the
set of equivalence classes of (non zero) representations (on
Kapranov and Voevodsky's 2-vector spaces) of a 2-group
$\G$, with classifying 3-cocycle $\alpha\in Z^3(\pi_0(\G),\pi_1(\G))$,
and the set of equivalence 
classes of quadruples $(n,\rho,\beta,c)$, with $n\geq
1$, $\rho:\pi_0(\G)\To S_n$ a morphism of groups,
$\beta:\pi_1(\G)\To(\Complex^*)^n_{\rho}$ a 
morphism of $\pi_0(\G)$-modules such that $[\beta\circ\alpha]=0$ in
$H^3(\pi_0(\G),(\Complex^*)^n_{\rho})$ and $c\in
C^2(\pi_0(\G),(\Complex^*)^n_{\rho})$ a normalized 2-cochain such that
$\partial c=\beta\circ\alpha$, two such quadruples $(n,\rho,\beta,c)$
and $(n',\rho',\beta',c')$ being equivalent if $n=n'$ and if there
exists a permutation $\sigma\in S_n$ such that
$\rho'=\sigma\rho\sigma^{-1}$, 
$\beta'=\sigma\cdot\beta$ and $[c']=[\sigma\cdot c]$ in
$\widetilde{H}^2(\pi_0(\G),(\Complex^*)^n_{\rho'})$ ($\widetilde{H}^n$
denotes the group of $n$-cochains modulo $n$-coboundaries).}

\vspace{0.4 truecm}
Note that for 2-groups $\G$ such that $[\alpha]=0$ (called {\it split} 2-groups
in this work), the datum $c$ above reduces to a 2-cocycle of
$\pi_0(\G)$ with values in $(\Complex^*)^n_{\rho}$. It also follows
that 1-dimensional representations of any $\G$ are
described by pairs $(\beta,c)$, with $\beta,c$ as above. 

\vspace{0.4 truecm}
{\bf Theorem.} {\it For any 2-group $\G$ and pairs $(\beta,c)$ and
  $(\beta',c')$ as above, the category of morphisms ${\bf
  Hom}(\I_{\beta,c},\I_{\beta',c'})$ is terminal if $\beta\neq\beta'$.
  Otherwise, it is equivalent to the category of projective
  representations of $\pi_0(\G)$ with cohomology class 
$[c'-c]\in H^2(\pi_0(\G),\Complex^*)$. In particular, for any
  $(\beta,c)$, we have an equivalence of {\sl monoidal} categories
  ${\bf End}(\I_{\beta,c})\simeq {\bf Rep}_{{\bf
      Vect}_{\Complex}}(\pi_0(\G))$, where ${\bf Rep}_{{\bf
      Vect}_{\Complex}}(\pi_0(\G))$ denotes the monoidal category of
  finite dimensional linear representations of $\pi_0(\G)$.      }
\vspace{0.4 truecm}

Here, the monoidal structure on ${\bf End}(\I_{\beta,c})$ is that
coming from the composition of 1-arrows and the horizontal composition
of 2-arrows. Observe that this category will inherit an additional monoidal
structure from the monoidal structure on the 2-category of
representations of $\G$. But both are likely to be compatible and
hence, by a categorical Eckmann-Hilton argument, equivalent.

\vspace{0.4 truecm}
{\bf Theorem.} {\it Let $\G$ be any 2-group. Then, the category ${\bf
Rep}_{{\bf FinSets}}(\pi_0(\G))$ of representations on finite
sets of its first homotopy group $ \pi_0(\G)$ is equivalent to a
subcategory of ${\sf Ho}(\repg_{{\bf 2Vect}_{\Complex}}(\G))$, the
homotopy category of the 2-category of representations of $\G$ on Kapranov and
Voevodsky's 2-vector spaces.
}

\vspace{0.4 truecm}
Let us mention some reasons by which
one can get interested in studying the representation theory of
2-groups. It is well known that, for any (reasonably well-behaved)
path-connected topological space $X$, there is an equivalence of
categories between the category of locally constant sheaves of sets on
$X$ and the category of (left) $\pi_1(X)$-sets. The equivalence
is given by the so called {\it monodromy functor}. Actually, this
turns out to be true for locally constant sheaves of objects in any category
$\Cc$ and the category of representations of $\pi_1(X)$ on $\Cc$. The
point is that this result has recently been extended by Polesello and
Waschkies \cite{PW04} to the 2-category setting. Thus, given any
2-category $\Cgg$, if $\repg_{\Cgg}(\Pi_2(X))$ denotes the 2-category
of representations on $\Cgg$ of the fundamental 2-group of $X$
(defined in \S~\ref{definicio_exemples}), these authors prove that
$\repg_{\Cgg}(\Pi_2(X))$ is equivalent (as a 2-category) to the
2-category of locally constant stacks on $X$ with values in
$\Cgg$. Since any 2-group can be shown to be equivalent to the
fundamental 2-group of some topological space, this gives a
geometric interpretation of the representation theory for 2-groups.

A second motivation, which was in fact our original motivation, is the
possibility that monoidal 
2-categories of representations of 2-groups, or suitable
deformations of them as defined in \cite{jE1}, may give rise to new
4-manifold invariants, in the same way as
suitable deformations of the categories of representations of some
Lie groups have been shown to give rise to new invariants of
3-manifolds, the so-called {\it quantum invariants} (see, or ex.,
\cite{vT94}). This is part of an old program which started with the
work by Crane and Frenkel \cite{CF94} on the so-called {\it
Hopf categories}, introduced as suitable algebraic structures
from which four-dimensional topological quantum field theories
could be constructed.

Finally, there are also reasons coming from mathematical physics. More
concretely, from attempts to construct a quantum theory of gravity.
Thus, it has been suggested \cite{CS03} that
``the fundamental symmetry to use to construct quantum gravity is the
Poincar\'e group action, but with 
the translation subgroup differentiated from the Lorentz group''.
This naturally leads to the search of a new representation theory
which respects this decomposition. Hopefully, while this
information is lost when treating the Poincar\'e group simply as a
group, this is not so when we think of it as a 2-group with
$\pi_0(\G)={\rm SO}(3,1)$ and $\pi_1(\G)=\R^4$ (see Section
\ref{2grups}). 

The outline of the paper is as follows. In Section 2 we recall the
necessary definitions from 2-category theory, mainly the notion of
pseudofunctor bicategory. Section 3 is devoted to reviewing basic
results in 2-group theory. In particular, we give an explicit
description of the isomorphism classes of morphisms between any two
given 2-groups in a way suitable for our purposes. In Section 4 we
give the 
general definition of the 2-category of representations of a
2-group on any 2-category $\Cgg$ and we unpack the corresponding
notions of representation, 1-intertwiner and 2-intertwiner and give a
description of its objects up to isomorphism. In Section 5, the two
notions of 2-vector space mentioned above are recalled and compared
and we introduce the 2-category of 2-matrices, the version of
${\dv}^{KV}_K$ on which we are going to represent 
2-groups. This version, denoted by ${\bf 2Mat}_K$, has been defined in
\cite{jE3} and here we just recall its definition. Finally, in Section
6, which constitutes the core of the paper, we consider in detail the
case $\Cgg={\bf 2Mat}_{\Complex}$. In particular, we relate the
resulting representation theory for an arbitrary 2-group to classical
representation theories of its first homotopy group $\pi_0(\G)$.

\vspace{0.3 truecm}
\noindent{{\bf Notation and terminology.}} For any $n\geq 1$, we denote by
${\bf e}_1,\ldots,{\bf e}_n\in\Natural^n$ the canonical basis of
$\Complex^n$ and by ${\bf I}_n$ the identity matrix of order $n$.
Given any $m\times n$ complex matrix {\bf R} and any vector ${\bf
x}\in\Complex^n$, ${\bf R}({\bf x})\in\Complex^m$ stands for the
usual action of {\bf R} on {\bf x}, and
${\bf R}_k$ denotes the vector in $\Complex^n$ with components the
entries in the $k^{th}$ row of {\bf R}. $\Complex^*$ denotes the
multiplicative group of the non zero complex numbers and $S_n$, $n\geq
1$, the symmetric group on $n$ letters. For any morphism of groups
$\rho:G\To G'$ and $G'$-module $M$, we write $M_{\rho}$ to denote the
abelian group $M$ equipped with the $G$-module structure induced via
$\rho$ by its $G'$-module structure.   

All over the paper, 2-category means a strict 2-category, while we use
the term bicategory for the weak version. We assume the reader is
familiar with both of these notions (see \cite{fB94}) as well as with
their one object versions, the strict and weak monoidal categories
(see \cite{sM98}, \cite{SR72}). Vertical and horizontal composition
laws for 2-morphisms are denoted by $\cdot$ and $\circ$, respectively,
and the associativity and unit constraints by $a,l,r$. We adopt the
convention that the associator $a_{h,g,f}$ and the left and right unit
constraints $l_f$ and $r_f$ go from $(h\circ g)\circ f$ to $h\circ
(g\circ f)$, from ${\id}_Y\circ f$ to $f$ and from $f\circ{\id}_X$ to
$f$, respectively. {\bf Cat} denotes the category of (small) 
categories and functors while $\underline{{\bf Cat}}$ stands for the
2-category with the same objects and morphisms and the natural
transformations as 2-morphisms.



\section{Pseudofunctor bicategories}
\label{seccio_bicategories_pseudofunctor}

Given arbitrary bicategories $\Cgg$
and $\Dgg$, there is a {\it pseudofunctor bicategory} $[\Dgg,\Cgg]$ analogous the usual functor categories
(cf. \cite{sM98}, \S II.4). Its objects, 1-morphisms and 2-morphisms are respectively the
pseudofunctors $\Ff:\Dgg\rightarrow\Cgg$, the pseudonatural transformations and the
modifications. For the sake of completeness and because there are various versions for the first two notions, we recall here the definitions.

\begin{defn} \label{def_pseudofunctor}
Let $\Cgg,\Dgg$ be any bicategories. A {\sl pseudofunctor} $\Ff:\Dgg\To\Cgg$ consists of the following set of data:
\begin{itemize}
\item[{\rm (D1)}] For every object $X$ of $\Dgg$, an object $\Ff(X)$ of $\Cgg$.
\item[{\rm (D2)}] For every ordered pair $(X,Y)$ of objects in $\Dgg$, a functor $\Ff_{X,Y}:\Dgg(X,Y)\To\Cgg(\Ff(X),\Ff(Y))$ (for short, we shall write $\Ff(*)$ instead of $\Ff_{X,Y}(*)$, for $*$ any object or morphism in $\Dgg(X,Y)$).
\item[{\rm (D3)}] For every ordered triple $(X,Y,Z)$ of objects in $\Dgg$, a natural isomorphism
$$
\Ff_2^{X,Y,Z}:c^{\Cgg}_{\Ff(X),\Ff(Y),\Ff(Z)}\circ(\Ff_{X,Y}\times\Ff_{Y,Z})\Rightarrow\Ff_{X,Z}\circ c^{\Dgg}_{X,Y,Z}
$$
(the $c$'s denote composition functors in the corresponding bicategory), whose $(f,g)$-component ($f:X\To Y$, $g:Y\To Z$) is denoted by $\Ff_2(g,f):\Ff(g)\circ\Ff(f)\Rightarrow\Ff(g\circ f)$.
\item[{\rm (D4)}] For any object $X$ of $\Dgg$, a 2-isomomorphism $\Ff_0(X):\Ff({\id}_X)\Rightarrow{\id}_{\Ff(X)}$.
\end{itemize}
Moreover, these data must satisfy the following axioms:
\begin{itemize}
\item[{\rm (A1)}] For every path $X\stackrel{f}{\To} Y\stackrel{g}{\To} Z\stackrel{h}{\To} T$ in $\Dgg$, the following diagram commutes:
\begin{equation} \label{axioma_composicio_pseudofunctor}
\xymatrix{
\Ff(h)\circ(\Ff(g)\circ \Ff(f))\ar[rrr]^{1_{\Ff(h)}\circ\Ff_2(g,f)} & & & \Ff(h)\circ\Ff(g\circ f)\ar[d]^{\Ff_2(h,g\circ f)} \\ (\Ff(h)\circ\Ff(g))\circ\Ff(f)\ar[u]^{a_{\Ff(h),\Ff(g),\Ff(f)}}\ar[d]_{\Ff_2(g,f)\circ 1_{\Ff(h)}} & & & \Ff(h\circ(g\circ f)) \\ \Ff(h\circ g)\circ\Ff(f)\ar[rrr]_{\Ff_2(h\circ g,f)} & & & \Ff((h\circ g)\circ f)\ar[u]_{\Ff(a_{h,g,f})} 
}
\end{equation}
\item[{\rm (A2)}] For any 1-morphism $f:X\To Y$, the following diagrams commute:
\begin{equation} 
\xymatrix{ 
\Ff(f)\circ\Ff({\id}_X)\ar[rr]^{\Ff_2(f,{\id}_X)}\ar[d]_{1_{\Ff(f)}\circ\Ff_0(X)} & & \Ff(f\circ{\id}_X)\ar[d]^{\Ff(r_f)}  \\ \Ff(f)\circ{\id}_{\Ff(X)}\ar[rr]_{r_{\Ff(f)}} & & \Ff(f) }\quad  
\xymatrix{
\Ff({\id}_Y)\circ\Ff(f)\ar[rr]^{\Ff_2({\id}_Y,f)}\ar[d]_{\Ff_0(Y)\circ 1_{\Ff(f)}} & & \Ff({\id}_Y\circ f)\ar[d]^{\Ff(l_f)}  \\ {\id}_{\Ff(Y)}\circ\Ff(f)\ar[rr]_{l_{\Ff(f)}} & & \Ff(f) }
\label{axioma_unitat_2_pseudofunctor}
\end{equation}
\end{itemize}
The pseudofunctor is called a {\sl 2-functor} when the 2-isomorphisms $\Ff_2(g,f)$ and $\Ff_0(X)$ are all identities.
\end{defn}

\begin{defn} \label{def_transf_pseudonatural}
Let $\Cgg,\Dgg$ be bicategories and $\Ff,\Ff':\Dgg\To\Cgg$ pseudofunctors. A {\sl pseudonatural transformation} $\xi:\Ff\To\Ff'$ consists of the following set of data:
\begin{itemize}
\item[{\rm (D1)}] For every object $X$ of $\Dgg$, a 1-morphism $\xi_X:\Ff(X)\To\Ff'(X)$ in $\Cgg$.
\item[{\rm (D2)}] For every ordered pair $(X,Y)$ of objects in $\Dgg$, a natural isomorphism
$$
\xi_2^{X,Y}:c^{\Cgg}_{\Ff(X),\Ff'(X),\Ff'(Y)}(\xi_X,-)\circ\Ff'_{X,Y}\Rightarrow c^{\Cgg}_{\Ff(X),\Ff(Y),\Ff'(Y)}(-,\xi_Y)\circ\Ff_{X,Y}
$$
whose $f$-component ($f:X\To Y$) is denoted by $\xi_2(f):\Ff'(f)\circ\xi_X\Rightarrow\xi_Y\circ\Ff(f)$.
\end{itemize}
Moreover, these data must satisfy the following axioms:
\begin{itemize}
\item[{\rm (A1)}] For every path $X\stackrel{f}{\To} Y\stackrel{g}{\To} Z$ in $\Dgg$, the following diagram commutes:
\begin{equation} \label{axioma_composicio_trans_pseudo}
\xymatrix{
\xi_Z\circ(\Ff(g)\circ\Ff(f))\ar[rrr]^{1_{\xi_Z}\circ\Ff_2(g,f)} & & & \xi_Z\circ\Ff(g\circ f) \\ (\xi_Z\circ\Ff(g))\circ\Ff(f)\ar[u]^{a_{\xi_Z,\Ff(g),\Ff(f)}} & & &  \Ff'(g\circ f)\circ\xi_X\ar[u]_{\xi_2(g\circ f)} \\ (\Ff'(g)\circ\xi_Y)\circ\Ff(f)\ar[u]^{\xi_2(g)\circ 1_{\Ff(f)}} & & & (\Ff'(g)\circ\Ff'(f))\circ\xi_X\ar[u]_{\Ff'_2(g,f)\circ 1_{\xi_X}}\ \\ \Ff'(g)\circ(\xi_Y\circ\Ff(f))\ar[u]^{a^{-1}_{\Ff'(g),\xi_Y,\Ff(f)}} & & & \Ff'(g)\circ(\Ff'(f)\circ\xi_X)\ar[lll]^{1_{\Ff'(g)}\circ\xi_2(f)}\ar[u]_{a^{-1}_{\Ff'(g),\Ff'(f),\xi_X}} 
}
\end{equation}
\item[{\rm (A2)}] For any object $X$ of $\Dgg$, the following diagram commutes:
\begin{equation} \label{axioma_unitat_trans_pseudo}
\xymatrix{
\xi_X\circ{\id}_{\Ff(X)}\ar[rr]^{l^{-1}_{\xi_X}\circ r_{\xi_X}} & & {\id}_{\Ff'(X)}\circ\xi_X  \\ \xi_X\circ\Ff({\id}_X)\ar[u]^{1_{\xi_X}\circ\Ff_0(X)} & & \Ff'({\id}_X)\circ\xi_X \ar[ll]^{\xi_2({\id}_X)}\ar[u]_{\Ff'_0(X)\circ 1_{\xi_X}} }
\end{equation}
\end{itemize}
When all {\sl structural 2-isomorphisms} $\xi_2(f)$ are identities, $\xi$ is called a {\sl 2-natural transformation}.
\end{defn}
Given pseudonatural transformations $\xi:\Ff\To\Ff'$ and $\xi':\Ff'\To\Ff''$, with $\Ff,\Ff',\Ff'':\Dgg\To\Cgg$, the composite $\xi'\circ\xi:\Ff\To\Ff''$ is defined by
\begin{align}
(\xi'\circ\xi)_X&=\xi'_X\circ\xi_X
\label{composicio_vertical_transformacions_pseudo_1}
\\ (\xi'\circ\xi)_2(f)&=a^{-1}_{\xi'_Y,\xi_Y,\Ff(f)}\cdot(1_{\xi'_Y}\circ\xi_2(f))\cdot a_{\xi'_Y,\Ff'(f),\xi_X}\cdot(\xi_2'(f)\circ 1_{\xi_X})\cdot a^{-1}_{\Ff''(f),\xi'_X,\xi_X}
\label{composicio_vertical_transformacions_pseudo_2} 
\end{align}

\begin{defn} \label{def_modificacio}
Let $\Cgg,\Dgg$ be bicategories, $\Ff,\Ff':\Dgg\To\Cgg$ pseudofunctors and $\xi,\zeta:\Ff\To\Ff'$ pseudonatural transformations. A {\sl modification} $\ngg:\xi\To\zeta$ consists of a set of 2-morphism $\{\ngg_X:\xi_X\Rightarrow\zeta_X\}$ in $\Cgg$, indexed by the objects $X$ of $\Dgg$, satisfying the {\sl naturality axiom}
\begin{equation} \label{axioma_naturalitat_modificacio}
\zeta_2(f)\cdot(1_{\Ff'(f)}\circ\ngg_X)=(\ngg_Y\circ 1_{\Ff(f)})\cdot\xi_2(f)
\end{equation}
for all 1-morphisms $f:X\To Y$ in $\Dgg$.
\end{defn}

Modifications can be composed in two ways by using the vertical and horizontal compositions in $\Cgg$. More explicitly, given modifications $\mgg:\xi\To\zeta$,
$\ngg:\zeta\To\eta$ and $\mgg':\xi'\To\zeta'$, with
$\xi,\zeta,\eta:\Ff\To\Ff'$ and $\xi',\zeta':\Ff'\To\Ff''$, the
{\it vertical composite} $\ngg\cdot\mgg:\xi\To\eta$ is defined by
\begin{equation} \label{composicio_modificacions_1}
(\mgg'\cdot\mgg)_X=\mgg'_X\cdot\mgg_X
\end{equation}
and the {\it horizontal composite}
$\mgg'\circ\mgg:\xi'\circ\xi\Rightarrow\zeta'\circ\zeta$ by
\begin{equation} \label{composicio_modificacions_2}
(\mgg'\circ\mgg)_X=\mgg'_X\circ\mgg_X
\end{equation}
for all objects $X$ in $\Dgg$. Then, we have the following

\begin{defn}
Given arbitrary bicategories $\Cgg,\Dgg$, the {\sl pseudofunctor bicategory} $[\Dgg,\Cgg]$ is the bicategory whose objects, 1-morphisms and 2-morphisms are respectively the pseudofunctors $\Ff:\Dgg\To\Cgg$, the pseudonatural transformations and the modifications, and whose composition laws for 1- and 2-morphisms are as defined above. 
\end{defn}
Recall that two bicategories $\Dgg$ and $\Cgg$ are called {\sl biequivalent} when there exists a {\sl biequivalence} between them, i.e., a pseudofunctor $\Ff:\Dgg\To\Cgg$ which is essentially surjective (for any object $U$ of $\Cgg$, there exists an object $X$ of $\Dgg$ such that $U$ is equivalent to $\Ff(X)$) and a local equivalence (each functor $\Ff_{X,Y}$ is an equivalence of categories). The following results are well known:

\begin{prop} \label{bicat_pseudof_equivalents}
Let $\Agg$, $\Bgg$, $\Cgg$ and $\Dgg$ be any bicategories. Then:
\begin{itemize}
\item[(i)] If $\Cgg$ is a 2-category, the same is true for $[\Dgg,\Cgg]$.
\item[(ii)] Any pseudofunctor $\Ff:\Dgg\To\Cgg$ induces pseudofunctors $\Ff\circ:[\Agg,\Dgg]\To[\Agg,\Cgg]$ and $\circ\Ff:[\Dgg,\Bgg]\To[\Cgg,\Bgg]$. If $\Ff$ is a biequivalence, the same is true for $\Ff\circ$ and $\circ\Ff$. 
\end{itemize}
\end{prop}

The proof of the next result is an easy exercise left to the reader: 

\begin{lem} \label{lema_2-categoria_2-functors}
Let $\Ff,\Ff':\Dgg\To\Cgg$ be arbitrary objects in $[\Dgg,\Cgg]$. Then, a 1-morphism $\xi:\Ff\To\Ff'$ is invertible (resp., an equivalence) if and only if $\xi_X$ is invertible (resp., an equivalence) in $\Cgg$ for each object $X$ of $\Dgg$.
\end{lem} 

Note that, because of the naturality axiom in the definition of modification, it may occur that two 1-morphisms $\xi,\zeta:\Ff\To\Ff'$ in $[\Dgg,\Cgg]$ are not 2-isomorphic even when the respective components $\xi_X$ and $\zeta_X$, for any object $X$ of $\Dgg$, are 2-isomorphic in $\Cgg$ (cf. \S\ref{subseccio_categories_de_morfismes}).



\section{2-groups}

\label{2grups}

Recall that an object $A$ in a monoidal category $(\Cc,\otimes,I,a,l,r)$ is called {\sl invertible} (or {\sl 2-regular}; cf. \cite{SR72}) if both ``translation functors'' $-\otimes A,A\otimes -:\Cc\To\Cc$ are equivalences. It is shown this is equivalent to the existence of an object $\alpha$ such that $A\otimes\alpha\cong I\cong \alpha\otimes A$ (if such an object exists, it is unique up to isomorphism). The object $A$ is said to be  {\sl strictly invertible} if $A\otimes\alpha=\alpha\otimes A=I$ for some $\alpha$.  

\subsection{Definition and some examples}

\label{definicio_exemples}
There are several ways of defining the notion of 2-group, depending on the amount of structure assumed on it. We shall take as definition the less structured one, corresponding to what Baez and Lauda call a {\sl weak 2-group} (see \cite{BL03}).  

\begin{defn}
{\rm (Category point of view)}\ A {\sl 2-group} is a monoidal category
$(\G,\otimes,I,a,l,r)$ which is a groupoid (i.e., all morphisms are isomorphisms)
and such that all objects are invertible. 
\end{defn}
To simplify notation, we shall often write $\G$ instead of $(\G,\otimes,I,a,l,r)$. It is shown (\cite{BL03}, Thm. 17) that this definition is indeed equivalent to the more structured one, called by Baez and Lauda a {\sl coherent 2-group}, and obtained by weakening the notion of a group object in {\bf Cat} via suitable coherent natural isomorphisms.

As pointed out in the introduction, a monoidal category $\Cc$ is the same thing as a bicategory with only one object: the 1- and 2-morphisms are
the objects and morphisms of $\Cc$, the vertical composition of
2-morphisms corresponds to the composition of
morphisms in $\Cc$ and the composition 
of 1-morphisms and horizontal composition of 2-morphisms
is given by the tensor product in $\Cc$ (see Benabou \cite{jB67}) . The bicategory with one object so defined is denoted by $\underline{\Cc}$. This 
leads to the following alternative way of thinking of a 2-group, completely equivalent to the previous one.

\begin{defn}
{\rm (Bicategory point of view)}\ A {\sl 2-group} is a bigroupoid with only one
object, i.e., a bicategory with only one object and such that all
1-morphisms are equivalences and all 2-morphisms are isomorphisms. 
\end{defn}

The following particular types of 2-groups appear in the sequel (as discussed in \S\ref{2-grups_modul_equivalencia}, both are general enough to include all 2-groups up to the appropriate notion of equivalence).

\begin{defn}
Let $\G$ be a 2-group. Then:
\begin{itemize}
\item[(i)]
$\G$ is called {\sl strict} if the underlying monoidal category is strict (i.e., all natural isomorphisms $a,l,r$ are identities) and every object is strictly invertible.  
\item[(ii)]
$\G$ is called s{\sl pecial} if the underlying category is skeletal (i.e., isomorphic objects are equal), the natural isomorphisms $l,r$ are identities and every object is strictly invertible.
\end{itemize}
\end{defn}

\begin{nota} \label{notacio_objectes_2-grup} {\rm 
If $\G$ is skeletal (in particular, special), its objects are denoted by $g$'s, with subscripts if necessary, and the unit object by $e$; otherwise, they are denoted by capital letters $A,B,\ldots$, with $I$ for the unit object. This notation is due to the obvious fact that the set of objects for such a 2-group is a group with the tensor product. }
\end{nota}

\begin{exs} \label{primers_exemples_2-grups} {\rm 
(1) The terminal category {\bf 1} with only one object and one (identity) morphism is a strict 2-group in the obvious way. It is called the {\sl trivial 2-group}.

(2) Among the simplest examples of non trivial 2-groups there are the {\sl
discrete 2-groups}, i.e., the groups
themselves viewed as discrete categories. The tensor product
$\otimes$ is given by the group law. For any group $G$, 
the corresponding discrete 2-group is strict and is denoted by $G[0]$ (the zero
indicates that the elements of $G$ correspond to the objects in the monoidal category point of view).

(3) Any abelian group $A$ viewed as a category with only one object also 
defines a strict 2-group, denoted $A[1]$ (the elements of $A$ now correspond to morphisms). Both the
tensor product and the composition coincide and are given by the group law in A (this requires the group to be abelian because otherwise the tensor product $\otimes$ will not be functorial).
}
\end{exs}

In general, given any bicategory $\Cgg$ and any object $X$ of
$\Cgg$, the category ${\sf Equiv}_{\Cgg}(X)$ with objects the autoequivalences $f:X\To X$ and morphisms the 2-isomorphisms between these is a monoidal subcategory of $\Cgg(X,X)$. It is clearly a 2-group (in fact, any 2-group is of this sort for some bicategory $\Cgg$ and object $X$ of $\Cgg$). Moreover, the full subcategory ${\sf
  Aut}_{\Cgg}(X)\subset {\sf Equiv}_{\Cgg}(X)$ with objects only the automorphisms of $X$ (i.e., the strictly invertible 1-morphisms $f:X\To X$) is a sub-2-group, which is strict when $\Cgg$ is a
2-category. This gives some more examples of 2-groups

\begin{exs} {\rm 
(1) Let $\Cgg=\underline{{\bf Cat}}$. Then, for any category $\Cc$, we have the corresponding 2-groups of symmetries ${\sf Equiv}_{\underline{{\bf Cat}}}(\Cc)$ or ${\sf Aut}_{\underline{{\bf Cat}}}(\Cc)$, which could be called, when $\Cc$ is finite, the {\it weak} and {\it strict symmetric 2-groups}, respectively. They include the usual symmetric groups as special cases (take $\Cc$ to be any finite set viewed as a discrete category).

(2) A particularly important case is when $\Cgg=\Pi_2(X)$, the fundamental bigroupoid of some topological space $X$ (see \cite{HKK01}). This is the bicategory with the points of $X$ as objects, the paths between points as 1-morphisms and the homotopy classes of homotopies between paths as 2-morphisms. Then, for any object $x\in X$, we have the corresponding 2-group of autoequivalences, usually denoted by $\Pi_2(X,x)$ and called the {\sl fundamental 2-group} of the pointed space $(X,x)$. This example turns out to be completely general. Any 2-group is equivalent, in a sense made precise later, to the fundamental 2-group of some pointed space. This is a consequence of a fundamental theorem relating the theory of 2-groups to that of homotopy 2-types, and according to which there is a biequivalence between the 2-category of 2-groups (defined below) and a suitably defined 2-category of connected pointed homotopy 2-types (see Tamsamani \cite{zT96}).   }
\end{exs}

\subsection{2-category of 2-groups}

2-groups constitute the objects of a 2-category {\sf 2Grp}, which is a (full) sub-2-category of {\sf MonCat}, the 2-category of monoidal categories. For the sake of completeness, we recall here the precise definitions.

\begin{defn}
Let $\G$, $\G'$ be 2-groups (viewed as monoidal categories). A {\sl morphism of 2-groups} $\F:\G\To\G'$ is a pair $\F=(F,F_2)$, where $F:\G\To\G'$ is a functor and $F_2$ (the {\sl monoidal structure} of $F$) is a collection of natural isomorphisms $F_2(A,B):F(A)\otimes' F(B)\stackrel{\cong}{\To} F(A\otimes B)$, for all objects $A,B$ in $\G$, making commutative the diagrams
\begin{equation} \label{axioma_hexagon}
\xymatrix{
F(A)\otimes'(F(B)\otimes' F(C))\ar[rrr]^{{\id}_{F(A)}\otimes' F_2(B,C)} & & & F(A)\otimes' F(B\otimes C)\ar[d]^{F_2(A,B\otimes C)} \\ (F(A)\otimes' F(B))\otimes' F(C)\ar[u]^{a'_{F(A),F(B),F(C)}}\ar[d]_{F_2(A,B)\otimes'{\id}_{F(C)}} & & & F(A\otimes(B\otimes C)) \\ F(A\otimes B)\otimes' F(C)\ar[rrr]_{F_2(A\otimes B,C)} & & & F((A\otimes B)\otimes C)\ar[u]_{F(a_{A,B,C})} 
}
\end{equation} 
for all objects $A,B,C$ in $\G$. The morphism $\F$ is called {\sl strict} when each of the isomorphisms $F_2(A,B)$ is an identity and $F(r_I)=r'_{I'}$ (equivalently \footnote{Recall that in any monoidal category $(\Cc,\otimes,I,a,l,r)$, the morphisms $r_I,l_I:I\otimes I\To I$ coincide.}, $F(l_I)=l'_{I'}$).
\end{defn}

Notice that, given a morphism of 2-groups $\F=(F,F_2):\G\To\G'$, there always exists a unique isomorphism $F_0:F(I)\stackrel{\cong}{\To} I'$ such that the diagrams
\begin{equation} \label{axiomes_unitat}
\xymatrix{
F(A)\otimes' F(I)\ar[rr]^{F_2(A,I)}\ar[d]_{{\id}_{F(A)}\otimes' F_0} & & F(A\otimes I)\ar[d]^{F(r_A)}  \\ F(A)\otimes' I'\ar[rr]_{r'_{F(A)}} & & F(A) }\quad \xymatrix{
F(I)\otimes' F(A)\ar[rr]^{F_2(I,A)}\ar[d]_{F_0\otimes' {\id}_{F(A)}} & & F(I\otimes A)\ar[d]^{F(l_
A)}  \\ I'\otimes' F(A)\ar[rr]_{l'_{F(A)}} & & F(A) }
\end{equation}
commute for all objects $A$ in $\G$. Indeed, $F(I)\otimes' -$ is an equivalence and hence, there is a unique morphism $F_0:F(I)\To I'$, which is necessarily an isomorphism, making commutative the diagram on the left with $A=I$. It is shown that this isomorphism makes all previous diagrams commute. Therefore, a morphism of 2-groups indeed is the same thing as a monoidal functor between the underlying monoidal categories or, in the bicategory point of view, a pseudofunctor $\Ff:\underline{\G}\To\underline{\G}'$. Strict morphisms correspond to strict monoidal functors or to 2-functors, respectively. It may further be shown that any morphism of 2-groups preserves inverses, i.e., for any object $A$, it is $F(\alpha)\cong\overline{F(A)}$, the isomorphisms satisfying the appropriate coherent laws (cf. \cite{BL03}).

Morphisms of 2-groups can be composed as follows. If $\F=(F,F_2):\G\To\G'$ and $\F'=(F',F'_2):\G'\To\G''$, the composite morphism $\F'\circ\F:\G\To\G''$ is $\F'\circ\F=(F'\circ F, F_2'\star F_2)$, with $(F'_2\star F_2)(A,B)=F'_2(F(A),F(B))\circ F'(F_2(A,B))$ for all $A,B$.

\begin{defn}
Let $\F,\widetilde{\F}:\G\To\G'$ be morphisms of 2-groups. A {\sl monoidal natural transformation} from $\F$ to $\widetilde{\F}$ is a natural transformation $\tau:\F\Rightarrow\widetilde{\F}$ such that the diagrams
\begin{equation} \label{monoidalitat_transf_natural}
\xymatrix{
F(A)\otimes' F(B)\ar[rr]^{F_2(A,B)}\ar[d]_{\tau_A\otimes'\tau_B} & & F(A\otimes B)\ar[d]^{\tau_{A\otimes B}}  \\ \widetilde{F}(A)\otimes' \widetilde{F}(B)\ar[rr]_{\widetilde{F}_2(A,B)} & & \widetilde{F}(A\otimes B) }\quad
\xymatrix{
F(I)\ar[rr]^{\tau_I}\ar[rd]_{F_0} & & \widetilde{F}(I)\ar[ld]^{\widetilde{F}_0} \\ & I' & }
\end{equation}
commute for all objects $A,B$ of $\G$.
\end{defn}

It is immediate to check that vertical and horizontal composites of monoidal natural transformations are still monoidal. Hence, we have the following 2-category of 2-groups:

\begin{defn}
The {\sl 2-category of 2-groups}, denoted ${\sf 2Grp}$, is the 2-category with objects, 1-morphisms and 2-morphisms the 2-groups, the morphisms of 2-groups and the monoidal natural transformations between these, respectively.
\end{defn}
Observe that all 2-morphisms in this 2-category are 2-isomorphisms, because any 2-group is a groupoid by definition.

\subsection{The 2-functors $\pi_0$ and $\pi_1$}

\label{2-functors_pi0_pi1}

Let ${\bf Grp}[0]$ and ${\bf Ab}[0]$ be the categories of groups and of abelian groups, respectively, thought of as 2-categories with only identity 2-arrows. Both are sub-2-categories of {\sf 2Grp} via the 2-functors $[0]:{\bf Grp}[0]\To{\sf 2Grp}$ and $[1]:{\bf Ab}[0]\To{\sf 2Grp}$ acting on objects are described in Examples~\ref{primers_exemples_2-grups}-(2),(3). Next result seems to be due to Sinh \cite{hxS75}:

\begin{thm}
There are 2-functors $\pi_0:{\sf 2Grp}\To{\bf Grp}[0]$ and $\pi_1:{\sf 2Grp}\To{\bf Ab}[0]$.
\end{thm}
They are defined as follows. For any 2-group $\G$, $\pi_0(\G)$ is the group of isomorphism classes of objects of $\G$ (it is a group with the product given by $[A][B]=[A\otimes B]$) and $\pi_1(\G)={\rm Aut}_{\G}(I)$, the group of automorphisms of the unit object (as the endomorphism monoid of the unit object of any monoidal category, it is indeed abelian; see \cite{SR72}, Chap. I, \S 1.3.3.1). For any morphism of 2-groups $\F:\G\To\G'$, $\pi_0(\F):\pi_0(\G)\To\pi_0(\G')$ is the group morphism defined by $\pi_0(\F)([A])=[F(A)]$ and $\pi_1(\F):\pi_1(\G)\To\pi_1(\G')$ is that given by $\pi_1(\F)(u)=F_0\circ F(u)\circ F_0^{-1}$ for all $u:I\To I$. $\pi_0(\G)$ and $\pi_1(\G)$ are usually called the {\it homotopy groups} of $\G$ because of the following important case.

\begin{ex} {\rm 
If $\G$ is the fundamental 2-group of a pointed space $(X,x)$, $\pi_0(\G)$ and $\pi_1(\G)$ are equal to the first and second homotopy groups of $(X,x)$, respectively. }
\end{ex}

Since $\pi_0$ and $\pi_1$ are 2-functors, isomorphic morphisms $\F,\widetilde{\F}:\G\To\G'$, for given 2-groups $\G$ and $\G'$, are mapped to isomorphic, hence equal, morphisms in ${\bf Grp}[0]$ and ${\bf Ab}[0]$. In other words, $\pi_0(\F)$ and $\pi_1(\F)$ are isomorphism invariants of $\F$. As discussed in \S~\ref{subseccio_morfismes_modul_isomorfisme}, however, they do not suffice to completely determine the isomorphism class of $\F$.

To be coherent with the notation introduced in \ref{notacio_objectes_2-grup}, elements in $\pi_0(\G)$ are denoted by $g$'s (possibly with subscripts), althought one should keep in mind that they are isomorphism classes of objects in $\G$.

\subsection{2-groups up to equivalence}

\label{2-grups_modul_equivalencia}

As in any 2-category, two objects $\G,\G'$ of {\sf 2Grp} are said to be {\sl equivalent} when there exists an {\sl equivalence} between them, i.e., a morphism of 2-groups $\F=(F,F_2):\G\To\G'$ which is invertible up to a 2-isomorphism. It is easy to see that this is equivalent to the existence of a morphism $\F=(F,F_2):\G\To\G'$ whose underlying functor $F$ is an equivalence of categories.

The following strictification/skeletization results are used in the sequel (see \cite{BL03}, Prop. 37 and 43):

\begin{thm} \label{teorema_estrictificacio}
Let $\G$ be an arbitrary 2-group. Then:
\begin{itemize}
\item[(i)]
$\G$ is equivalent to a strict 2-group.
\item[(ii)]
$\G$ is equivalent to a special 2-group.
\end{itemize}
\end{thm}

Althought both homotopy groups $\pi_0(\G)$ and $\pi_1(\G)$ are invariants (up to equivalence) of $\G$, they do not suffice to classify $\G$. The complete classification additionally involves an action of $\pi_0(\G)$ on $\pi_1(\G)$ together with a cohomology class $[\alpha]\in H^3(\pi_0(\G),\pi_1(\G))$. It seems this also was first proved in Sinh's thesis \cite{hxS75} (a proof can also be found in \cite{BL03}, Sec. 8.3). More concretely, if $\pi_0({\sf 2Grp})$ denotes the set of equivalence classes of 2-groups, we have the following:

\begin{thm}
There is a 1-1 correspondence between $\pi_0({\sf 2Grp})$ and the set of isomorphism classes of triples $(G,M,[\alpha])$, with $G$ a group, $M$ a $G$-module and $[\alpha]\in H^3(G,M)$, two such triples $(G,M,[\alpha])$ and $(G',M',[\alpha'])$ being isomorphic iff there exists an isomorphism of groups $\varphi:G\To G'$ and an isomorphism of $G$-modules $\psi:M\To M'_{\varphi}$ such that $[\psi\circ \alpha]=[\alpha'\circ\varphi^3]\in H^3(G,M'_{\varphi})$. 
\end{thm}
If $(G,M,[\alpha])$ is a classifying triple of $\G$, any representative $\alpha$ of $[\alpha]$ is called a {\sl classifying 3-cocycle} of $\G$.

For later use, let us recall how this correspondence works (see \cite{hxS75} for more details). Suppose $\G$ is an arbitrary 2-group, with unit object $I$. For any object $A$ of $\G$, the maps $\delta_A,\gamma_A:\pi_1(\G)\To{\rm Aut}_{\G}(A)$ defined by
\begin{equation} \label{isos_delta_gamma}
\delta_A(u)=r_A\circ ({\id}_A\otimes u)\circ r_A^{-1},\quad  
\gamma_A(u)=l_A\circ (u\otimes{\id}_A)\circ l_A^{-1}\qquad \forall\ u\in\pi_1(\G)
\end{equation}
are group isomorphisms \footnote{In fact, the maps $\delta_A$ and $\gamma_A$ can be defined for any monoidal category (see Saavedra Rivano \cite{SR72}, Chap. I, \S 1.3.3.3., where they are denoted $\rho'_A$ and $\rho_A$, respectively), and they provide a transitive system of canonical isomorphisms between the endomorphism monoids of all invertible objects (see {\it loc. cit.}, Chap. I,\S 2.5.2.).}. Then, a classifying triple $(G,M,[\alpha])$ of $\G$ is obtained as follows:

\begin{itemize}
\item[(a)]
$G=\pi_0(\G)$.
\item[(b)]
$M=\pi_1(\G)$ with the $\pi_0(\G)$-module structure given by $g\cdot u=\gamma_A^{-1}(\delta_A(u))$
for any representative $A\in g$. In particular, if $\G$ is strict, this reduces to $g\cdot u={\id}_A\otimes u\otimes {\id}_{A^{-1}}$.
\item[(c)] 
To obtain a classifying 3-cocycle $\alpha\in Z^3(\pi_0(\G),\pi_1(\G))$, let us choose for each $g\in\pi_0(\G)$ a representative $A_g\in g$ with $A_e=I$, and for any other object $B\in g$, a morphism $\iota_B:B\To A_g$ with $\iota_{A_g}={\id}_{A_g}$ and $\iota_{I\otimes A_g}=l_{A_g}$ (if $I\otimes A_g\neq A_g$), $\iota_{A_g\otimes I}=r_{A_g}$ (if $A_g\otimes I\neq A_g$). This amounts to choosing a skeleton $\G_{sk}$ of $\G$ together with a specific equivalence of categories $E:\G\To\G_{sk}$, given by $E(f)=\iota_C\circ f\circ \iota^{-1}_B$ when $f:B\To C$. Then, for any $g_1,g_2,g_3\in\pi_0(\G)$, let $\widetilde{a}_{g_1,g_2,g_3}\in{\rm Aut}_{\G}(A_{g_1g_2g_3})$ be defined by
$$
\widetilde{a}_{g_1,g_2,g_3}=\iota_{A_{g_1}\otimes A_{g_2g_3}}\circ ({\id}_{A_{g_1}}\otimes\iota_{A_{g_2}\otimes A_{g_3}})\circ a_{A_{g_1},A_{g_2},A_{g_3}}\circ(\iota_{A_{g_1}\otimes A_{g_2}}^{-1}\otimes{\id}_{A_{g_3}})\circ\iota_{A_{g_1g_2}\otimes A_{g_3}}^{-1}
$$
It gives the associator of the monoidal structure on $\G_{sk}$ induced by that on $\G$ via the equivalence $E$. In particular, it satisfies the pentagon axiom. Then, a classifying 3-cocycle $\alpha$ for $\G$ is given by
\begin{equation} \label{3-cocicle}
\alpha(g_1,g_2,g_3)=\gamma^{-1}_{A_{g_1g_2g_3}}(\widetilde{a}_{g_1,g_2,g_3})\in\pi_1(\G)
\end{equation}
(the pentagon axiom on $\widetilde{a}$ together with the restrictions made in choosing the representatives $A_g$ and morphisms $\iota_B$ ensure that this is a normalized 3-cocycle and it is shown that its cohomology class is independent of the representatives $A_g$ and morphisms $\iota_B$ started with). For special 2-groups, all $\iota$'s are identities, $\widetilde{a}=a$ and (\ref{3-cocicle}) reduces to
\begin{equation} \label{a_vs_abarra}
\alpha(g_1,g_2,g_3)=\gamma^{-1}_{g_1g_2g_3}(a_{g_1,g_2,g_3})=a_{g_1,g_2,g_3}\otimes{\id}_{(g_1g_2g_3)^{-1}},\quad g_1,g_2,g_3\in\pi_0(\G)
\end{equation}
(following the convention in \ref{notacio_objectes_2-grup}, we just write $g$ instead of $A_g$).
\end{itemize}
Conversely, given a triple $(G,M,[\alpha])$, with $[\alpha]\in H^3(G,M)$, a (special) 2-group $\G$ is built as follows:
\begin{itemize}
\item
${\rm Obj}(\G)=G$.
\item
${\rm Mor}(\G)=G\times M$, a pair $(g,m)$ being a morphism $(g,m):g\To g$.
\item
Composition and tensor product are given by
\begin{align} 
(g,m')\circ(g,m)&=(g,m'+m) \label{composicio_morfismes_2-grup_associat} \\
g_1\otimes g_2&=g_1g_2 \label{producte_tensorial_objectes_2-grup_associat} \\
(g_1,m_1)\otimes(g_2,m_2)&=(g_1g_2,m_1+(g_1\cdot m_2)) \label{producte_tensorial_morfismes_2-grup_associat}
\end{align}
(in particular, the set of morphisms ${\rm Mor}(\G)$ equipped with the tensor product is equal to the semidirect product group $G\times M$).
\item
Left and right unit constraints are trivial while the associator is defined in terms of a representative $\alpha$ of $[\alpha]$ by
\begin{equation} \label{associador_2-grup_especial}
a_{g_1,g_2,g_3}=(g_1g_2g_3,\alpha(g_1,g_2,g_3))
\end{equation}
(the 3-cocycle condition on $\alpha$ ensures that $a$ satisfies the pentagon axiom and it is shown that the equivalence class of the 2-group so defined is independent of the representative $\alpha$ of $[\alpha]$ used).
\end{itemize}

\subsection{Morphisms of 2-groups up to isomorphism}

\label{subseccio_morfismes_modul_isomorfisme}
Very useful for our purposes, besides the previous description of the set $\pi_0({\sf 2Grp})$, is a description of the sets $\pi_0({\bf Hom}_{{\sf 2Grp}}(\G,\G'))$ of isomorphism classes of morphisms between any 2-groups $\G,\G'$. These can be described again in cohomological terms. The result seems to be due to Joyal and Street \cite{JS86} (cf. \cite{BL03}), althought probably not stated in this way. Explicitly, for any cochain complex of abelian groups $C^{\bullet}$, let $\widetilde{H}^{\bullet}(C^{\bullet})$ be the ``total cohomology'' of $C^{\bullet}$, i.e., $\widetilde{H}^n(C^{\bullet})=C^n/B^n$ for all $n\in\Integer$ (we take $n$-cochains, instead of $n$-cocycles, modulo $n$-coboundaries). To be coherent with the standard notation, we shall write $\widetilde{H}^{\bullet}(G,M)$ in case $C^{\bullet}$ is the complex of (normalized) cochains of a group $G$ with values in a $G$-module $M$. Then, we have the following:

\begin{thm} \label{morfismes_modul_isomorfisme}
Let $\G$ and $\G'$ be arbitrary 2-groups and $\alpha$, $\alpha'$ classifying 3-cocycles of $\G$, $\G'$ respectively. Then, there is a 1-1 correspondence between $\pi_0({\bf Hom}_{{\sf 2Grp}}(\G,\G'))$ and the set of triples $(\rho,\beta,[c])$, with $\rho:\pi_0(\G)\To\pi_0(\G')$ a morphism of groups, $\beta:\pi_1(\G)\To\pi_1(\G')_{\rho}$ a morphism of $\pi_0(\G)$-modules such that $[\beta\circ\alpha]=[\alpha'\circ\rho^3]$ in $H^3(\pi_0(\G),\pi_1(\G')_{\rho})$ and $[c]\in\widetilde{H}^2(\pi_0(\G),\pi_1(\G')_{\rho})$ such that \footnote{Here and in what follows, we use additive notation for the composition law in $\pi_1(\G')$ because this is an abelian group. But the reader should keep in mind that it is given by the composition of morphisms. Thus, the last equality actually means that, for all $g_1,g_2,g_3\in\pi_0(\G)$, $(\partial c)(g_1,g_2,g_3)$ is the automorphism of $I'$ given by $\beta(\alpha(g_1,g_2,g_3))\circ\alpha'(\rho(g_1),\rho(g_2),\rho(g_3))^{-1}$.} $\partial c=\beta\circ\alpha-\alpha'\circ\rho^3$.
\end{thm}
Because of the importance of this theorem for what follows, let us give the proof in some detail.

\begin{proof}
Note first that choosing specific classifying 3-cocycles $\alpha$ of $\G$ and $\alpha'$ of $\G'$ amounts to choosing special 2-groups $\G_{\alpha}$ and $\G'_{\alpha'}$ equivalent to $\G$ and $\G'$, respectively. Namely, those constructed from the triples $(\pi_0(\G),\pi_1(\G),\alpha)$ and $(\pi_0(\G'),\pi_1(\G'),\alpha')$ as described in \S\ref{2-grups_modul_equivalencia}. Then, for any equivalences $\EE_{\alpha}:\G\To\G_{\alpha}$ and $\EE'_{\alpha'}:\G'_{\alpha'}\To\G'$, we have a 1-1 correspondence $\pi_0({\bf Hom}_{{\sf 2Grp}}(\G_{\alpha},\G'_{\alpha'}))\To\pi_0({\bf Hom}_{{\sf 2Grp}}(\G,\G'))$ defined by $[\F]\mapsto[\EE'_{\alpha'}\circ\F\circ\EE_{\alpha}]$. We can now use the following general result concerning 1-arrows in {\sf MonCat}:

\begin{lem} \label{lema_functors_monoidals}
For $\Cc,\Cc'$ arbitrary monoidal categories, any monoidal functor $(F,F_2,F_0):\Cc\To\Cc'$ is monoidally isomorphic to a monoidal functor $(\overline{F},\overline{F}_2,\overline{F}_0):\Cc\To\Cc'$ such that $\overline{F}_0:\overline{F}(I)\To I'$ is an identity.
\end{lem}
\begin{proof}
Define $\overline{F}(X)=F(X)$ for all objects $X\neq I$ and $\overline{F}(I)=I'$, and for any morphism $f:X\To Y$, set
$$
\overline{F}(f)=\left\{ \begin{array}{ll} F(f), & \mbox{if $X,Y\neq I$} \\ F(f)\circ f_0^{-1},  & \mbox{if $X=I$, $Y\neq I$} \\ F_0\circ F(f), & \mbox{if $X\neq I$, $Y=I$} \\ F_0\circ F(f)\circ F_0^{-1}, & \mbox{if $X=Y=I$}.
\end{array}\right.
$$
Furthermore, for $X,Y\neq I$, let $\overline{F}_2(X,Y):\overline{F}(X)\otimes'\overline{F}(Y)\To\overline{F}(X\otimes Y)$ be equal to $F_2(X,Y)$ if $X\otimes Y\neq I$ or $F_0\circ F_2(X,Y)$ if $X\otimes Y=I$. Otherwise, set 
\begin{align*}
\overline{F}_2(I,Y)&=\mbox{$F_2(I,Y)\circ(F_0\otimes{\id}_{F(Y)})^{-1}$ or $F_0\circ F_2(I,Y)\circ(F_0\otimes{\id}_{F(Y)})^{-1}$} \\ \overline{F}_2(X,I)&=\mbox{$F_2(X,I)\circ({\id}_{F(X)}\otimes F_0)^{-1}$ or $F_0\circ F_2(X,I)\circ({\id}_{F(X)}\otimes F_0)^{-1}$} \\ \overline{F}_2(I,I)&=\mbox{$F_2(I,I)\circ(F_0\otimes F_0)^{-1}$ or $F_0\circ F_2(I,I)\circ(F_0\otimes F_0)^{-1}$},
\end{align*}
the two possibilities respectively correponding to the cases where the tensor product of both objects is equal to $I$ or not. Finally, set $\overline{F}_0={\id}_{I'}$. Then, it is a routine exercise left to the reader to check that this indeed defines a monoidal functor monoidally equivalent to $(F,F_2,F_0)$.  
\end{proof}

When applied to the morphisms between $\G_{\alpha}$ and $\G'_{\alpha'}$, we get a 1-1 correspondence $\pi_0({\bf Hom}^s_{{\sf 2Grp}}(\G_{\alpha},\G'_{\alpha'}))\To\pi_0({\bf Hom}_{{\sf 2Grp}}(\G_{\alpha},\G'_{\alpha'}))$, where ${\bf Hom}^s_{{\sf 2Grp}}(\G_{\alpha},\G'_{\alpha'})$ denotes the full subcategory of ${\bf Hom}_{{\sf 2Grp}}(\G_{\alpha},\G'_{\alpha'})$ with objects the so called {\sl special morphisms}, i.e., the morphisms $\F:\G_{\alpha}\To\G'_{\alpha'}$ for which $F_0$ is an identity (cf. \cite{BL03}). Next, we can use the following parametrization of the special morphisms between special 2-groups:

\begin{lem}
If $\G_s,\G'_s$ are special 2-groups and $\alpha,\alpha'$ the corresponding distinguished classifying 3-cocycles given by (\ref{a_vs_abarra}), there is a 1-1 correspondence between the set ${\rm Hom}^s_{{\sf 2Grp}}(\G_s,\G'_s)$ and the set of triples $(\rho,\beta,c)$, with $\rho:\pi_0(\G_s)\To\pi_0(\G'_s)$ a morphism of groups, $\beta:\pi_1(\G_s)\To\pi_1(\G'_s)_{\rho}$ a morphism of $\pi_0(\G_s)$-modules such that $[\beta\circ\alpha]=[\alpha'\circ\rho^3]$ in $H^3(\pi_0(\G_s),\pi_1(\G'_s)_{\rho})$ and $c\in C^2(\pi_0(\G_s),\pi_1(\G_s')_{\rho})$ a normalized 2-cochain such that $\partial c=\beta\circ\alpha-\alpha'\circ\rho^3$. Moreover, two such triples $(\rho_1,\beta_1,c_1)$ and $(\rho_2,\beta_2,c_2)$ correspond to isomorphic morphisms iff $\rho_1=\rho_2$, $\beta_1=\beta_2$ and $[c_1]=[c_2]$ in $\widetilde{H}^2(\pi_0(\G'),\pi_1(\G'_s)_{\rho})$.
\end{lem}

\begin{proof}
Given a special morphism $\F=(F,F_2):\G_s\To\G'_s$, the corresponding triple is $(\pi_0(\F),\pi_1(\F),c(\F))$, where $\pi_0(\F), \pi_1(\F)$ are as defined in \S\ref{2-functors_pi0_pi1} and $c(\F):\pi_0(\G_s)\times\pi_0(\G_s)\To\pi_1(\G'_s)$ is given by
\begin{equation} \label{2-cocadena}
c(\F)(g_1,g_2)=\gamma_{F(g_1g_2)}^{-1}(F_2(g_1,g_2))=F_2(g_1,g_2)\otimes{\id}_{F(g_1g_2)^{-1}}
\end{equation}
for all $g_1,g_2\in\pi_0(\G_s)$. The coherence condition on the maps $F_2(g_1,g_2)$ ensures that $c$ has the right coboundary. Conversely, for any triple $(\rho,\beta,c)$ as in the lemma, the corresponding special morphism $\F^{(\rho,\beta,c)}:\G_s\To\G'_s$ is that whose underlying functor $F^{(\rho,\beta,c)}$ acts on objects as $\rho$ and on morphisms by
\begin{equation} \label{accio_sobre_morfismes_F_(ro,beta,c)}
F^{(\rho,\beta,c)}(\varphi)=(\gamma_{\rho(g)}\circ\beta\circ\gamma^{-1}_g)(\varphi)=\beta(\varphi\otimes{\id}_{g^{-1}})\otimes{\id}_{\rho(g)}
\end{equation} 
for any morphism $\varphi:g\To g$ of $\G_s$, and whose monoidal structure is given by
\begin{equation} \label{estructura_monoidal_F_(ro,beta,c)}
F^{(\rho,\beta,c)}_2(g_1,g_2)=\gamma_{\rho(g_1g_2)}(c(g_1,g_2))=c(g_1,g_2)\otimes{\id}_{\rho(g_1g_2)}
\end{equation}
for all $g_1,g_2\in\pi_0(\G_s)$. It is not hard to check that this really defines a special morphism of 2-groups and that $\F^{(\pi_0(\F),\pi_1(\F),c(\F))}=\F$ and that $(\pi_0(\F^{(\rho,\beta,c)}),\pi_1(\F^{(\rho,\beta,c)}),c(\F^{(\rho,\beta,c)}))=(\rho,\beta,c)$.

Suppose now that $\F^{(\rho_1,\beta_1,c_1)}$ and $\F^{(\rho_2,\beta_2,c_2)}$ are monoidally isomorphic and let $\tau$ be a monoidal isomorphism, with components $\tau_g:\rho_1(g)\To\rho_2(g)$. Since $\G_s'$ is skeletal, this implies $\rho_1(g)=\rho_2(g)$ for all $g$ and hence, $\rho_1=\rho_2$. Furthermore, naturality in $g$ says that $\tau_g\circ F^{(\rho_1,\beta_1,c_1)}(\varphi)=F^{(\rho_2,\beta_2,c_2)}(\varphi)\circ\tau_g$ for all morphisms $\varphi:g\To g$ and, since $\pi_1(\G'_s)$ is abelian, the $\tau$'s cancel out each other. Hence, naturality implies no condition on the components of $\tau$ but requires, for such a $\tau$ to exist, that $F^{(\rho_1,\beta_1,c_1)}(\varphi)=F^{(\rho_2,\beta_2,c_2)}(\varphi)$ for all $\varphi$. In particular, this is true for all $u\in\pi_1(\G_s)$, which gives $\beta_1=\beta_2$. Finally, the monoidality conditions on $\tau$ (cf. Equations~(\ref{monoidalitat_transf_natural})) say that
\begin{align*}
\tau_e&={\id}_{e'} \\ \tau_{g_1g_2}\circ\gamma_{\rho(g_1g_2)}(c_1(g_1,g_2))&=\gamma_{\rho(g_1g_2)}(c_2(g_1,g_2))\circ(\tau_{g_1}\otimes\tau_{g_2})
\end{align*}
(we write $\rho_1=\rho_2=\rho$). If $\overline{\tau}:\pi_0(\G_s)\To\pi_1(\G'_s)$ is defined by $\overline{\tau}(g)=\gamma_{\rho(g)}^{-1}(\tau_g)$, the first condition says that $\overline{\tau}$ is a normalized 1-cochain, while the second says that $\partial\overline{\tau}=c_1-c_2$. Indeed, in terms of $\overline{\tau}$, this condition says that
$$
\gamma_{\rho(g_1g_2)}(\overline{\tau}(g_1g_2))\circ\gamma_{\rho(g_1g_2)}(c(g_1,g_2))=\gamma_{\rho(g_1g_2)}(c_2(g_1,g_2))\circ(\gamma_{\rho(g_1)}(\overline{\tau}(g_1))\otimes\gamma_{\rho(g_2)}(\overline{\tau}(g_2)))
$$
But in any 2-group $\G$, the $\gamma$'s satisfy $\gamma_A(u)\otimes\gamma_B(v)=\gamma_{A\otimes B}(u\circ ([A]\cdot v))$ for all objects $A,B$ and morphisms $u,v\in\pi_1(\G)$. Applying this to the right hand side and using that $\gamma_{\rho(g_1g_2)}$ is a group isomorphism, we obtain
$$
\overline{\tau}(g_1g_2)\circ c_1(g_1,g_2)=c_2(g_1,g_2)\circ \overline{\tau}(g_1)\circ [g_1\cdot(\overline{\tau}(g_2)]
$$
Hence, $c_1$ and $c_2$ are cohomologous. Conversely, if $c_1,c_2$ are cohomologous and $\overline{\tau}:\pi_0(\G_s)\To\pi_1(\G_s')$ is a normalized 1-cochain such that $c_1-c_2=\partial\overline{\tau}$, the reader may easily check that $\F^{(\rho,\beta,c_1)}$ and $\F^{(\rho,\beta,c_2)}$ are monoidally isomorphic, a monoidal isomorphism being defined by $\tau_g=\gamma_{\rho(g)}(\overline{\tau}(g))$.
\end{proof}

By applying this Lemma to $\G_s=\G_{\alpha}$ and $\G'_s=\G'_{\alpha'}$, we conclude there is a bijection between the set $Q(\G,\G';\alpha,\alpha')$ of triples $(\rho,\beta,[c])$ as in the statement of the theorem and the set $\pi_0({\bf Hom}^s_{{\sf 2Grp}}(\G_{\alpha},\G'_{\alpha'}))$. The desired correspondence is then the composite $Q(\G,\G';\alpha,\alpha')\To\pi_0({\bf Hom}^s_{{\sf 2Grp}}(\G_{\alpha},\G'_{\alpha'}))\To\pi_0({\bf Hom}_{{\sf 2Grp}}(\G_{\alpha},\G'_{\alpha'}))\To\pi_0(({\bf Hom}_{{\sf 2Grp}}(\G,\G'))$ and hence, it is given by
\begin{align}
(\rho,\beta,[c])&\longrightarrow[\EE'_{\alpha'}\circ\F^{(\rho,\beta,c)}\circ\EE_{\alpha}] \label{correspondencia1} \\ (\pi_0(\widetilde{\F}_{\alpha,\alpha'}),\pi_1(\widetilde{\F}_{\alpha,\alpha'}),[c(\widetilde{\F}_{\alpha,\alpha'})])&\longleftarrow [\F], \label{correspondencia2}
\end{align}
where $c$ is any normalized 2-cochain in the class $[c]$ and $\widetilde{\F}_{\alpha,\alpha'}:\G_{\alpha}\To\G'_{\alpha'}$ is the special morphism given by $\widetilde{\F}_{\alpha,\alpha'}=\overline{\hat{\EE}'_{\alpha'}\circ\F\circ\hat{\EE}_{\alpha}}$, for $\hat{\EE}_{\alpha},\hat{\EE}'_{\alpha'}$ some (quasi)inverses of $\EE_{\alpha}$ and $\EE'_{\alpha'}$, respectively.
\end{proof}
Note that, unlike the description of $\pi_0({\sf 2Grp})$, this description of $\pi_0({\bf Hom}_{{\sf 2Grp}}(\G,\G'))$ is non canonical, because it depends on the fixed classifying 3-cocycles $\alpha,\alpha'$ of $\G$ and $\G'$. However, this only concerns the cohomology class $[c]$, because $\pi_0(\F)$ and $\pi_1(\F)$ only depend on the isomorphism class of $\F$. Thus, if the isomorphism class of a morphism $\F:\G\To\G'$ is mapped to the triple $(\rho,\beta,[c])$ when we choose cocycles $\alpha,\alpha'$, by choosing different cocycles $\alpha_1,\alpha'_1$ it is mapped to a triple with the same $\rho$ and $\beta$ but a possibly different class $[c_1]$.

\subsection{Split 2-groups and examples}

A type of 2-groups, including those of the form $G[0]$ and $A[1]$, for which the representation theory simplifies a bit is the following:

\begin{defn}
A 2-group $\G$ is called {\sl split} if its classifying cohomology class is $[\alpha]=0$.
\end{defn}

Note that, if $\G$ is split, it is equivalent to the skeletal strict 2-group constructed by the method described in \S\ref{2-grups_modul_equivalencia} with $G=\pi_0(\G)$ and $M=\pi_1(\G)$. Conversely, if $\G$ is equivalent to a skeletal strict 2-group $\G'$, it easily follows from the above procedure of computing a classifying 3-cocycle (cf. \S\ref{2-grups_modul_equivalencia}) that a classifying 3-cocycle of $\G'$ is the trivial one and hence, that $\G$ is split. Therefore, we have the following alternative characterization of split 2-groups.

\begin{lem}
A 2-group $G$ is split if and only if it is equivalent to a skeletal strict 2-group.
\end{lem}

Since the classifying cohomology class of split 2-groups is trivial, these 2-groups can be identified with $G$-modules for some group $G$ or, equivalently, with the semidirect products of groups one of whose factors (that on which the other one acts) is abelian. Actually, there is a notion of action of a 2-group on another one and a corresponding notion of semidirect product for 2-groups (see \cite{GI01}) so that split 2-groups are (up to equivalence) those of the form $\G=G[0]\times_{\triangleleft} A[1]$ for some group $G$ and some abelian group $A$ ($\lhd$ denotes an action of $G$ on $A$; such an action indeed induces an action of $G[0]$ on $A[1]$). Its homotopy groups are $\pi_0(\G)=G$ and $\pi_1(\G)=A$.

\begin{exs} \label{exemples_2-grups_escindits} {\rm 
(1) Any linear representation
$\rho:G\To{\rm GL}(V)$ of an arbitrary group $G$ defines a split 2-group $\G$, with $\pi_0(\G)=G$ and $\pi_1(\G)=V$ (the underlying abelian group of the vector space $V$), $\pi_0(\G)$ acting on it according to $\rho$. These are called {\sl representational 2-groups} by Crane and Yetter \cite{CY03}. In particular, groups as important in geometry and physics as the euclidian 3-dimensional group $\EE(3)$ or the Poincar\' e group $\PP$ can be thought as split 2-groups, with $\pi_0(\EE(3))=SO(3)$, $\pi_1(\EE(3))=\R^3$ and the usual action of $SO(3)$ on $\R^3$, and $\pi_0(\PP)=SO(3,1)$, $\pi_1(\PP)=\R^4$ and the action of $SO(3,1)$ on $\R^4$ by rotations and boosts.

(2) The dihedral groups $D_{2m}$ ($m\geq 1$) and the symmetric $S_4$ and alternating $A_4$ groups are finite split 2-groups, with homotopy groups $\pi_0$ and $\pi_1$ respectively given by $\Integer_2$ and $\Integer_m$, $S_3$ and
$\Integer_2\times\Integer_2$ and $\Integer_3$ and
$\Integer_2\times\Integer_2$, and with the actions $\lhd$ defined by
$\overline{1}\lhd\overline{k}=-\overline{k}$ ($\overline{k}\in\Integer_m$) 
in the case of $D_{2m}$, by the obvious action of $S_3$ on the 3-element set
$\{(\overline{0},\overline{1}),(\overline{1},\overline{0}),(\overline{1},\overline{1})\}\subset\Integer_2\times\Integer_2$ 
in the case of $S_4$ and by the restriction of this action to
$\Integer_3\cong\{{\id},(123),(132)\}\subset S_3$ in the case of $A_4$. }
\end{exs}

\subsection{Strict 2-groups and crossed modules of groups}

\label{2-grups_estrictes}

Let us finish this section by recalling an alternative way of thinking of a 2-group when it is strict. Namely, as a {\it crossed module of groups}. Recall that this is a pair of groups $E$ and $N$ (none of them necessarily abelian), respectively called {\sl base} and {\sl principal} or {\sl fiber group}, together with 
an action $\lhd:E\times N\To N$ of $E$ on $N$ and a
group morphism $\partial:N\To E$, called the {\sl boundary}, such that
\begin{align}
\partial(e\lhd n)&=e\partial(n)e^{-1} \label{modul_creuat_1} \\
\partial(n)\lhd n'&=nn'n^{-1} \label{modul_creuat_2}
\end{align}
for all $e\in E$ and $n,n'\in N$ (in particular, ${\rm Im}\partial$ is a normal subgroup of $E$). For an introduction to crossed modules of groups, see
\cite{BRup}. A morphism between two crossed modules $(E,N,\partial,\lhd)$ and $(E',N',\partial',\lhd')$ is a commutative square
\begin{equation} \label{morfisme_modul_creuat_1}
\xymatrix{
N\ar[r]^{\partial}\ar[d]_{\psi} & E\ar[d]^{\phi} \\ N'\ar[r]_{\partial'} & E'
}
\end{equation}
in which the vertical arrows are compatible with the actions of $E$ and $E'$ on $N$ and $N'$, i.e., such that
\begin{equation} \label{morfisme_modul_creuat_2}
\phi(e\lhd n)=\phi(e)\lhd'\psi(n)
\end{equation}
for all $(e,n)\in E\times N$. This defines the category of crossed modules of groups, denoted ${\bf CMGrp}$. Then, we have the following well known result \footnote{The first published proof of this result seems to be that of Brown and Spencer \cite{BS76}, but the result was apparently known long time before.}:

\begin{thm} \label{equivalencia_Brown_Spencer}
Let ${\bf 2Grp}_s$ be the category with objects strict 2-groups and morphisms {\sl strict} morphisms of 2-groups. Then, there is an equivalence of categories ${\bf CMGrp}\simeq{\bf 2Grp}_s$.
\end{thm} 

On objects, the equivalence is defined as follows. Let $\G$ be an arbitrary strict 2-group and let $\G_0$ and $\G_1$ respectively denote its sets of objects and morphisms. It is easily seen that they are both groups with the tensor product. Then, $\G$ is mapped to the following crossed module $(E,N,\partial,\lhd)$. Let $E=\G_0$, let $N$ be the (normal) subgroup of $\G_1$ defined by the morphisms whose source is the unit object of $\G$, let $\partial:N\To E$ be the restriction of the target map $t:\G_1\To\G_0$ to $N$ and set
\begin{equation} \label{accio_G_sobre_H}
A\lhd (I_{\G}\stackrel{f}{\To} B)={\id}_A\otimes f\otimes {\id}_{A^{-1}}
\end{equation}
for all $A\in E$ and $f\in N$. Conversely, given a crossed module $(E,N,\partial,\lhd)$, it is mapped to the strict 2-group $\G$ with $\G_0=E$, $\G_1=E\times N$, a pair $(e,n)$ being thought of as a morphism $(e,n):e\To\partial(n)e$, and with the composition and tensor product defined by
\begin{equation} \label{comp_morf_2-grup_associat_modul_creuat}
(e',n')\circ(e,n)=(e,n'n)
\end{equation}
whenever $e'=\partial(n)e$, and
\begin{align} 
e\otimes e'&=ee' \label{prod_tens_obj_2-grup_associat_modul_creuat} \\
(e',n')\otimes(e,n)&=(e'e,n'(e'\lhd n)) \label{prod_tens_morf_2-grup_associat_modul_creuat}
\end{align}
(cf. Equations (\ref{composicio_morfismes_2-grup_associat})-(\ref{producte_tensorial_morfismes_2-grup_associat})). For more details, see \cite{FB02}.

\begin{exs} {\rm
(1) Any normal subgroup $N$ of a group $G$ defines a crossed module with $\partial:N\To G$ the inclusion and $\lhd$ the action of $G$ on $N$ by conjugation. Such crossed modules correspond to the strict 2-groups $\G$ with $\pi_1(\G)=1$. In particular, if $N=1$, one recovers the trivial 2-groups $G[0]$.

(2) For any group $G$, a crossed module is given by $E={\rm Aut}(G)$, $N=G$, $\partial:G\To{\rm Aut}(G)$ the morphism which maps any $g\in G$ to the corresponding inner automorphism and $\lhd$ the obvious one. The associated strict 2-group is ${\sf Aut}_{{\bf Cat}}(\underline{G})$, where $\underline{G}$ denotes the groupoid with only one object and $G$ as set of morphisms. }
\end{exs}

Strict 2-groups can also be classified starting with the corresponding crossed module. Thus, given a crossed module $(E,G,\partial,\lhd)$, the homotopy groups of the associated strict 2-group $\G$ are given by
\begin{align} 
\pi_0(\G)&=E/{\rm Im}\ \partial \label{pi0_via_modul_creuat} \\ 
\pi_1(\G)&={\rm Ker}\ \partial, \label{pi1_via_modul_creuat},
\end{align}
with the first acting on the second by
\begin{equation} \label{accio_pi0_pi1_via_modul_creuat}
[e]\cdot n=e\lhd n
\end{equation}
for all $n\in{\rm Ker}\ \partial$, $[e]\in E/{\rm Im}\ \partial$. A classifying 3-cocycle can be obtained as follows. Start with any set section $s:\pi_0(\G)\To E$ and let $\alpha_s:\pi_0(\G)\times\pi_0(G)\To {\rm Im}\ \partial$ be given by
$$
\alpha_s([e],[e'])=s([e])s([e'])s([ee'])^{-1}
$$
(it measures how far the section $s$ is from a group morphism). Associativity of the product in $E$ implies that $\alpha_s$ satisfies the non abelian cocycle condition
$$
(s([e])\cdot\alpha_s([e'],[e'']))\alpha_s([e]),[e'e''])=\alpha_s([e],[e'])\alpha_s([ee'],[e'']),
$$
with $E$ acting on ${\rm Im}\ \partial$ by conjugation. Then, choose any lift $\widetilde{\alpha}_s:\pi_0(\G)\times\pi_o(\G)\To N$ of $\alpha_s$. It turns out that $\widetilde{\alpha}_s$ no longer satisfies the previous cocycle condition, but the following weakened one:
$$
(s([e])\cdot\widetilde{\alpha}_s([e'],[e'']))\widetilde{\alpha}_s([e]),[e'e''])=\alpha([e],[e'],[e''])\widetilde{\alpha}_s([e],[e'])\widetilde{\alpha}_s([ee'],[e''])
$$
for some function $\alpha:\pi_0(\G)\times\pi_0(\G)\times\pi_0(\G)\To\pi_1(\G)$. It is shown that this function is a 3-cocycle and that its cohomology class is independent of the choices of $s$ and $\widetilde{\alpha}_s$.

Note that the strict 2-group $\G$ associated to a crossed module $(E,N,\partial,\lhd)$ is split if there exists a section $s$ as above which is a group morphism. Clearly, this is the case if the crossed module is {\sl trivial}, i.e., such that ${\rm Ker}\ \partial=N$. Strict 2-groups corresponding to such crossed modules are called {\sl automorphic 2-groups} by Crane and Yetter \cite{CY03}.



\section{2-category of representations of a 2-group}

\label{2cat_repr}

\subsection{Definition and unfolding}

Recall that a linear representation of a group $G$ can be defined either as a group morphism $\rho:G\To GL(V)$ for some vector space $V$ or equivalently, as a functor $F:\underline{G}\To {\bf Vect}$, where $\underline{G}$ denotes the groupoid with only one object and $G$ as set of morphisms and {\bf Vect} is the category of vector spaces. More generally, representations of $G$ in categories $\Cc$ other than {\bf Vect} also correspond to functors $F:\underline{G}\To\Cc$. Similarly, a morphism between two representations $\rho:G\To {\rm Aut}_{\Cc}(X)$ and $\rho':G\To {\rm Aut}_{\Cc}(X')$ can be defined either as a morphism $f:X\To X'$ in $\Cc$ which commutes with the actions of $G$ or as a natural transformation between the corresponding functors $F:\underline{G}\To\Cc$ and $F':\underline{G}\To\Cc$. Furthermore, in the categorical language, the composition of two such morphisms is given by the vertical composition of natural transformations. In other words, for any category $\Cc$, the category of representations of a group $G$ on $\Cc$ is nothing but the functor category $[\underline{G},\Cc]$. This together with the fact that a 2-group is a bigroupoid with only one object naturally suggests the following definition for the bicategory of representations of a 2-group on a bicategory $\Cgg$.

\begin{defn} \label{definicio_2-categoria_representacions}
For any bicategory $\Cgg$, the {\sl bicategory of representations} of a 2-group $\G$ on $\Cgg$, denoted $\repg_{\Cgg}(\G)$, is the pseudofunctor bicategory
$$
\repg_{\Cgg}(\G):=[\underline{\G},\Cgg]
$$
Objects in $\repg_{\Cgg}(\G)$ are called {\sl
representations} of $\G$ on $\Cgg$, while 1- and 2-morphisms are
respectively called {\sl 1-} and {\sl 2-intertwiners}. A representation $\Ff:\underline{\G}\To\Cgg$ (resp. a 1-intertwiner $\xi:\Ff\To\Ff'$) is called {\sl strict} when the pseudofunctor $\Ff$ is a 2-functor (resp. when the pseudonatural transformation $\xi$ is a 2-natural transformation).
\end{defn}

It immediately follows from Proposition~\ref{bicat_pseudof_equivalents}-(ii) that, for any bicategory $\Cgg$, the corresponding bicategory of representations $\repg_{\Cgg}(\G)$ is an invariant of $\G$ in the following sense:

\begin{lem}
For any bicategory $\Cgg$, if $\G$ and $\G'$ are equivalent 2-groups, the bicategories of representations $\repg_{\Cgg}(\G)$ and $\repg_{\Cgg}(\G')$ are biequivalent.
\end{lem}

Thus, when studying representations of a 2-group $\G$, and whenever convenient, one can assume without loss of generality that $\G$ is either strict or special (cf. Theorem~\ref{teorema_estrictificacio}).

\begin{rem} {\rm
In a previous version of this work, attention was focused on the {\it semistrict} representation theory, i.e., the full sub-bicategory of $\repg_{\Cgg}(\G)$ with objects only the 2-functors $\underline{\G}\To\Cgg$. However, B. Toen pointed out to me that this is not the right object to study, because it is not an invariant of $\G$ in the above sense. }
\end{rem}

Similarly, the replacement of the bicategory $\Cgg$ by another biequivalent bicategory $\Cgg'$ gives essentially the same representation theory, in the sense that $\repg_{\Cgg}(\G)$ and $\repg_{\Cgg'}(\G)$ are also biequivalent (it also follows from Proposition~\ref{bicat_pseudof_equivalents}-(ii)). Hence, the strictification theorem for bicategories (\cite{GPS95}, Thm. 1.4) further allows one to assume that $\Cgg$ is a 2-category, an assumption which is made from now on.  

The following more explicit descriptions of the notions of representation and 1- and 2-intertwiner readily follow from Definitions~\ref{def_pseudofunctor}, \ref{def_transf_pseudonatural} and \ref{def_modificacio} and the fact that $\underline{\G}$ has only one object. Observe that representations and 1-intertwiners have descriptions completely analogous to the familiar ones in the context of groups, except that a 1-intertwiner preserves now the actions of $\G$ only up to suitable (i.e., natural and coherent) 2-isomorphisms. 

\begin{prop} \label{representacio_cas_general}
Let $\G$ be a 2-group and $\Cgg$ an arbitrary 2-category. Then:
\begin{itemize}
\item[(i)]
A representation $\Ff:\underline{\G}\To\Cgg$ is a pair $(X,\F)$, where $X$ is an object of $\Cgg$ and $\F=(F,F_2):\G\To{\sf Equiv}_{\Cgg}(X)$ is a morphism of 2-groups.
\item[(ii)]
Given two representations
$\Ff,\Ff':\underline{\G}\To\Cgg$, with $\Ff=(X,\F)$ and
$\Ff'=(X',\F')$, a 1-intertwiner $\xi:\Ff\To\Ff'$ is a pair $(f,\Phi)$, where $f:X\To X'$ is a morphism in $\Cgg$ and $\Phi=\{\Phi(A):F'(A)\circ f\stackrel{\cong}{\Rightarrow} f\circ F(A)\}$ is a family of 2-isomorphisms also in $\Cgg$, indexed by the objects $A$ of $\G$, and such that 
\begin{itemize}
\item[$\bullet$]
$(1_{f}\circ F(\varphi))\cdot\Phi(A)=\Phi(B)\cdot(F'(\varphi)\circ 1_{f})$, for each morphism $\varphi:A\To B$ of $\G$ ({\sl naturality condition}).
\item[$\bullet$]
$\Phi(A\otimes B)\cdot(F'_2(A,B)\circ 1_f)=(1_f\circ F_2(A,B))\cdot(\Phi(A)\circ 1_{F(B)})\cdot(1_{F'(A)}\circ\Phi(B))$, for all objects $A,B$ of $\G$ ({\sl first coherence condition}).
\item[$\bullet$]
$F'_0\circ 1_f=(1_{f}\circ F_0)\cdot\Phi(I)$ ({\sl second coherence condition}).
\end{itemize}
\item[(iii)]
Given two 1-intertwiners $\xi,\zeta:\Ff\To\Ff'$, with $\xi=(f,\Phi)$ and
$\zeta=(g,\Psi)$, a 2-intertwiner $\ngg:\xi\Rightarrow\eta$ is a 2-morphism
$\tau:f\Rightarrow g$ in $\Cgg$ such that
\begin{equation} \label{axioma_2-intertwiner}
\Psi(A)\cdot(1_{F'(A)}\circ\tau)=(\tau\circ 1_{F(A)})\cdot\Phi(A)
\end{equation}
for every object $A$ of $\G$.
\end{itemize} 
\end{prop}

\subsection{Example}

\label{exemple_teoria_representacions} 

For any category $\Cc$, let $\Cc[0]$ be the corresponding {\sl locally discrete} 2-category (i.e., that with the same objects and morphisms as $\Cc$ and with only identity 2-arrows). Then, we have the following

\begin{prop} \label{exemple_mes_simple}
For any 2-group $\G$ and category $\Cc$, $\repg_{\Cc[0]}(\G)={\bf Rep}_{\Cc}(\pi_0(\G))[0]$,
where ${\bf Rep}_{\Cc}(\pi_0(\G))$ denotes the category of representations of $\pi_0(\G)$ on $\Cc$. More generally, if $\Cgg$ is a 2-category biequivalent (not necessarily equal) to the locally discrete 2-category $\Cc[0]$, there is a biequivalence $\repg_{\Cgg}(\G)\simeq{\bf Rep}_{\Cc}(\pi_0(\G))[0]$.
\end{prop}
\begin{proof}
For any object $X$ of $\Cc[0]$, it is easy to see that ${\sf Equiv}_{\Cc[0]}(X)={\rm Aut}_{\Cc}(X)[0]$. It follows (item $(i)$ in Proposition~\ref{representacio_cas_general}) that an object $\Ff$ of $\repg_{\Cc[0]}(\G)$ is a pair $(X,\rho)$ with $\rho:\pi_0(\G)\To{\rm Aut}_{\Cc}(X)$ a group morphism (note that isomorphic objects in $\G$ should be sent to isomorphic, hence equal, objects in ${\rm Aut}_{\Cc}(X)[0]$). Moreover, according to item $(ii)$, a 1-intertwiner $\xi:(X,\rho)\To(X',\rho')$ consists of a morphism $f:X\To X'$ in $\Cc$ together with 2-isomorphisms $\Phi(A):f\circ\rho([A])\cong\rho'([A])\circ f$, which are necessarily identities because $\Cc[0]$ is discrete. By the same reason, there are no non trivial 2-intertwiners, and Equation~(\ref{composicio_vertical_transformacions_pseudo_1}) implies that composition of 1-intertwiners really corresponds to composing the $f$'s.
\end{proof}

For instance, let ${\sf 2Grp}_1$ be the full sub-2-category of {\sf 2Grp} with objects only those of type $A[1]$ for some abelian group $A$. Then, it is easy to check that ${\sf 2Grp}_1$ is biequivalent to its sub-2-category ${\sf 2Grp}_1^s$ with only the strict monoidal functors as 1-arrows, which is in turn equal to ${\sf 2Grp}_1^s={\bf Ab}[0]$. Hence, it is $\repg_{{\sf 2Grp}_1}(\G)\simeq{\bf Rep}_{\bf Ab}(\pi_0(\G))[0]$.

Clearly, the cases considered in the Proposition~\ref{exemple_mes_simple} constitute the coarsest possible representation theories for 2-groups. In these cases, all 2-groups ${\rm Equiv}_{\Cgg}(X)$ are discrete (up to equivalence). Finer theories are obtained by considering target 2-categories $\Cgg$ such that the 2-groups ${\rm Equiv}_{\Cgg}(X)$ are split, but non discrete. As discussed below, this is the case when $\Cgg={\bf 2Mat}_{\Complex}$.

\subsection{Describing the objects of $\repg_{\Cgg}(\G)$ up to isomorphism}

The proof of the next result is an easy exercise left to the reader.

\begin{lem} \label{lema_representacions_isomorfes}
Let $X$ be any object of a 2-category $\Cgg$ and let $\F,\F':\G\To{\sf Equiv}_{\Cgg}(X)$ be two representations of a 2-group $\G$ as autoequivalences of $X$. Then, if $\F$ and $\F'$ are isomorphic 1-morphisms in {\sf 2Grp}, the pairs $(X,\F)$ and $(X,\F')$ are isomorphic (hence equivalent) objects in $\repg_{\Cgg}(\G)$.
\end{lem}

Observe that the converse statement is false. There may exists pairs $(X,\F),(X,\F')$ with $\F$ and $\F'$ non monoidally isomorphic which are however isomorphic (or at least equivalent) as representations of $\G$. This is because, even when $X=X'$ in $(ii)$, a generic 1-morphism between $(X,\F)$ and $(X,\F')$ is not always a monoidal natural transformation between the monoidal functors $\F$ and $\F'$, due to the possibly non trivial morphism $f:X\To X$. This becomes clearer below (cf \S~\ref{classes_equivalencia_repr}).

Since one is ultimately interested in the representations up to equivalence, a first useful consequence of the previous Lemma is that it allows us to restrict to representations $(X,\F)$ for which the isomorphism $F_0$ is equal to $1_{{\id}_X}$ (cf. Lemma~\ref{lema_functors_monoidals}). For 1-intertwiners $(f,\Phi)$ between two such representations, the second coherence condition on the family $\Phi$ simplifies to $\Phi(I)=1_f$.

The previous Lemma further allows one to give a more concrete description, up to isomorphism, of the objects of $\repg_{\Cgg}(\G)$ when $\Cgg$ is any target 2-category.

\begin{prop} \label{repr_modul_isomorfisme_cas_general}
Up to isomorphism, a representation of a 2-group $\G$ on a 2-category $\Cgg$ is given by a quadruple $(X,\rho,\beta,c)$, with $X$ an object of $\Cgg$, $\rho:\pi_0(\G)\To\pi_0({\sf Equiv}_{\Cgg}(X))$ a group morphism, $\beta:\pi_1(\G)\To\pi_1({\rm Equiv}_{\Cgg}(X))_{\rho}$ a morphism of $\pi_0(\G)$-modules such that $[\beta\circ\alpha]=[\alpha_X\circ\rho^3]$ in $H^3(\pi_0(\G),\pi_1({\rm Equiv}_{\Cgg}(X))_{\rho})$ and $c\in C^2(\pi_0(\G),\pi_1({\rm Equiv}_{\Cgg}(X))_{\rho})$ a normalized 2-cochain such that $\partial c=\beta\circ\alpha-\alpha_X\circ\rho^3$, where $\alpha$, $\alpha_X$ are fixed classifying 3-cocycles of $\G$ and ${\rm Equiv}_{\Cgg}(X)$, respectively. Furthermore, any such quadruple determines a representation of $\G$ on $\Cgg$ (up to isomorphism).
\end{prop}
\begin{proof}
The result is an immediate consequence of Theorem~\ref{morfismes_modul_isomorfisme}, Proposition~\ref{representacio_cas_general}-$(i)$ and Lemma~\ref{lema_representacions_isomorfes}.
\end{proof}
Note that the isomorphism class of the representation determined by $(X,\rho,\beta,c)$ uniquely depends on the cohomology class of $c$, so that this gives a non faithful parametrization of the set of isomorphism classes of objects in $\repg_{\Cgg}(\G)$. Let us also remark that the really interesting thing is not this set of isomorphism classes, which is not an invariant (under biequivalence) of $\repg_{\Cgg}(\G)$, but the set of equivalences classes. Later on, this set is described explicitly when $\Cgg={\bf 2Mat}_{\Complex}$ (in this case, it turns out to coincide with the set of isomorphism classes). Finally, observe that when $\Cgg=\Cc[0]$, the quadruple $(X,\rho,\beta,c)$ reduces to a representation of $\pi_0(\G)$ on the group of automorphisms of some object $X$ of $\Cc$, in accordance with Example~\ref{exemple_teoria_representacions}.



\section{2-vector spaces and 2-matrices}

\label{2cat_2vect}

As in the classical setting of groups, the goodness of a representation theory for 2-groups is measured by the amount of information on the 2-group hidden in the corresponding 2-category of representations. The ideal situation occurs when the 2-group can be reconstructed (up to equivalence) from the corresponding 2-category of representatins. This is obviously not the case for all choices of $\Cgg$. Thus, as discussed in the previous section, representing 2-groups on locally discrete 2-categories implies loosing all information on the 2-group except, at most, its $\pi_0$.

In this paper, we focuse the attention on an analog of the classical {\sl linear} representation theory of groups, i.e., we consider representations on a certain categorical analog of a vector space. Such structures are usually called {\sl 2-vector spaces}.

\subsection{Notions of 2-vector space}

Let $K$ be a fixed field. To our knowledge, two notions of 2-vector space over $K$ have been proposed until now:

\begin{itemize}
\item Kapranov and Voevosdky's 2-vector spaces, whose definition was introduced in \cite{KV94} to give a conceptual framework for the Zamolodchikov tetrahedra equations, and
\item Baez and Crans's 2-vector spaces, a much more recent definition, introduced in \cite{BC03} in order to define a categorified notion of a Lie algebra.
\end{itemize}
Both definitions are intended to be ``categorifications'' of the usual notion of vector space over $K$. The simplest one is that by Baez and Crans, which reads as follows.

\begin{defn}
A {\sl Baez and Crans 2-vector space} over $K$ is a category in the category ${\bf Vect}_K$ of vector spaces over $K$.
\end{defn}

In other words, it is a category $\V$ whose sets of objects and morphisms, respectively denoted $V_0$ and $V_1$, are vector spaces over $K$ and such that the source, target, identity assigning and composition maps are $K$-linear. By appropriately internalising the notions of functor and natural transformation, Baez and Crans organize these 2-vector spaces into a 2-category $\dv^{BC}_K$ and show that $\dv^{BC}_K$ is is biequivalent to the 2-category ${\bf Comp}_1(K)$ with objects the length 1 complexes of vector spaces over $K$ (i.e., linear maps $d:V_1\To V_0$), and with morphisms of complexes as 1-arrows and chain homotopies as 2-arrows (cf. \cite{BC03}).

In this paper, however, we wish to focus attention on the representation theory on Kapranov and Voevodsky's 2-vector spaces, defined as follows.

\begin{defn}
A {\sl Kapranov and Voevodsky (finite dimensional) 2-vector space} over $K$ is a ${\bf Vect}_K$-module category which is ${\bf Vect}_K$-module equivalent to ${\bf Vect}_K^n$ for some $n\geq 0$, called the {\sl rank}.   
\end{defn}

Unlike the previous definition, completely unpacking this definition will take some pages. Hence, let us just say that a ${\bf Vect}_K$-module category is a symmetric monoidal category $\Cc$ (analog of the underlying abelian group of a vector space) equipped with a (left) action of ${\bf Vect}_K$ (analog of the action of $K$ on the abelian group), i.e., a bifunctor $\odot:{\bf Vect}_K\times\V\To\V$ together with natural isomorphisms arising from the weakening of the usual axioms of a module over a ring (for ex., isomorphisms of the type $u_{V,V',X}:V\odot(V'\odot X)\cong(V\otimes V')\odot X$, $v_{V,V',X}:(V\oplus V')\odot X\cong(V\odot X)\overline{\oplus}(V'\odot X)$, etc., where $\overline{\oplus}$ denotes the symmetric tensor product on $\Cc$). For more details, the reader is refered to Kapranov and Voevodsky's paper \cite{KV94}.

Actually, the previous two notions of 2-vector space are particular instances of the more general notion of {\sl 2-(semi)modules over a 2-semiring}. Roughly, a {\sl 2-semiring}, as a 2-group, is the appropriately weakened version of the notion of a semiring object in {\bf Cat} (i.e., a category equipped with two functors, called sum and product, satisfying the usual axioms of a semiring). For example, the categories $K[0]$ and ${\bf Vect}_K$ are 2-semirings (the first one is even a 2-ring), with the sum and product functors given by the sum and product of $K$ and the direct sum and tensor product of vector spaces, respectively. Given a 2-semiring $\Ss$, a {\sl 2-(semi)module} over $\Ss$ is then a symmetric 2-group or more generally, a symmetric monoidal category (analogous to the abelian group in the classical context) equipped with a (left) action of $\Ss$ in the same sense as before. In this language, a Baez and Crans 2-vector space over $K$ turns out to be equivalent to a 2-module over the 2-ring $K[0]$, while Kapranov and Voevodsky 2-vector spaces are 2-(semi)modules over the 2-semiring ${\bf Vect}_K$. Thus, the main difference between both notions, and actually the main point in the categorification of the notion of vector space over $K$, is what one takes as analog of the field $K$, either $K[0]$ or ${\bf Vect}_K$, or anything else intended to be some sort of categorification of $K$. In this sense, it is worth noting that Kapranov and Voevosky's choice ${\bf Vect}_K$ is in fact a categorification of the set $\Natural$ of natural numbers, not of $K$.

As Baez and Crans, Kapranov and Voevodsky introduce suitable notions of 1- and 2-morphism between their 2-vector spaces and organize them into a 2-category $\dv_K^{KV}$. To avoid introducing these notions, one may take as model of this 2-category the following simplified version (see \cite{dYnp}). 

\begin{thm}
$\dv_K^{KV}$ is biequivalent to the 2-category $\dv_K$ whose objects, 1- and 2-morphisms are the $K$-linear categories ${\bf Vect}_K^n$ for $n\geq 0$, the $K$-linear functors and the natural transformations, respectively, with the usual composition laws and identity arrows.
\end{thm}

\subsection{The 2-category of 2-matrices}

In developping the representation theory, we shall take neither $\dv_K^{KV}$ nor $\dv_K$ as model for the 2-category of Kapranov and Voevodsky's 2-vector spaces. Instead, and in view of possible future applications of the theory involving the monoidal structure on these 2-categories (cf. also Remark~\ref{remarca_fidelitat_representacions}), we shall take as model the {\sl 2-category of 2-matrices} ${\bf 2Mat}_K$ introduced in \cite{jE3}. This 2-category is to be thought of as the analog in this setting of ${\bf Mat}_K$, the category with objects the natural numbers and morphisms from $n$ to $m$ the $m\times n$ matrices with entries in $K$, which appears as a coordinatized version of ${\bf Vect}_K$. The interesting feature of ${\bf Mat}_K$ is that it provides a strict monoidal version of the weak monoidal category ${\bf Vect}_K$. Similarly, ${\bf 2Mat}_K$ has the advantage over the versions ${\bf 2Vect}_K^{KV}$ and $\dv_K$ of being a semistrict monoidal 2-category, instead of a weak one. It is worth emphasizing, however, that the choice of model is somewhat irrelevant because, as pointed out before, the corresponding representation theories are the same in all three cases (up to biequivalence, of course).

Here, we just recall the definition of ${\bf 2Mat}_K$, and refer the reader to \cite{jE3} for more details. 

\begin{defn}
The {\sl 2-category of 2-matrices} over $K$, denoted ${\bf 2Mat}_K$, is the 2-category defined by the following data:
\begin{itemize}
\item
objects: the natural numbers $n\geq 0$.
\item
1-morphisms: if $n,m\geq 1$, a 1-morphism $n\To m$ is a pair
$({\bf R},s)$, where
\begin{itemize}
\item[(i)]
${\bf R}=(R_{ij})$ is an $m\times n$ matrix with entries
in $\Natural$, and
\item[(ii)]
$s$ is a collection of $m$ fields $s_1,\ldots,s_m$ of invertible matrices over $\Natural^n$, with $s_i({\bf a})\in {\rm GL}({\bf R}({\bf
  a})_i)$ for all ${\bf a}\in\Natural^n$ such that ${\bf
  R}({\bf a})_i\neq 0$ and 
$s_i({\bf a})=1$ otherwise, and satisfying the {\sl normalization
condition}
\begin{equation} \label{condicio_normalitzacio_gauge}
s_i({\bf e}_j)={\bf I}_{R_{ij}},\quad i=1,\ldots,m,\ j=1,\ldots,n
\end{equation}
whenever ${\bf R}({\bf e}_j)_i\neq 0$ (otherwise, it is equal to 1).
\end{itemize}
In particular, for any $n,m\geq 1$, there is a {\sl zero 1-morphism}, denoted ${\bf 0}_{n,m}$ and defined by the $m\times n$ zero matrix and the unique possible gauge. {\bf R} is called the {\sl rank matrix} and $s$
the {\sl gauge} of the 1-morphism. When $n=0$ and/or $m=0$, there
is only one 1-morphism $n\To m$, also denoted by ${\bf 0}_{n,m}$.
\item
2-morphisms: given objects $n,m\geq 1$ and 1-morphisms $({\bf
  R},s),({\bf R}',s'):n\To m$, with ${\bf R}=(R_{ij})$ and ${\bf
  R}'=(R'_{ij})$, a 2-morphism $({\bf R},s)\Rightarrow({\bf R}',s')$ is an
$m\times n$ matrix {\sf T} whose $(i,j)$ entry ${\sf T}_{ij}$ is an
$R'_{ij}\times R_{ij}$ complex matrix if both
$R_{ij},R'_{ij}\neq 0$ and empty otherwise. If $n=0$ and/or $m=0$,
both 1-morphisms are necessarily equal and there is a unique
2-morphism between them denoted $1_{{\bf 0}_{n,m}}$ (obviously, it
is the identity 2-morphism of ${\bf 0}_{n,m}$).
\item
composition of 1-morphisms: if $({\bf R},s):n\To m$ and
$(\widetilde{{\bf
  R}},\widetilde{s}):m\To p$, with $n,m,p\geq 1$, their composite
$(\widetilde{{\bf R}},\widetilde{s})\circ({\bf R},s):n\To p$ is the
pair
\begin{equation} \label{composicio_1-morfismes}
(\widetilde{{\bf R}},\widetilde{s})\circ({\bf
R},s)=(\widetilde{{\bf R}}{\bf R},\widetilde{s}\ast s)
\end{equation}
where $\widetilde{{\bf R}}{\bf R}$ denotes the usual matrix product
and $\widetilde{s}\ast s$ is defined by
\begin{align} \label{seccions_functor_composicio}
(\widetilde{s}\ast s)_k({\bf a})=\widetilde{s}_k({\bf R}({\bf a}))&
\left(\bigoplus_{i=1}^m {\bf I}_{\widetilde{R}_{ki}}\otimes
s_i({\bf a})\right){\bf P}(\widetilde{{\bf R}}_k,{\bf R},{\bf
  a}) \nonumber \\ &\hspace{1.3
truecm}\left(\bigoplus_{j=1}^n\ \widetilde{s}_k({\bf R}({\bf
e}_j))^{-1}\otimes {\bf I}_{a_j}\right)
\end{align}
for all ${\bf a}\in\Natural^n$ and $k=1,\ldots,p$. Here, ${\bf I}_0$
denotes the empty matrix (in particular, any tensor product ${\bf I}_0\otimes
s_i({\bf a})$ is again the empty matrix and it makes no contribution to
the previous direct sum) and ${\bf P}(\widetilde{{\bf R}}_k,{\bf R},{\bf a})$ a suitable
permutation matrix of order
$\sum_{j=1}^n\sum_{i=1}^m\widetilde{R}_{ki}R_{ij}a_j$. In case one of the
numbers $n,m,p$ is zero, the composite is the corresponding zero
map.
\item
vertical composition of 2-morphisms: given 2-morphisms ${\sf
T}:({\bf
  R},s)\Rightarrow({\bf R}',s')$ and
${\sf T}':({\bf R}',s')\Rightarrow({\bf R}'',s'')$, with $({\bf
R},s),({\bf R}',s'),({\bf R}'',s''):n\To m$ and $n,m\geq 1$, the
vertical composite ${\sf T}'\cdot{\sf T}$ is the matrix obtained by
multiplying componentwise both matrices {\sf T} and ${\sf T}'$,
i.e.
\begin{equation} \label{composicio_vertical}
({\sf T}'\cdot{\sf T})_{ij}={\sf T}'_{ij}{\sf T}_{ij}, \quad
    i=1,\ldots,m,\ j=1,\ldots,n
\end{equation}
where we agree that the product of a matrix by the empty matrix
is the empty matrix or the appropriate zero matrix when convenient (for
example, 
it can be $R'_{ij}=0$, so that ${\sf T}_{ij}={\sf T}'_{ij}=\emptyset$,
but $R_{ij},R''_{ij}\neq 0$; the matrix $({\sf T}'\cdot{\sf T})_{ij}$
is then the $R''_{ij}\times R_{ij}$ zero matrix rather than the empty one).

\item
horizontal composition of 2-morphisms: given 2-morphisms ${\sf
T}:({\bf
  R},s)\Rightarrow({\bf R}',s'):n\To m$ and $\widetilde{{\sf T}}:(\widetilde{{\bf R}},\widetilde{s})\Rightarrow(\widetilde{{\bf R}}',\widetilde{s}'):m\To
p$, with $n,m,p\geq 1$, the horizontal composite $\widetilde{{\sf
T}}\circ{\sf T}$ is the $p\times n$ matrix
with entries defined by
\begin{equation} \label{composicio_horitzontal}
(\widetilde{{\sf T}}\circ{\sf T})_{kj}=\widetilde{s}'_k({\bf
R}'({\bf e}_j))\left(\bigoplus_{i=1}^m\ \widetilde{{\sf
T}}_{ki}\otimes{\sf T}_{ij}\right)\widetilde{s}_k({\bf R}({\bf
e}_j))^{-1}
\end{equation}
for all $j=1,\ldots,n$ and $k=1,\ldots,p$, where we agree again
that the tensor product of any matrix by the empty matrix is the
empty matrix.
\item
identity 1-morphisms: for $n\geq 1$, ${\id}_n$ is the pair $({\bf I}_n,{\bf I})$, where {\bf
I} denotes the trivial gauge, namely, the gauge $s$ such that, for
any ${\bf a}\in\Natural^n$ and $i=1,\ldots,n$, $s_i({\bf a})\in{\rm
GL}(a_i)$ is the identity matrix if $a_i\neq 0$ and equal to 1
otherwise.
\item
identity 2-morphisms: for any 1-morphism $({\bf R},s):n\To m$, ${\bf
1}_{({\bf R},s)}$ is the $m\times n$ matrix given by $({\bf 1}_{({\bf R},s)})_{ij}={\bf
I}_{R_{ij}}$ if $R_{ij}\neq 0$ and empty otherwise.
\end{itemize}
\end{defn}

In what follows, all it is needed to know about the permutation matrices appearing in (\ref{seccions_functor_composicio}) are the the following {\sl normalization conditions} (cf. \cite{jE3}):

\begin{itemize}
\item[(i)] ${\bf P}({\bf e}_i,{\bf R},{\bf a})={\bf
  I}_{\sum_{j=1}^nR_{ij}a_j}$, for all $i=1,\ldots,m$;
\item[(ii)]
${\bf P}(\widetilde{{\bf R}}_k,{\bf R},{\bf e}_j)=
{\bf I}_{\sum_{i=1}^m\widetilde{R}_{ki}R_{ij}}$, for all
$j=1,\ldots,n$;
\item[(iii)]
${\bf P}(\widetilde{{\bf R}}_k,{\bf I}_n,{\bf a})={\bf
I}_{\sum_{j=1}^n\widetilde{R}_{kj}a_j}$, and
\item[(iv)]
${\bf P}(\widetilde{{\bf R}}_k,{\bf r},a)= {\bf
I}_{a\sum_{i=1}^m\widetilde{R}_{ki}r_i}$ when {\bf R} is an $m\times
  1$ matrix {\bf r}, so that {\bf a} reduces to a number.
\end{itemize}
In particular, (ii) ensures that $\widetilde{s}\ast s$ indeed satisfies the required normalization condition (\ref{condicio_normalitzacio_gauge}).

\begin{rem} \label{remarca_2-morfismes_identitat} {\rm
Note that a 2-morphism {\sf T} may have all its nonempty entries
equal to identity matrices even in case it is not an identity 2-morphism. Just
think of 2-morphisms {\sf T} between two 1-morphisms having the same rank matrix but
different gauges. }
\end{rem}

The following properties concerning 1- and 2-morphisms in ${\bf 2Mat}_K$
readily follow from the previous definition and its proof is left to the reader:

\begin{lem} \label{morfismes_invertibles}
Let $({\bf R},s),({\bf R}',s'):n\To m$ be 1-morphisms in ${\bf
  2Mat}_K$. Then:
\begin{itemize}
\item[(i)]
$({\bf R},s)$ is invertible if and only if $n=m$ and its rank matrix {\bf R} is a permutation
matrix.
\item[(ii)]
$({\bf R},s)$ and $({\bf R}',s')$ are 2-isomorphic if and only if
${\bf R}={\bf R}'$, and a 2-morphism ${\sf T}:({\bf
R},s)\Rightarrow({\bf R},s')$ is a 2-isomorphism if and only if all
nonempty entries in {\sf T} are non singular complex matrices.
\item[(iii)]
$({\bf R},s)$ is an equivalence if and only if it is a 1-isomorphism. 
\end{itemize}
\end{lem}

In what follows, we shall take $K=\Complex$, the field of complex numbers.



\section{Representation theory on ${\bf 2Mat}_{\Complex}$}

\label{seccio_repr_2-e.v}

From now on, when speaking of a representation of a 2-group $\G$, we shall always mean a representation on the 2-category ${\bf 2Mat}_{\Complex}$. To emphasize this, we shall use the name {\sl (complex) 2-matrix representation}. Moreover, since any 2-group is equivalent to a special 2-group and representation theories of equivalent 2-groups are equivalent, {\it in all what follows we make the simplifying assumption that $\G$ is special}. We recall this means that $\G$ is skeletal (hence, the set of objects can be identified with $\pi_0(\G)$), that all objects of $\G$ are strictly invertible and that the unit object $I$ is a strict unit for the tensor product.

\subsection{Kapranov and Voevodsky's general linear 2-groups ${\sf Equiv}_{{\bf 2Mat}_{\Complex}}(n)$}

According to Proposition~\ref{repr_modul_isomorfisme_cas_general}, a representation of a 2-group $\G$ on ${\bf 2Mat}_{\Complex}$ is determined, up to isomorphism, by a natural number $n\geq 0$ together with some more data involving the groups $\pi_0({\sf Equiv}_{{\bf 2Mat}_{\Complex}}(n))$ and $\pi_1({\sf Equiv}_{{\bf 2Mat}_{\Complex}}(n))$, the action of the first onto the second and a classifying 3-cocycle of ${\sf Equiv}_{{\bf 2Mat}_{\Complex}}(n)$. Hence, the first step in studying the representation theory on ${\bf 2Mat}_{\Complex}$ is to classify the 2-groups ${\sf Equiv}_{{\bf 2Mat}_{\Complex}}(n)$ for all $n\geq 0$, which play in this setting the role usual general linear groups ${\rm GL}(n)$ play in matrix representation theory of groups.

Clearly, ${\sf Equiv}_{{\bf 2Mat}_{\Complex}}(0)\simeq {\bf 1}$. Hence, let us consider the case $n\geq 1$. By item $(iii)$ in Lemma~\ref{morfismes_invertibles}, it is
${\sf Equiv}_{{\bf 2Mat}_{\Complex}}(n)={\sf Aut}_{{\bf 2Mat}_{\Complex}}(n)$.
Consequently, ${\sf Equiv}_{{\bf 2Mat}_{\Complex}}(n)$ is a strict 2-group and it can be identified with a crossed module $(E_n,N_n,\partial_n,\lhd_n)$ via Theorem~\ref{equivalencia_Brown_Spencer}. This crossed module can be computed explicitly and hence, used to classify ${\sf Equiv}_{{\bf 2Mat}_{\Complex}}(n)$ by the method described in \S~\ref{2-grups_estrictes}. Althought the classification can be done directly, without using this crossed module (cf. \S~\ref{2-grups_modul_equivalencia}), we compute it explicitly for later reference (cf. Remark~\ref{remarca_fidelitat_representacions}).

Let us first introduce some notation. Given any permutation $\sigma\in S_n$, let ${\bf P}(\sigma)$ be the $n\times n$ permutation matrix defined by
\begin{equation} \label{matriu_permutacio_sigma}
{\bf P}(\sigma)_{ij}=\delta_{i,\sigma(j)},\quad i,j=1,\ldots,n
\end{equation}
Given two permutations $\sigma,\sigma'\in S_n$, we shall denote by $\sigma'\sigma$ the permutation defined by $(\sigma'\sigma)(i)=\sigma'(\sigma(i))$ for all
$i\in\{1,\ldots,n\}$. In this way, we have ${\bf P}(\sigma'\sigma)={\bf P}(\sigma'){\bf
  P}(\sigma)$. Finally, let us define
$$
{\bf GL}_n:=\prod_{{\bf a}\in(\Natural^n)^*}\ {\rm
GL}(a_1)\times\cdots\times{\rm GL}(a_n)\quad (n\geq 1)
$$
where $(\Natural^n)^*=\Natural^n\backslash\{{\bf e}_1,\ldots,{\bf
e}_n\}$ (when $a_i=0$, we agree that ${\rm GL}(0)=1$).

Let us think of ${\bf GL}_n$ as a group with the direct product group structure. Then, there is a natural left action $\rhd$ of $S_n$ on ${\bf GL}_n$ by
group automorphisms defined as follows. If $\sigma\in S_n$ and 
$A=((A_1({\bf a}),\ldots,A_n({\bf a})))_{{\bf
a}\in(\Natural^n)^*}\in{\bf GL}_n$, with $A_i({\bf a})\in{\rm GL}(a_i)$, let
\begin{equation} \label{accio_S_n_GL_n}
(\sigma\rhd A)_i({\bf a})=A_{\sigma^{-1}(i)}(\sigma^{-1}\cdot{\bf
a}),\quad {\bf a}\in(\Natural^n)^*,\ i=1,\ldots,n
\end{equation}
where $\cdot$ denotes the usual left action of $S_n$ on
$\Complex^n$, i.e., for any ${\bf a}=(a_1,\ldots,a_n)$, 
\begin{equation} \label{wreath_product}
\sigma\cdot{\bf a}:=P(\sigma)({\bf
  a})=(a_{\sigma^{-1}(1)},\ldots,a_{\sigma^{-1}(n)}) 
\end{equation}
Notice that the
$\sigma^{-1}(i)$-entry of $\sigma^{-1}\cdot{\bf a}$ is
$a_{\sigma(\sigma^{-1}(i))}=a_i$, so that
$A_{\sigma^{-1}(i)}(\sigma^{-1}\cdot{\bf a})\in{\rm GL}(a_i)$ as required. We leave to the reader to check that (\ref{accio_S_n_GL_n}) indeed defines a left action by
automorphisms of ${\bf GL}_n$. We denote by $S_n\times_{\triangleright} {\bf GL}_n$ the corresponding semidirect product.

\begin{prop} \label{modul_creuat_n}
There is an isomorphism of crossed modules $(E_n,N_n,\partial_n,\lhd_n)\cong(S_n\times_{\triangleright}{\bf GL}_n,(\Complex^*)^n\times{\bf
GL}_n,\partial_n,\lhd_n)$, with $\partial_n$ and $\lhd_n$ defined by
\begin{align} \label{morfisme_parcial}
\partial_n(\mbox{\boldmath$\lambda$},A)&=({\id}_n,A)
\\ \label{accio}
(\sigma,A)\lhd_n(\mbox{\boldmath$\lambda$},A')&=
(\sigma\cdot\mbox{\boldmath$\lambda$},A(\sigma\rhd A')A^{-1})
\end{align}
for all $\sigma\in S_n$, $A,A'\in{\bf GL}_n$ and
$\mbox{\boldmath$\lambda$}\in(\Complex^*)^n$ ($\cdot$ and $\rhd$ respectively denote the actions {\rm (\ref{wreath_product})} and {\rm (\ref{accio_S_n_GL_n})}).
\end{prop}

\begin{proof}
Let $\Phi:S_n\times_{\triangleright}{\bf GL}_n\To E_n$ and $\Psi:E_n\To
S_n\times_{\triangleright}{\bf GL}_n$ be the maps defined by
\begin{align} \label{definicio_Phi}
\Phi(\sigma,A)&=({\bf P}(\sigma),s^{(\sigma,A)})
\\ \label{definicio_Psi}
\Psi({\bf P},s)&=(\sigma({\bf P}),A^{({\bf P},s)})
\end{align}
where $s^{(\sigma,A)}$ denotes the (normalized) gauge such that
\begin{equation} \label{s_A}
s_i^{(\sigma,A)}({\bf a})=A_i(\sigma\cdot{\bf a}),\quad {\bf
a}\in(\Natural^n)^*,\ i=1,\ldots,n,
\end{equation}
$\sigma({\bf P})$ is the unique permutation $\sigma\in S_n$ such
that ${\bf P}={\bf P}(\sigma)$ and $A^{({\bf P},s)}$ is given by
\begin{equation} \label{A_s}
A_i^{({\bf P},s)}({\bf a})=s_i({\bf P}^{-1}({\bf a})),\quad {\bf
a}\in(\Natural^n)^*,\ i=1,\ldots,n.
\end{equation}
It is immediate to check that $\Psi\circ\Phi$ and
$\Phi\circ\Psi$ are identities. Furthermore, for any $(\sigma,A),(\sigma',A')\in S_n\times_{\triangleright}{\bf GL}_n$, we have
\begin{align*}
\Phi((\sigma',A')(\sigma,A))&=({\bf P}(\sigma'\sigma),s^{(\sigma'\sigma,A'(\sigma'\rhd A))}) \\ 
\Phi(\sigma',A')\circ\Phi(\sigma,A)&=({\bf P}(\sigma'\sigma),s^{(\sigma',A')}\ast s^{(\sigma,A)})
\end{align*}
and it is easily seen that both gauges $s^{(\sigma'\sigma,A'(\sigma'\rhd A))}$ and $s^{(\sigma',A')}\ast s^{(\sigma,A)}$ coincide and are given by (cf. (\ref{seccions_functor_composicio}))
$$
s^{(\sigma'\sigma,A'(\sigma'\rhd A))}_i({\bf a})=
(s^{(\sigma',A')}\ast s^{(\sigma,A)})_i({\bf
  a})=A'_i((\sigma'\sigma')\cdot{\bf a})\ A_{(\sigma')^{-1}(i)}(\sigma\cdot{\bf
  a}) 
$$
Hence, $\Phi$ is an isomorphism of groups. In particular, notice that, if $(\sigma,A)\in
S_n\times_{\triangleright}{\bf GL}_n$ is the pair describing the
1-automorphism $({\bf P},s):n\To n$ in ${\bf 2Mat}_{\Complex}$, the
pair describing the inverse 1-automorphism $({\bf P},s)^{-1}$ is
\begin{equation} \label{automorfisme_invers_1}
\Psi(({\bf P},s)^{-1})=(\sigma,A)^{-1}=(\sigma^{-1},\sigma^{-1}\rhd A^{-1})
\end{equation}
Hence, we may write
\begin{equation} \label{automorfisme_invers_2}
({\bf P},s)^{-1}=({\bf P}^{-1},s^{-1})
\end{equation}
if we define $s^{-1}$ as the normalized gauge given by
$$
s^{-1}_i({\bf a})=s_i^{(\sigma^{-1},\ \sigma^{-1}\rhd A^{-1})}({\bf
a})=(\sigma^{-1}\rhd A^{-1})_i(\sigma^{-1}\cdot{\bf
a})=A^{-1}_{\sigma(i)}({\bf a})
$$
i.e.,
\begin{equation} \label{automorfisme_invers_3}
s^{-1}_i({\bf a})=s_{\sigma(i)}({\bf P}^{-1}({\bf a}))^{-1}
\end{equation}
for all ${\bf a}\in(\Natural^n)^*$ and $i=1,\ldots,n$.

To determine the group $N_n$, recall that it is the set of (iso)morphisms in ${\bf Aut}_{{\bf 2Mat}_{\Complex}}(n)$ with domain the identity
1-morphism ${\rm id}_n:n\To n$ in ${\bf 2Mat}_{\Complex}$. Hence, by Lemma~\ref{morfismes_invertibles}, the codomain 1-automorphism has rank matrix necessarily equal to ${\bf I}_n$. An element
in $N_n$ reduces then to a gauge $s$ specifying
this codomain, with $s_i({\bf a})\in{\rm GL}(a_i)$ for
all ${\bf a}\in\Natural^n$ and $i=1,\ldots,n$, together with a vector
$\mbox{\boldmath$\lambda$}=(\lambda_1,\ldots,\lambda_n)\in\Complex^n$
defining the 2-morphism ${\sf T}:({\bf I}_n,{\bf
I})\Rightarrow({\bf I}_n,s)$, with $\lambda_i\neq 0$ for all $i=1,\ldots,n$ for {\sf T} to be a 2-isomorphism. Because of the normalization condition on the
gauge, this gives a bijection of sets $\Theta:H_n\To(\Complex^*)^n\times{\bf GL}_n$ such that
\begin{equation} \label{identificacio_H_n}
\Theta\left(({\bf I}_n,{\bf I})\stackrel{{\sf T}}{\Rightarrow}({\bf I}_n,s)\right)
=(\mbox{\boldmath$\lambda$}^{{\sf T}},A^{({\bf I}_n,s)})
\end{equation}
where $\mbox{\boldmath$\lambda$}^{{\sf T}}$ is given by
\begin{equation} \label{iso_Theta}
\lambda_i^{{\sf T}}={\sf T}_{ii},\quad i=1,\ldots,n
\end{equation}
and $A^{({\bf I}_n,s)}$ is defined by {\rm (\ref{A_s})}. To see that this is a group morphism, recall that the group structure on $N_n$ is given by the
horizontal composition of 2-morphisms in ${\bf 2Mat}_{\Complex}$. Then, for any ${\sf T}:({\bf I}_n,{\bf I})\Rightarrow({\bf I}_n,s)$ and
$\widetilde{{\sf T}}:({\bf I}_n,{\bf I})\Rightarrow({\bf
I}_n,\widetilde{s})$, we have
$$
\Theta(\widetilde{{\sf T}}\circ{\sf T})=
(\mbox{\boldmath$\lambda$}^{\widetilde{{\sf T}}\circ{\sf
T}},A^{({\bf I}_n,\widetilde{s}\ast s)})
$$
But Equations~(\ref{iso_Theta}) and (\ref{composicio_horitzontal}) give that
$$
\lambda_i^{\widetilde{{\sf T}}\circ{\sf T}}
=(\widetilde{{\sf T}}\circ{\sf T})_{ii}
=\widetilde{{\sf T}}_{ii}\otimes{\sf T}_{ii}
=\lambda_i^{\widetilde{{\sf T}}}\lambda_i^{{\sf T}},\quad i=1,\ldots,n
$$
(the gauges are normalized and all nondiagonal
entries in $\widetilde{{\sf T}}$ and {\sf T} are empty), while (\ref{seccions_functor_composicio}) and the fact that all involved rank matrices are identities give
$$
A_i^{({\bf I}_n,\widetilde{s}\ast s)}({\bf a})=(\widetilde{s}\ast
s)_i({\bf a})=\widetilde{s}_i({\bf a})s_i({\bf a})=A_i^{({\bf
I}_n,\widetilde{s})}({\bf a})A_i^{({\bf I}_n,s)}({\bf a}),\quad {\bf
a}\in(\Natural^n)^*
$$
(in the second equality, we use the normalization
conditions on both the gauges and on the permutation matrices ${\bf
P}(\widetilde{{\bf R}}_k,{\bf R},{\bf a})$). Therefore
$$
\Theta(\widetilde{{\sf T}}\circ{\sf T})=(\mbox{\boldmath$\lambda$}^{\widetilde{{\sf
T}}}\mbox{\boldmath$\lambda$}^{{\sf T}},A^{({\bf
I}_n,\widetilde{s})}A^{({\bf
I}_n,s)})=(\mbox{\boldmath$\lambda$}^{\widetilde{{\sf T}}},A^{({\bf
I}_n,\widetilde{s})})(\mbox{\boldmath$\lambda$}^{{\sf T}},A^{({\bf
I}_n,s)})=\Theta(\widetilde{{\sf T}})\Theta({\sf T})
$$
Formula~(\ref{morfisme_parcial}) readily follows now from the commutative diagram (\ref{morfisme_modul_creuat_1}). It remains to be checked that the action $\lhd_n$ corresponds by these isomorphisms to that given by (\ref{accio}). Recall that the action of an object $({\bf
P},s)$ in $E_n$ on a morphism ${\sf T}:({\bf I}_n,{\bf
I})\Rightarrow({\bf I}_n,s')$ in $N_n$ is defined by
$$
({\bf P},s)\lhd_n{\sf T}=1_{({\bf P},s)}\circ{\sf T}\circ 1_{({\bf
P},s)^{-1}}
$$
Let us write $\widetilde{{\sf T}}=1_{({\bf P},s)}\circ{\sf T}\circ
1_{({\bf P},s)^{-1}}$ and let $\widetilde{s}$ be the gauge of the
codomain 1-automorphism of $\widetilde{{\sf T}}$, i.e.,
$\widetilde{s}=s\ast s'\ast s^{-1}$, where $s^{-1}$ is defined by
(\ref{automorfisme_invers_3}). Then, we need to see that
\begin{align}
\mbox{\boldmath$\lambda^{\widetilde{{\sf
        T}}}$}&=\sigma\cdot\mbox{\boldmath$\lambda^{{\sf T}}$} \label{fet1} \\ 
A^{({\bf I}_n,\widetilde{s})}&=A^{({\bf P},s)}(\sigma\rhd A^{({\bf
I}_n,s')})\left(A^{({\bf P},s)}\right)^{-1} \label{fet2}
\end{align}
where, for short, we have written $\sigma$ instead of $\sigma({\bf
P})$. To see (\ref{fet1}), notice that all involved rank
matrices are permutation matrices so that all gauge
terms appearing in (\ref{composicio_horitzontal}) are trivial. Hence, Equations~(\ref{iso_Theta}), (\ref{composicio_horitzontal}) and
(\ref{automorfisme_invers_2}) give
$$
\lambda_i^{\widetilde{{\sf T}}}=
({\bf 1}_{({\bf P},s)}\circ{\sf T}\circ {\bf 1}_{({\bf
P},s)^{-1}})_{ii}=\bigoplus_{j=1}^n\left(\bigoplus_{j'=1}^n\ ({\bf
1}_{({\bf P},s)})_{ij'}\otimes{\sf T}_{j'j}\right)\otimes({\bf
1}_{({\bf P}^{-1},s^{-1})})_{ji}
$$
But $({\bf 1}_{({\bf P}^{-1},s^{-1})})_{ji}$ is empty except in
case $P^{-1}_{ji}=\delta_{j\sigma^{-1}(i)}$ is non zero, in which
case it is equal to 1, and similarly with the term $({\bf 1}_{({\bf
P},s)})_{ij'}$. It follows that
$$
\lambda^{\widetilde{{\sf T}}}_i=\bigoplus_{j'=1}^n\ ({\bf 1}_{({\bf P},s)})_{ij'}\otimes{\sf
T}_{j'\sigma^{-1}(i)}={\sf
T}_{\sigma^{-1}(i)\sigma^{-1}(i)}=\lambda^{{\sf
T}}_{\sigma^{-1}(i)}
$$
Since this is true for all $i=1,\ldots,n$, we conclude that
$\mbox{\boldmath$\lambda^{\widetilde{{\sf T}}}$}
=\sigma\cdot\mbox{\boldmath$\lambda^{{\sf T}}$}$, as required.
Let us now prove (\ref{fet2}). Indeed, by Equation~(\ref{A_s}) and
(\ref{seccions_functor_composicio}) and because of the
normalization conditions on both the gauge $s$ and the permutation
matrices ${\bf P}(\widetilde{{\bf R}}_k,{\bf R},{\bf a})$ appearing
in (\ref{seccions_functor_composicio}), it is
$$
A^{({\bf I}_n,\widetilde{s})}_i({\bf a})=(s\ast s'\ast
s^{-1})_i({\bf a})=s_i({\bf P}^{-1}({\bf
a}))\left(\bigoplus_{j=1}^n\ {\bf I}_{P_{ij}}\otimes (s'\ast
s^{-1})_j({\bf a})\right)
$$
But $P_{ij}=\delta_{i\sigma(j)}$, so that the direct sum of
matrices reduces to the term $j=\sigma^{-1}(i)$ and we get
$$
A^{({\bf I}_n,\widetilde{s})}_i({\bf a})=s_i({\bf P}^{-1}({\bf a}))
(s'\ast s^{-1})_{\sigma^{-1}(i)}({\bf a})
$$
Now, using again (\ref{seccions_functor_composicio}) and the
normalization conditions, we have
$$
(s'\ast s^{-1})_{\sigma^{-1}(i)}({\bf a})=s'_{\sigma^{-1}(i)}({\bf
P}^{-1}({\bf a}))s_{\sigma^{-1}(i)}^{-1}({\bf
a})=s'_{\sigma^{-1}(i)}({\bf P}^{-1}({\bf a}))s_i({\bf P}^{-1}({\bf
a}))^{-1}
$$
where we have used (\ref{automorfisme_invers_3}) in the last
equality. Thus
\begin{align*}
A^{({\bf I}_n,\widetilde{s})}_i({\bf a})&=s_i({\bf P}^{-1}({\bf
a})) s'_{\sigma^{-1}(i)}({\bf P}^{-1}({\bf a}))s_i({\bf
P}^{-1}({\bf a}))^{-1} \\ &=A_i^{({\bf P},s)}({\bf
a})A_{\sigma^{-1}(i)}^{({\bf I}_n,s')}(\sigma^{-1}\cdot{\bf
a})A_i^{({\bf P},s)}({\bf a})^{-1} \\ &=\left(A^{({\bf
P},s)}(\sigma\rhd A^{({\bf I}_n,s')})(A^{({\bf
P},s)})^{-1}\right)_i({\bf a})
\end{align*}
and this is true for all ${\bf a}\in(\Natural^n)^*$ and
$i=1,\ldots,n$, as we wanted to prove.
\end{proof}

\begin{ex} \label{exemple_n=1} {\rm
For $n=1$, it is $E_1={\bf GL}_1=\prod_{q\geq 2}{\rm GL}(q)$ and $N_1=\Complex^*\times{\bf GL}_1$, with ${\bf GL}_1$ acting on
$\Complex^*\times{\bf GL}_1$ by componentwise conjugation on the second factor.
 }
\end{ex}

\begin{prop} \label{G(n)_E(n)}
If $n\geq 1$, ${\sf Equiv}_{{\bf 2Mat}_{\Complex}}(n)$ is equivalent to the split 2-group $\G(n)=S_n[0]\times (\Complex^*)^n[1]$, with $S_n$ acting on $(\Complex^*)^n$ according to Equation~(\ref{wreath_product}). An explicit equivalence $\EE(n):\G(n)\To{\sf Equiv}_{{\bf 2Mat}_{\Complex}}(n)$ is defined as follows: for all $\sigma,\sigma'\in S_n$ and $\mbox{\boldmath$\lambda$}\in(\Complex^*)^n$,
\begin{align} 
&\EE(n)(\sigma)=({\bf P}(\sigma),{\bf I})\ \mbox{(action on objects)} \label{accio_E(n)_sobre_objectes}
\\ &\EE(n)(\sigma,\mbox{\boldmath$\lambda$})_{\sigma(j),j}=\lambda_{\sigma(j)},\quad j=1,\ldots,n\ \mbox{(action on morphisms)}\ \label{accio_E(n)_sobre morfismes}
\\ &\EE(n)_2(\sigma,\sigma')=1_{({\bf P}(\sigma\sigma'),{\bf I})}\ \mbox{(structural isomorphisms)} \label{isos_estructurals_E(n)}
\end{align}
where {\bf I} denotes the trivial gauge (recall that, by definition of 2-morphisms in ${\bf 2Mat}_{\Complex}$, the entry $(i,j)$ of $\EE(n)(\sigma,\mbox{\boldmath$\lambda$}):({\bf P}(\sigma),{\bf I})\Rightarrow({\bf P}(\sigma),{\bf I})$ is empty except when $i=\sigma(j)$).
\end{prop}
\begin{proof}
Using (\ref{pi0_via_modul_creuat})-(\ref{accio_pi0_pi1_via_modul_creuat}) and the above description of the crossed module associated to ${\sf Equiv}_{{\bf 2Mat}_{\Complex}}(n)$, it easily follows that $\pi_0({\sf Equiv}_{{\bf 2Mat}_{\Complex}}(n))\cong S_n$ and $\pi_1({\sf Equiv}_{{\bf 2Mat}_{\Complex}}(n))\cong(\Complex^*)^n$, with the action of the first onto the second given by (\ref{wreath_product}). Furthermore, the 3-cocycle of ${\sf Equiv}_{{\bf 2Mat}_{\Complex}}(n)$ is cohomologically trivial because the section $s:S_n\To S_n\times_{\triangleright}{\bf GL}_n$ defined by $\sigma\mapsto(\sigma,{\bf I})$ is a group morphism. Hence, by the classification theorem for 2-groups, ${\sf Equiv}_{{\bf 2Mat}_{\Complex}}(n)$ is indeed equivalent to $\G(n)$. The fact that $\EE(n)$ defines an equivalence of 2-groups is an easy check left to the reader. 
\end{proof}

\subsection{Objects of ${\repg}_{{\bf 2Mat}_{\Complex}}(\G)$}

\label{subseccio_representacions_lineals}

Objects of ${\repg}_{{\bf 2Mat}_{\Complex}}(\G)$ are naturally classified by the following notion of dimension.

\begin{defn}
A 2-matrix representation of $\G$ is called of {\sl dimension} $n$
($n\geq 0$) when the objects of $\G$ are represented as
autoequivalences of the object $n$ of ${\bf 2Mat}_{\Complex}$. 
\end{defn}
 
There is only one zero-dimensional 2-matrix representation, given by
the unique functor $\G\To{\sf Equiv}_{{\bf 2Mat}_{\Complex}}(0)\simeq
{\bf 1}$. It will be denoted by $\OO$ and called the {\sl zero
  2-matrix representation}. For dimensions $n\geq 1$, there are
various (non equivalent) possibilities. The following description
readily follows from
Propositions~\ref{repr_modul_isomorfisme_cas_general} and
\ref{G(n)_E(n)}. 

\begin{thm} \label{representacions_modul_isomorfisme}
Let $\G$ be an arbitrary 2-group and $\alpha\in
Z^3(\pi_0(\G),\pi_1(\G))$ a fixed classifying 3-cocycle of $\G$. Then,
up to isomorphism, a non zero 2-matrix representation of $\G$ is
completely determined by a quadruple $(n,\rho,\beta,c)$, where $n\geq
1$ is the dimension of the representation, $\rho:\pi_0(\G)\To S_n$ a
morphism of groups, $\beta:\pi_1(\G)\To(\Complex^*)^n_{\rho}$ a
morphism of $\pi_0(\G)$-modules such that $[\beta\circ\alpha]=0$ in
$H^3(\pi_0(\G),(\Complex^*)^n_{\rho})$ and $c\in
C^2(\pi_0(\G),(\Complex^*)^n_{\rho})$ a normalized 2-cochain such that
$\partial c=\beta\circ\alpha$. Furthermore, any such quadruple
determines a non zero 2-matrix representation of $\G$ (up to
isomorphism). 
\end{thm} 

As pointed out before, this gives a non faithful description of the
set of isomorphism classes of representations. Thus, if $c,c'$ are
cohomologous, with $c\neq c'$, the quadruples $(n,\rho,\beta,c)$ and
$(n,\rho,\beta,c')$ determine different but isomorphic
representations. As discussed below
(cf. \S~\ref{classes_equivalencia_repr}), this is not the unique
source of non faithfulness of this description. 

\begin{rem} \label{remarca_2-cocadenes} {\rm  
As already pointed out, the piece of
$C^2(\pi_0(\G),(\Complex^*)^n_{\rho'})$ involved in this description
is non canonical. In particular, when $\G$ is split, we can always
assume the involved 2-cochains $c$ are actually 2-cocycles by choosing
as classifying 3-cocycle $\alpha$ the trivial one.     }
\end{rem}

For later use, observe that an explicit representation $\Ff(n,\rho,\beta,c):\underline{\G}\To{\bf 2Mat}_{\Complex}$ corresponding to the quadruple $(n,\rho,\beta,c)$ is that defined by the pair $(n,\F(\rho,\beta,c))$, where $\F(\rho,\beta,c)=(F(\rho,\beta,c),F(\rho,\beta,c)_2):\G\To{\sf Equiv}_{{\bf 2Mat}_{\Complex}}(n)$ is the monoidal functor given by the composite $\EE(n)\circ\F^{(\rho,\beta,c)}$, with $\EE(n)$ the equivalence of Proposition~\ref{G(n)_E(n)} (cf. Equation~(\ref{correspondencia1})). Explicitly:
\begin{itemize}
\item
For any object $g\in\pi_0(\G)$, $F(\rho,\beta,c)(g):n\To n$ is the autoequivalence of $n$
\begin{equation}
F(\rho,\beta,c)(g)=({\bf P}(\rho(g)),{\bf I}) \label{accio_representacio_sobre_objectes} 
\end{equation}
\item
For any morphism $\varphi\in{\rm Aut}_{\G}(g)$, $F(\rho,\beta,c)(\varphi):({\bf P}(\rho(g)),{\bf I})\Rightarrow({\bf P}(\rho(g)),{\bf I})$ is the 2-automorphism whose non empty entries are
\begin{equation} 
F(\rho,\beta,c)(\varphi)_{\rho(g)(j),j}=\beta(\gamma_g^{-1}(\varphi))_{\rho(g)(j)},\quad j=1,\ldots,n, \label{accio_representacio_sobre_morfismes}
\end{equation}
\item
For all objects $g_1,g_2\in\pi_0(\G)$, $F(\rho,\beta,c)_2(g_1,g_2):({\bf P}(\rho(g_1g_2)),{\bf I})\Rightarrow({\bf P}(\rho(g_1g_2)),{\bf I})$ is the 2-automorphism with non empty entries given by 
\begin{equation}
F(\rho,\beta,c)_2(g_1,g_2)_{\rho(g_1g_2)(j),j}=c(g_1,g_2)_{\rho(g_1g_2)(j)} \label{isos_estructurals_representacio}
\end{equation}
for all $j=1,\ldots,n$. 
\end{itemize}

We shall refer to these representations as the {\sl gauge trivial} 2-matrix representations, because all 1-morphisms $F(\rho,\beta,c)(g):n\To n$, for $g\in\pi_0(\G)$, are of trivial gauge. Clearly, composing the morphisms $\F^{(\rho,\beta,c)}$ with equivalences $\EE'$ other than the $\EE(n)$ of Proposition~\ref{G(n)_E(n)}, we would obtain equivalent but non gauge trivial representations. The representation will be called {\sl strict} when the 2-cochain $c$ is trivial (hence, when the associated pseudofunctor $\Ff(n,\rho,\beta,c)$ is in fact a 2-functor).

\begin{exs} \label{exemples_representacions} {\rm 
(1) For any 2-group $\G$ and given any classifying 3-cocycle $\alpha$ of $\G$, there is a 1-dimensional 2-matrix representation of $\G$ for each pair $(\beta,c)$, with $\beta:\pi_1(\G)\To\Complex^*$ any $\pi_0(\G)$-invariant character of $\pi_1(\G)$ such that $[\beta\circ\alpha]=0$ in $H^3(\pi_0(\G),\Complex^*)$ and $c\in C^2(\pi_0(\G),\Complex^*)$ any normalized 2-cochain such that $\partial c=\beta\circ\alpha$ (here, $\Complex^*$ stands for the trivial $\pi_0(\G)$-module). Conversely, any 1-dimensional 2-matrix representation of $\G$ is isomorphic to one of this form for some pair $(\beta,c)$. Such representations are denoted by $\I_{\beta,c}$. If $\beta$ is trivial (i.e., $\beta(u)=1$ for all $u\in\pi_1(\G)$), so that $c$ is necessarily a 2-cocycle $z$, we shall write $\I_z$ and just $\I$ if $z$ is further the trivial 2-cocycle. $\I$ will be called the {\sl trivial representation} of $\G$. Note that in case $\G$ is split, by choosing as classifying 3-cocycle the trivial one, the pairs $(\beta,c)$ become pairs $(\beta,z)$ with $\beta:\pi_1(\G)\To\Complex^*$ any $\pi_0(\G)$-invariant character of $\pi_1(\G)$ and $z\in Z^2(\pi_0(\G),\Complex^*)$ any normalized 2-cocycle.

(2) For any dimension $n\geq 1$ and any 2-group $\G$, we have the {\sl purely permutational 2-matrix representations}, defined by quadruples $(n,\rho,\beta,c)$ with $\rho:\pi_0(\G)\To S_n$ any group morphism, $\beta$ the constant map $\beta(u)=(1,\stackrel{n)}{\ldots},1)$ for all $u\in\pi_1(\G)$ and $c$ the trivial 2-cocycle. Such representations are denoted by $\R_{n,\rho}$ or just by $\R_n$ when $\rho$ is trivial. They include the trivial representation $\I$ (it is $\I=\R_1$).

(3) For any dimension $n\geq 1$ and any 2-group $\G$, we have the {\sl purely cocyclic 2-matrix representations}, defined by quadruples $(n,\rho,\beta,c)$ with $\rho:\pi_0(\G)\To S_n$ the trivial group morphism, $\beta$ the constant map $\beta(u)=(1,\stackrel{n)}{\ldots},1)$ for all $u\in\pi_1(\G)$ and $c=z$ any normalized 2-cocycle of $\pi_0(\G)$ with values in the trivial $\pi_0(\G)$-module $(\Complex^*)^n$. Such representations are denoted by $\R_{n,z}$. They include the representations $\I_z$ and $\R_n$, corresponding respectively to $\R_{1,z}$ and $\R_{n,z}$ with $z$ trivial.
}
\end{exs}

\begin{rem} \label{remarca_fidelitat_representacions} {\rm
By considering strict representations of a strict 2-group $\G$, Barrett and Mackaay pointed out the paucity of examples of such representations except for a $\G$ whose group of objects $\G_0$ is profinite. In particular, as mentioned in the introduction, they noticed that, except for such $\G$, the set of these representations is not even collectively faithful (with respect to objects). It is worth being careful at this point, however, because such a lack of collective faithfulness is only true in the weak sense, i.e., in the sense that there exists $g_1,g_2\in\pi_0(\G)$, with $g_1\neq g_2$, such that the functors $F(g_1)$ and $F(g_2)$ are always {\it isomorphic} (not necessarily equal) for any representation $(n,\F)$. Thus, using non gauge trivial strict representations, it turns out that there exists even a 1-dimensional 2-matrix strict representation $\F:\G\To{\rm Equiv}_{{\bf 2Mat}_{\Complex}}(1)$ which is faithful on objects in the strict sense whenever $\G$ is such that $\G_0$ is any compact Lie group, not necessarily profinite. Indeed, it is easy to see from the description of ${\rm Equiv}_{{\bf 2Mat}_{\Complex}}(n)$ in the previous section that a 1-dimensional 2-matrix strict representation of $\G$ is equivalent to a pair $(\widetilde{\rho},\chi)$, with $\widetilde{\rho}:\G_0\To{\bf GL}_1=\prod_{q>1}\ {\rm GL}(q)$ a group morphism and $\chi$ a character of the fiber group $H\subset\G_1$ which is constant over the $\G_0$-orbits. The strict monoidal functor $F$ corresponding to such a pair is given on objects by $F(g)=({\id}_1,s^{\widetilde{\rho}(g)})$, with $s^{\widetilde{\rho}(g)}$ the gauge given by $s^{\widetilde{\rho}(g)}(q)=\widetilde{\rho}_q(g)$, $q\geq 2$. Peter-Weyl's theorem then ensures the existence of the above mentioned ``faithful'' representation. }
\end{rem}

In all what follows, by an $n$-dimensional 2-matrix representation ($n\geq 1$) we shall always mean one of the previous gauge trivial 2-matrix representations. As it has been argued, this implies no loss of generality while it considerably simplifies some equations. Such a representation is denoted by $(n,\rho,\beta,c)$.

\subsection{Categories of morphisms}

\label{subseccio_categories_de_morfismes}

We wish now to describe the categories of morphisms between two arbitrary (gauge trivial) 2-matrix representations. To simplify notation,  we shall write ${\bf Hom}(\Ff,\Ff')$ instead of ${\bf Hom}_{\repg_{{\bf 2Mat}_{\Complex}}(\G)}(\Ff,\Ff')$ to denote the category of morphisms between two given representations $\Ff,\Ff':\underline{\G}\To{\bf 2Mat}_{\Complex}$.

Suppose $\Ff,\Ff'$ are of dimensions $n$ and $n'$, respectively, $n,n'\geq 0$. According to Proposition~\ref{representacio_cas_general}
\begin{itemize}
\item[(a)]
a 1-intertwiner $\xi:\Ff\To\Ff'$ consists of
\begin{itemize}
\item
a 1-morphism $({\bf R},s):n\To n'$ in ${\bf 2Mat}_{\Complex}$, and
\item
a family of 2-isomorphisms $\Phi=\{\Phi(g):\Ff'(g)\circ ({\bf R},s)\Rightarrow({\bf R},s)\circ\Ff(g)\}$ in ${\bf 2Mat}_{\Complex}$ indexed by the objects of $\G$, natural in $g$ and coherent,
\end{itemize}
\item[(b)]
a 2-intertwiner $\ngg:\xi\Rightarrow\overline{\xi}$ between two such 1-intertwiners $\xi=({\bf R},s,\Phi)$ and $\overline{\xi}=(\overline{{\bf R}},\overline{s},\overline{\Phi})$ consists of a 2-morphism ${\sf T}:({\bf R},s)\Rightarrow(\overline{{\bf R}},\overline{s})$ in ${\bf 2Mat}_{\Complex}$ satisfying a naturality axiom.
\end{itemize}

By the {\sl support} of a 1-intertwiner $({\bf R},s,\Phi):\Ff\To\Ff'$ we shall mean the support of its rank matrix {\bf R}, defined by
$$
{\rm Sup}({\bf R}):=\{ (i',i)\in\{1,\ldots,n'\}\times\{1,\ldots,n\}\ |\
R_{i,i'}\neq 0\} 
$$
${\bf Hom}(\Ff,\Ff')$ contains at least the {\sl zero 1-intertwiner}, characterized by the condition ${\rm Sup}({\bf R})=\emptyset$ and given by $({\bf R},s)={\bf 0}_{n,n'}$ and $\Phi(g)$ equal to the $n'\times n$ matrix with all entries equal to the empty matrix. Similarly, for any two 1-intertwiners $\xi=({\bf R},s,\Phi)$ and $\xi'=({\bf R}',s',\Phi')$, we always have the corresponding {\sl zero 2-intertwiner}, given by the matrix all of whose non empty entries are zero matrices. These are the only possible 1- and 2-intertwiners when one of the representations $\Ff,\Ff'$ is the zero representation $\OO$. More precisely: 

\begin{prop}
Let $(n,\rho,\beta,c)$ be any non zero 2-matrix representation. Then, the categories ${\bf Hom}(\OO,\OO)$, ${\bf Hom}(\OO,(n,\rho,\beta,c))$ and ${\bf Hom}((n,\rho,\beta,c),\OO)$ are terminal.
\end{prop}
The proof is an easy check left to the reader. In fact, we shall see in a moment that there are other categories of morphisms which are terminal (cf. Corollary~\ref{cas_trivial}).

Let now $(n,\rho,\beta,c)$, $(n',\rho',\beta',c')$ be any two non zero 2-matrix representations, $n,n'\geq 1$. The set of possible 1-intertwiners between them turns out to depend on some data easily computable from both quadruples. Concretely, given the pair $n,n'$, we have the cartesian product set $M(n',n)\equiv\{1,\ldots,n'\}\times\{1,\ldots,n\}$. Then, morphisms $\rho,\rho'$ induce a right action of $\pi_0(\G)$ on this set given by
$$
(i',i)\cdot g=(\rho'(g)^{-1}(i'),\rho(g)^{-1}(i)),\quad (i',i)\in M(n',n),\ g\in\pi_0(\G).
$$
This gives rise to a partition of $M(n',n)$ into $\pi_0(\G)$-orbits. In general, only a subset of these orbits, determined by $\beta$ and $\beta'$, will be able to support a 1-intertwiner. Such orbits are defined as follows. For any maps $\lambda:\pi_1(\G)\To(\Complex^*)^n$ and $\lambda':\pi_1(\G)\To(\Complex^*)^{n'}$, set
$$
I(\lambda,\lambda'):=\{(i',i)\in M(n',n)\ |\ \beta'_{i'}=\beta_i\}
$$
Let us call $(i',i)\in M(n',n)$ an {\sl intertwining point} for $\lambda,\lambda'$ when $(i',i)\in I(\lambda,\lambda')$. When $\lambda,\lambda'$ are the maps $\beta,\beta'$ in the above quadruples, we have the following:

\begin{lem}
If $(i',i)\in M(n',n)$ is an intertwining point for $\beta,\beta'$, the same is true for all points in its $\pi_0(\G)$-orbit (i.e., all points $(\rho'(g)(i'),\rho(g)(i))$, $g\in\pi_0(\G)$).
\end{lem}   
\begin{proof}
Suppose $(i',i)\in I(\beta,\beta')$ so that $\beta'_{i'}=\beta_i$. Since $\beta,\beta'$ are morphisms of $\pi_0(\G)$-modules when $\pi_0(\G)$ acts on $(\Complex^*)^n$ and $(\Complex^*)^{n'}$ via $\rho$ and $\rho'$, respectively, we have
$$
\beta'_{\rho'(g)(i')}(u)=(g^{-1}\cdot\beta'(u))_{i'}=\beta'_{i'}(g^{-1}\cdot u)=\beta_i(g^{-1}\cdot u)=(g^{-1}\cdot\beta(u))_i=\beta_{\rho(g)(i)}(u)
$$
for all $u\in\pi_1(\G)$ and $g\in\pi_0(\G)$. Hence, $(\rho'(g)(i'),\rho(g)(i))\in I(\beta,\beta')$ for any $g\in\pi_0(\G)$.
\end{proof}
This shows that $I(\beta,\beta')$ is exactly the union of some of the $\pi_0(\G)$-orbits of $M(n',n)$, and this leads to the following definition.
\begin{defn}
Given representations $(n,\rho,\beta,c)$ and $(n',\rho',\beta',c')$ as above, a $\pi_0(\G)$-orbit $\Oo$ of $M(n',n)$ is called {\sl intertwining} when $\Oo\subseteq I(\beta,\beta')$. 
\end{defn}

Finally, given any intertwining orbit $\Oo$, a 2-cocycle of $\pi_0(\G)$ can be built from the last pair $c,c'$ of the quadruple as follows. Let $\Oo=\{(i'_1,i_1),\ldots,(i'_k,i_k)\}$ and let $\Ff(\Oo,\Complex^*)_{\rho,\rho'}$ be the multiplicative abelian group of all functions $\lambda:\Oo\To\Complex^*$ equipped with the $\pi_0(\G)$-module structure
$$
(g\cdot\lambda)(i'_j,i_j)=\lambda((i'_j,i_j)\cdot g)=\lambda(\rho'(g)^{-1}(i'_j),\rho(g)^{-1}(i_j))
$$
for all $j=1,\ldots,k$, $g\in\pi_0(\G)$. Then, define $z_{\Oo}:\pi_0(\G)\times\pi_0(\G)\To\Ff(\Oo,\Complex^*)_{\rho,\rho'}$ by
\begin{equation} \label{definicio_cocicles}
z_{\Oo}(g_1,g_2)(i'_j,i_j)=\frac{c'(g_1,g_2)_{i'_j}}{c(g_1,g_2)_{i_j}}, \quad g_1,g_2\in\pi_0(\G),\ j=1,\ldots,k
\end{equation}

\begin{lem}
If the $\pi_0(\G)$-orbit $\Oo$ is intertwining, it is $z_{\Oo}\in Z^2(\pi_0(\G),\Ff(\Oo,\Complex^*)_{\rho,\rho'})$. Furthermore, the cohomology class of $z_{\Oo}$ depends only on the cohomology classes of $c$ and $c'$.
\end{lem}
\begin{proof}
For any $j=1,\ldots,k$, we have
\begin{align*}
(\partial z_{\Oo})(g_1,g_2,g_3)(i'_j,i_j)&=\left(\frac{c'(g_2,g_3)_{\rho'(g_1)^{-1}(i'_j)}}{c(g_2,g_3)_{\rho(g_1)^{-1}(i_j)}}\right) \left(\frac{c'(g_1g_2,g_3)^{-1}_{i'_j}}{c(g_1g_2,g_3)^{-1}_{i_j}}\right) \left(\frac{c'(g_1,g_2g_3)_{i'_j}}{c(g_1,g_2g_3)_{i_j}}\right) \left(\frac{c'(g_1,g_2)^{-1}_{i'_j}}{c(g_1,g_2)^{-1}_{i_j}}\right) \\ 
&=\frac{(\partial c')(g_1,g_2,g_3)_{i'_j}}{(\partial c)(g_1,g_2,g_3)_{i_j}}
\end{align*}
Furthermore, $c,c'$ are such that $\partial c=\beta\circ\alpha$ and $\partial c'=\beta'\circ\alpha$. Hence
\begin{align*}
(\partial c)(g_1,g_2,g_3)_i&=\beta_i(\alpha(g_1,g_2,g_3)),\quad i\in\{1,\ldots,n\}
\\ (\partial c')(g_1,g_2,g_3)_{i'}&=\beta'_{i'}(\alpha(g_1,g_2,g_3)),\quad i'\in\{1,\ldots,n'\}
\end{align*}
for all $g_1,g_2,g_3\in\pi_0(\G)$. But $\beta'_{i'_j}=\beta_{i_j}$ for all $j=1,\ldots,k$ because $\Oo$ is intertwining. Therefore
$$
(\partial c)(g_1,g_2,g_3)_{i_j}=(\partial c')(g_1,g_2,g_3)_{i'_j},\quad j=1,\ldots,k,
$$
from which it follows that $(\partial z_{\Oo})(g_1,g_2,g_3)(i'_j,i_j)=1$ for all $j=1,\ldots,k$, as required. To prove the last assertion, let $\tilde{c}=c+\partial x$ and $\tilde{c}'=c'+\partial x'$ for some 1-cochains $x,x'$ (althought $\Ff(\Oo,\Complex^*)_{\rho,\rho'}$ is a multiplicative group, we keep using additive notation for the group law between cochains). We leave to the reader to check that the new 2-cocycle $\tilde{z}_{\Oo}$ obtained using $\tilde{c}$ and $\tilde{c}'$ is $\tilde{z}_{\Oo}=z_{\Oo}+\partial(x'-x)$. 
\end{proof}

We can now describe the set of 1-intertwiners between any two non zero 2-matrix representations. The simplest situation occurs when the representations are such that no $\pi_0(\G)$-orbit of $M(n',n)$ is intertwining, i.e., when $I(\beta,\beta')=\emptyset$. In such cases, there is no 1-intertwiner except the zero one. This is a consequence of the next result.

\begin{prop} \label{cond_admissibilitat_matriu_rangs}
Let $(n,\rho,\beta,c)$, $(n',\rho',\beta',c')$ be any two non zero 2-matrix representations. Then, if $({\bf R},s):n\To n'$ is the 1-morphism of a 1-intertwiner $({\bf R},s,\Phi):(n,\rho,\beta,c)\To(n',\rho',\beta',c')$, the following holds:
\begin{itemize}
\item[(i)]
{\bf R} is $(\rho,\rho')$-{\sl invariant}, i.e.
\begin{equation} \label{invariancia_R}
R_{\rho'(g)(i'),\rho(g)(i)}=R_{i'i}, \quad \forall\ i'=1,\ldots,n',\
i=1,\ldots,n
\end{equation}
\item[(ii)]
${\rm Sup}({\bf R})$ is entirely made of intertwining $\pi_0(\G)$-orbits, i.e.,
\begin{equation} \label{condicio_sobre_suport}
{\rm Sup}({\bf R})\subseteq I(\beta,\beta')
\end{equation}
\end{itemize}
In particular, the support of any 1-intertwiner is a union of some (or all) intertwining $\pi_0(\G)$-orbits of $M(n',n)$.
\end{prop} 
\begin{proof}
Suppose there exists a 1-intertwiner whose 1-morphism is $({\bf R},s):n\To n'$. In particular, this implies that for any $g\in\pi_0(\G)$ there is a 2-isomorphism $\Phi(g):(\rho'(g)){\bf R},s)\Rightarrow({\bf P}({\bf RP}(\rho(g)),s)$ in ${\bf 2Mat}_{\Complex}$ (cf. Equation~(\ref{accio_representacio_sobre_objectes})). By Lemma~\ref{morfismes_invertibles}-$(ii)$, this is possible if and only if
${\bf RP}(\rho(g))={\bf P}(\rho'(g)){\bf R}$ for all $g\in\pi_0(\G)$. Using (\ref{matriu_permutacio_sigma}), the reader may easily check that
this is equivalent to the conditions
$$
R_{\rho'(g)^{-1}(i'),i}=R_{i',\rho(g)(i)}
$$
for all $i=1,\ldots,n$ and $i'=1,\ldots,n'$, which are indeed equivalent to conditions (\ref{invariancia_R}). This proves (i). To prove (ii), notice that, by hypothesis, the 2-isomorphisms $\Phi(g):({\bf RP}(\rho(g)),s)\Rightarrow({\bf P}(\rho'(g)){\bf R},s)$ are natural in $g$, i.e., they make commutative the diagrams
$$
\xymatrix{
({\bf P}(\rho'(g)){\bf R},s)\ar[r]^{\Phi(g)}\ar[d]_{F(n',\rho',\beta',c')(\varphi)\circ 1_{({\bf R},s)}} & ({\bf RP}(\rho(g)),s)\ar[d]^{1_{({\bf R},s)}\circ F(n,\rho,\beta,c)(\varphi)} \\ 
({\bf P}(\rho'(g)){\bf R},s)\ar[r]_{\Phi(g)} & ({\bf RP}(\rho(g)),s)
}
$$
for all morphisms $\varphi:g\To g$. An easy computation using Equations~(\ref{composicio_vertical}), (\ref{composicio_horitzontal}) and (\ref{accio_representacio_sobre_morfismes}) shows that this is equivalent to the conditions
$$
\Phi(g)_{i'i}\beta(\gamma_g^{-1}(\varphi))_{\rho(g)(i)}\ {\bf I}_{R_{i',\rho(g)(i)}}=\beta'(\gamma^{-1}_g(\varphi))_{i'}\Phi(g)_{i'i}\ {\bf I}_{R_{\rho'(g)^{-1}(i'),i}}
$$
for all pairs $(i',i)\in\{1,\ldots,n'\}\times\{1,\ldots,n\}$ such that $R_{i',\rho(g)(i)}\neq 0$ (otherwise, the corresponding entries of the involved 2-morphisms are empty). Since $(\Complex^*)^n$ is abelian, both terms $\Phi(g)_{i'i}$ cancel out each other and we get the conditions
$$
\beta(\gamma_g^{-1}(\varphi))_{\rho(g)(i)}=\beta'(\gamma^{-1}_g(\varphi))_{i'}
$$ 
In particular, if $\varphi=u:e\To e$, this reduces to $\beta(u)_i=\beta'(u)_{i'}$ for all pairs $(i',i)$ such that $R_{i'i}\neq 0$, which is exactly condition (\ref{condicio_sobre_suport}). Finally, last assertion is an immediate consequence of $(i)$ and $(ii)$.
\end{proof}

\begin{cor} \label{cas_trivial}
For any non zero 2-matrix representations $(n,\rho,\beta,c)$ and $(n',\rho',\beta',c')$  such that $I(\beta,\beta')=\emptyset$ the category ${\bf Hom}((n,\rho,\beta,c), (n',\rho',\beta',c'))$ is terminal.
\end{cor}
Before considering the non trivial case $I(\beta,\beta')\neq\emptyset$, let us introduce a few definitions.

\begin{defn}
Given triples $(n,\rho,\beta)$ and $(n',\rho',\beta')$ as before, by an $\left(\begin{array}{c} n\ \rho\ \beta \\ n'\rho'\beta'\end{array}\right)$-{\sl admissible} rank matrix we mean a rank matrix ${\bf R}\in{\rm Mat}_{n'\times n}(\Natural)$ satisfying conditions (\ref{invariancia_R}) and (\ref{condicio_sobre_suport}). Moreover, if ${\bf R}\in{\rm Mat}_{n'\times n}(\Natural)$ is a $(\rho,\rho')$-invariant rank matrix and $g\in\pi_0(\G)$, a pair $(i',i)\in\{1,\ldots,n'\}\times\{1,\ldots,n\}$ is called $({\bf R},g)$-{\sl admissible} if $R_{i',\rho(g)(i)}\neq 0$ (equivalently, if $R_{\rho'(g)^{-1}(i'),i}\neq 0$). The set of such points is denoted by ${\rm Sup}({\bf R},g)$.
\end{defn}
Observe that, according to this definition, if $({\bf R},s,\Phi)$ is the triple describing a 1-intertwiner $\xi:(n,\rho,\beta,c)\To(n',\rho',\beta',c')$, the entry $\Phi(g)_{i'i}$ of the matrix $\Phi(g)$ is non empty if and only if $(i',i)\in {\rm Sup}({\bf R},g)$. The following result is an immediate consequence of the definitions.

\begin{lem}
Let {\bf R} be a $(\rho,\rho')$-invariant rank matrix. Then, the following conditions on a pair $(i',i)$ are equivalent:
\begin{itemize}
\item[(i)] $(i',i)\in{\rm Sup}({\bf R})$.
\item[(ii)] $(i',\rho(g)^{-1}(i))\in {\rm Sup}({\bf R},g)$.
\item[(iii)] $(\rho'(g)(i'),i)\in {\rm Sup}({\bf R},g)$.
\end{itemize}
\end{lem}  
A first description of the set of 1-intertwiners between a pair of non zero 2-matrix representations having a non empty set of intertwining $\pi_0(\G)$-orbits is then as follows.

\begin{prop} \label{primera_descripcio_1-intertwiners}
Let $(n,\rho,\beta,c)$ and $(n',\rho',\beta',c')$ be two non zero 2-matrix representations such that $I(\beta,\beta')\neq\emptyset$. Then, there is a 1-1 correspondence between the set of non zero 1-intertwiners $\xi:(n,\rho,\beta,c)\To(n',\rho',\beta',c')$ and the set of triples $({\bf R},s,{\sf S})$, where {\bf R} is a non zero $\left(\begin{array}{c} n\ \rho\ \beta \\ n'\rho'\beta'\end{array}\right)$-admissible rank matrix, $s$ is a (normalized) gauge for this rank matrix, i.e.
$$
s_{i'}({\bf a})\in{\rm GL}({\bf R}({\bf a})_{i'}),\quad {\bf
a}\in\Natural^n\setminus\{{\bf e}_1,\ldots,{\bf e}_n\},\ i'=1,\ldots,n',
$$ 
and ${\sf S}=\{{\sf S}_{i',i}:\pi_0(\G)\To{\rm GL}(R_{i',i}),\ (i',i)\in{\rm Sup}({\bf R})\}$ is a collection of maps such that 
\begin{align}
{\sf S}_{i',i}(e)&={\bf I}_{R_{i'i}} \label{condicio_normalitzacio}
\\
{\sf S}_{i',i}(g_1g_2)&=\frac{c(g_1,g_2)_{i}}{c'(g_1,g_2)_{i'}}\ {\sf S}_{i',i}(g_1)\
{\sf
S}_{\rho'(g_1)^{-1}(i'),\rho(g_1)^{-1}(i)}(g_2) \label{condicio_morfisme_sobre_T} 
\end{align}
for all $(i',i)\in{\rm Sup}({\bf R})$, $g_1,g_2\in\pi_0(\G)$.
\end{prop}

\begin{proof}
Let a triple $({\bf R},s,{\sf S})$ as above correspond to the 1-intertwiner $({\bf R},s,\Phi)$ with $\Phi(g)$, for any $g\in\pi_0(\G)$, defined by
\begin{equation} \label{def_phi_vs_S}
\Phi(g)_{i',\rho(g)^{-1}(i)}={\sf S}_{i',i}(g),\quad \forall\ (i',i)\in {\rm Sup}({\bf R})
\end{equation}
These matrices define a 2-isomorphism $\Phi(g):({\bf P}(\rho'(g)){\bf R},s)\Rightarrow({\bf RP}(\rho(g)),s)$ which is natural in $g$ because of the admissibility condition (\ref{condicio_sobre_suport}) on {\bf R}. Indeed, we have seen before (cf. proof of Proposition~\ref{cond_admissibilitat_matriu_rangs}) that naturality in $g$ is equivalent to conditions
$$
\beta_{\rho(g)(i)}(\gamma_g^{-1}(\varphi))=\beta'_{i'}(\gamma^{-1}_g(\varphi))
$$
for all $(i',i)\in {\rm Sup}({\bf R},g)$, all morphisms $g\stackrel{\varphi}{\To} g$ and all $g\in\pi_0(\G)$. But for any $g\in\pi_0(\G)$, it is $(i',i)\in {\rm Sup}({\bf R},g)$ if and only if $(i',\rho(g)(i))\in{\rm Sup}({\bf R})$. Since ${\rm Sup}({\bf R})\subseteq I(\beta,\beta')$, it follows that $\beta'_{i'}=\beta_{\rho(g)(i)}$ for all such pairs $(i',i)\in{\rm Sup}({\bf R},g)$. It remains to see that the coherence conditions on the 2-isomorphisms $\Phi(g)$ (cf. Proposition~\ref{representacio_cas_general}) are equivalent to conditions (\ref{condicio_normalitzacio}) and (\ref{condicio_morfisme_sobre_T}) on the maps ${\sf S}_{i',i}$. Let us consider the first coherence condition, which in this case reads
\begin{align}
\Phi(g_1g_2)\cdot&(F(n',\rho',\beta',c')_2(g_1,g_2)\circ 1_{({\bf R},s)})= \nonumber\\ 
&=(1_{({\bf R},s)}\circ F(n,\rho,\beta,c)_2(g_1,g_2))\cdot(\Phi(g_1)\circ 1_{({\bf P}(\rho(g_2)),{\bf I})})\cdot(1_{({\bf P}(\rho'(g_1)),{\bf I})}\circ\Phi(g_2)) \label{coherencia1}
\end{align}
Actually, this is a set of as many equations as elements $(i',i)$ there are in ${\rm Sup}({\bf R})$. Explicitly, let us identify the pair $(i',i)\in {\rm Sup}({\bf R})$ with the pair $(i',\rho(g_1g_2)^{-1}(i))\in{\rm Sup}({\bf R},g_1g_2)$. Then, the equation corresponding to this pair is the following. From (\ref{composicio_horitzontal}), we have
$$
(F(n',\rho',\beta',c')_2(g_1,g_2)\circ 1_{({\bf R},s)})_{i',\rho(g_1g_2)^{-1}(i)}=\bigoplus_{j=1}^n \left[F(n',\rho',\beta',c')_2(g_1,g_2)_{i',j}\otimes(1_{({\bf R},s)})_{j,\rho(g_1g_2)^{-1}(i)}\right]
$$
(all gauge terms are identities because of the normalization condition). But all entries $F(n',\rho',\beta',c')_2(g_1,g_2)_{i',j}$ are empty except when $j=\rho'(g_1g_2)^{-1}(i')$, in which case it is equal to $c'(g_1,g_2)_{i'}$ (cf. Equation~(\ref{isos_estructurals_representacio})). Hence
$$
(F(n',\rho',\beta',c')_2(g_1,g_2)\circ 1_{({\bf R},s)})_{i',\rho(g_1g_2)^{-1}(i)}=c'(g_1,g_2)_{i'}\ {\bf I}_{R_{\rho'(g_1g_2)^{-1}(i'),\rho(g_1g_2)^{-1}(i)}}
$$
Similarly, it is
\begin{eqnarray*}
&(1_{({\bf R},s)}\circ F(n,\rho,\beta,c)_2(g_1,g_2))_{i',\rho(g_1g_2)^{-1}(i)}=c(g_1,g_2)_{i}\ {\bf I}_{R_{i',i}}
\\
&(\Phi(g_1)\circ 1_{({\bf P}(\rho(g_2)),{\bf I})})_{i',\rho(g_1g_2)^{-1}(i)}=\Phi(g_1)_{i',\rho(g_1)^{-1}(i)}
\\ &(1_{({\bf P}(\rho'(g_1)),{\bf I})}\circ\Phi(g_2))_{i',\rho(g_1g_2)^{-1}(i)}=\Phi(g_2)_{\rho'(g_1)^{-1}(i'),\rho(g_1g_2)^{-1}(i)}
\end{eqnarray*}
Putting all these equalities together according to (\ref{composicio_vertical}), we obtain for the $(i',\rho(g_1g_2)^{-1}(i))$- component of (\ref{coherencia1}) the equation
$$
c'(g_1,g_2)_{i'}\ \Phi(g_1g_2)_{i',\rho(g_1g_2)^{-1}(i)}=c(g_1,g_2)_i\ \Phi(g_1)_{i',\rho(g_1)^{-1}(i)}\ \Phi(g_2)_{\rho'(g_1)^{-1}(i'),\rho(g_1g_2)^{-1}(i)}
$$
which is indeed (\ref{condicio_morfisme_sobre_T}) when written in terms of the maps ${\sf S}_{i',i}$. As regards the second coherence condition, in this case simplifies to
$\Phi(e)=1_{({\bf R},s)}$ because both morphisms $\F(\rho,\beta,c)$ and $\F(\rho',\beta',c')$ are special, and this is clearly equivalent to (\ref{condicio_normalitzacio}) when written in terms of the maps ${\sf S}_{i',i}$.
\end{proof}

\begin{defn} \label{definicio_1-intertwiner_estricte}
The data $({\bf R},s,{\sf S})$ defining a 1-intertwiner
$\xi:(n,\rho,\beta,c)\To(n',\rho',\beta',c')$ will respectively be called {\sl rank matrix},
{\sl gauge} and {\sl weakening maps} of $\xi$. The 1-intertwiner is called
{\sl gauge trivial} when $s={\bf I}$, the trivial gauge defined in
Section \ref{2cat_2vect}, and it is called {\sl strict} when all its weakening maps
are trivial (i.e., ${\sf S}_{i',i}(g)={\bf I}_{R_{i',i}}$ for all $g\in\pi_0(\G)$ and $(i',i)\in{\rm Sup}({\bf R})$; such 1-intertwiners correspond to 2-natural transformations between the involved pseudofunctors as defined in Definition~\ref{def_transf_pseudonatural}).
\end{defn}

Note that there is no condition on the gauge $s$ in the previous description of 1-intertwiners. This is not surprising because, according to
Lemma~\ref{morfismes_invertibles}, changing $s$ for any other gauge
$s'$ just give rise to a new 1-morphism $({\bf R},s'):n\To n'$ in
${\bf 2Mat}_{\Complex}$ which is 2-isomorphic to $({\bf R},s)$. It follows that if a triple $({\bf R},s,{\sf S})$ defines a 1-intertwiner, any other triple of the form $({\bf R},s',{\sf S})$ defines a new 1-intertwiner 2-isomorphic to the original one. We conclude the following

\begin{cor} \label{cor_1-intertwiners_modul_isomorfisme}
Any non zero 1-intertwiner $\xi:(n,\rho,\beta,c)\To(n',\rho',\beta',c')$ (in case $I(\beta,\beta')\neq\emptyset$) is isomorphic to one which is gauge trivial and hence, it is completely determined, up to isomorphism, by a pair $({\bf R},{\sf S})$, with {\bf R} and {\sf S} as in Proposition~\ref{primera_descripcio_1-intertwiners}.
\end{cor}

\begin{rem} {\rm
Because of Corollary~\ref{cor_1-intertwiners_modul_isomorfisme}, one may be tempted to consider only gauge trivial 1-intertwiners. It is worth being careful at this point, however, because such 1-intertwiners do not define a locally full sub-2-category of ${\repg}_{{\bf 2Mat}_{\Complex}}(\G)$. Thus, composing gauge trivial 1-intertwiners may give rise to non gauge trivial ones, as it is shown below, Proposition~\ref{composicio_1-intertwiners} (basically this is because the same thing already occurs when composing 1-morphisms in ${\bf 2Mat}_{\Complex}$; cf. Equation~(\ref{seccions_functor_composicio})). }
\end{rem}
The pair $({\bf R},{\sf S})$ can be interpreted geometrically as follows \footnote{The idea of thinking of 1-intertwiners as some sort of vector bundles already appears in the work by Crane and Yetter \cite{CY03} on the representation theory of automorphic 2-groups on Yetter's 2-category of measurable categories. Similar ideas can also be found in Mackaay and Barret \cite{BM04}.}. Let us think of the rank matrix {\bf R} as the (nonhomogeneous) complex vector bundle $p:E({\bf R})\To M(n',n)$ whose fiber over the point $(i',i)\in M(n',n)$ is $E({\bf R})_{i',i}=\Complex^{R_{i'i}}$. It follows from conditions (\ref{invariancia_R}) and (\ref{condicio_sobre_suport}) that the involved vector bundles are not arbitrary. Instead, they are entirely supported on the intertwining orbits of $M(n',n)$ and homogeneous over each such orbit, althought of different dimensions for different orbits in general. Then, the collection {\sf S} is equivalent to a projective right action of $\pi_0(\G)$ on $E({\bf R})$ making $p$ a $\pi_0(\G)$-equivariant map. More explicitly, given {\sf S}, define $\Theta({\sf S}):E({\bf R})\times\pi_0(\G)\To E({\bf R})$ by
\begin{equation} \label{definicio_accio_projectiva}
\Theta({\sf S})(x,g)={\sf S}_{i',i}(g)^{-1}(x)\in E({\bf R})_{\rho'(g)^{-1}(i'),\rho(g)^{-1}(i)}, \quad  x\in E({\bf R})_{i',i}
\end{equation}
for all $(i',i)\in{\rm Sup}({\bf R})$ (this makes sense because of (\ref{invariancia_R})). Then, it is immediate that (\ref{condicio_normalitzacio}) translates into the condition
\begin{equation} \label{accio_projectiva_1}
\Theta({\sf S})(x,e)=x
\end{equation}
for all $x\in E({\bf R})$ while (\ref{condicio_morfisme_sobre_T}), equivalent to
$$
{\sf S}_{i',i}(g_1g_2)^{-1}=\frac{c'(g_1,g_2)_{i'}}{c(g_1,g_2)_{i}}\ {\sf
S}_{\rho'(g_1)^{-1}(i'),\rho(g_1)^{-1}(i)}(g_2)^{-1}\ {\sf S}_{i',i}(g_1)^{-1}
$$
is nothing but the condition
\begin{equation} \label{accio_projectiva_2}
\Theta({\sf S})(x,g_1g_2)=\frac{c'(g_1,g_2)_{i'}}{c(g_1,g_2)_{i}}\ \Theta({\sf S})(\Theta({\sf S})(x,g_1),g_2)
\end{equation}
for all $x\in E({\bf R})_{i',i}$ and all $(i',i)\in{\rm Sup}({\bf R})$. Hence, $\Theta({\sf S})$ indeed defines a projective right action covering the action of $\pi_0(\G)$ on $M(n',n)$. Conversely, given any map $\Theta:E({\bf R})\times\pi_0(\G)\To E({\bf R})$ such that $\Theta(-,g)$ maps linearly the fiber $E({\bf R})_{i',i}$ onto the fiber $E({\bf R})_{\rho'(g)^{-1}(i'),\rho(g)^{-1}(i)}$ and satisfying (\ref{accio_projectiva_1}) and (\ref{accio_projectiva_2}), the reader may easily check that a collection ${\sf S}(\Theta)$ as above is obtained by taking
$$
{\sf S}(\Theta)_{i',i}(g)=\left[\Theta(-,g)_{|E({\bf R})_{i',i}}\right]^{-1}
$$
where $\left[\Theta(-,g)_{|E({\bf R})_{i',i}}\right]$ denotes the matrix in canonical basis of the restriction map $\Theta(-,g):E({\bf R})_{i',i}\To E({\bf R})_{\rho'(g)^{-1}(i'),\rho(g)^{-1}(i)}$. This leads to the following more suggestive description (up to isomorphism) of the non zero 1-intertwiners (cf. \cite{BM04}).

\begin{cor} \label{1-intertwiners_com_fibrats}
Let $(n,\rho,\beta,c)$ and $(n',\rho',\beta',c')$ be such that $I(\beta,\beta')\neq\emptyset$. Then, up to isomorphism, a non zero 1-intertwiner $\xi:(n,\rho,\beta,c)\To(n',\rho',\beta',c')$ is determined by a collection of (complex) vector bundles $\{p_{\Oo}:E_{\Oo}\To\Oo\}$, one for each intertwining $\pi_0(\G)$-orbit $\Oo$ of $M(n',n)$ and not all of them zero, together with projective right actions $\Theta_{\Oo}$ on $E_{\Oo}$ covering the action of $\pi_0(\G)$ on $\Oo$, with associated respective cocycles $z_{\Oo}$.  
\end{cor}

\begin{ex} {\rm
For any cocycles $z\in Z^2(\pi_0(\G),(\Complex^*)^n)$ and $z'\in Z^2(\pi_0(\G),(\Complex^*)^{n'})$, a 1-intertwiner $\R_{n,z}\To\R_{n',z'}$ is determined, up to isomorphism, by a collection of projective (anti)representation of $\pi_0(\G)$, one for each point $(i',i)\in M(n',n)$, with respective cocycles $z_{i',i}=z'_{i'}/z_i$. In particular, endomorphisms of the trivial representation $\I$ basically correspond to usual linear representations of $\pi_0(\G)$.
 }
\end{ex}

\begin{rem} {\rm 
More generally, for any non zero representations $(n,\rho,\beta,c), (n',\rho',\beta',c')$ with $I(\beta,\beta')\neq\emptyset$, implicit in the weakenings maps ${\sf S}_{i',i}$ of any non zero 1-intertwiner $({\bf R},s,{\sf S})$ or equivalently, in the projective right actions $\Theta_{\Oo}$, there is a projective (anti)representation of the stabilizer $\pi_0(\G)_{(i',i)}\subseteq \pi_0(\G)$ of the point $(i',i)$, for each point $(i',i)\in {\rm Sup} {\bf R}$, with cocycle $z_{i',i}\in Z^2(\pi_0(\G)_{i',i},\Complex^*)$ defined by $z_{i',i}(g_1,g_2)=z_{\Oo}(g_1,g_2)(i',i)$ if $(i',i)\in\Oo$. }
\end{rem} 

To finish our description of the categories of morphisms ${\bf Hom}((n,\rho,\beta,c),(n',\rho',\beta',c'))$ in case $I(\beta,\beta')\neq\emptyset$, it remains to determine what 2-intertwiners are and how they are composed. By Corollary~\ref{cor_1-intertwiners_modul_isomorfisme}, it would be enough to consider the case of two arbitrary gauge trivial 1-intertwiners, but the description turns out to be the same in all cases.

\begin{prop} \label{descripcio_2-intertwiners}
Let $\xi,\overline{\xi}:(n,\rho,\beta,c)\To(n',\rho',\beta',c')$ be two gauge trivial 1-intertwiners, respectively described by triples $({\bf R},s,{\sf S})$ and $(\overline{{\bf R}},\overline{s},\overline{{\sf S}})$. Then, there is a 1-1 correspondence between the set of 2-intertwiners $\ngg:\xi\Rightarrow\overline{\xi}$ and the set of $n'\times n$ matrices {\sf T} with entries ${\sf T}_{i',i}\in{\rm Mat}_{\overline{R}_{i',i}\times R_{i',i}}(\Complex)$ if $\overline{R}_{i',i},R_{i',i}\neq 0$ and empty otherwise and such that
\begin{equation}
{\sf T}_{i',i}\ {\sf S}_{i',i}(g)=\overline{{\sf S}}_{i',i}(g)\ {\sf T}_{\rho'(g)^{-1}(i'),\rho(g)^{-1}(i)}
\label{condicio_2-intertwiner}
\end{equation}
for all $(i',i)\in{\rm Sup}({\bf R})\cap{\rm Sup}(\overline{{\bf R}})$ and $g\in\pi_0(\G)$. Furthermore, under this correspondence, the (vertical) composite of two 2-intertwiners ${\sf T}:({\bf R},s,{\sf S})\Rightarrow(\overline{{\bf R}},\overline{s},\overline{{\sf S}})$ and $\overline{{\sf T}}:(\overline{{\bf R}},\overline{s},\overline{{\sf S}})\Rightarrow(\overline{\overline{{\bf R}}},\overline{\overline{s}},\overline{\overline{{\sf S}}})$ is given by the $n'\times n$ matrix $\overline{{\sf T}}\cdot{\sf T}$ with entries
\begin{equation} \label{matriu_2-intertwiner_composicio_1}
(\overline{{\sf T}}\cdot{\sf T})_{i',i}=\overline{{\sf T}}_{i',i}\ {\sf
T}_{i',i}
\end{equation}
for all $(i',i)\in M(n',n)$.
\end{prop}

\begin{proof}
As pointed out before, a 2-intertwiner $\ngg:\xi\Rightarrow\overline{\xi}$ is given by a 2-morphism 
$({\bf R},{\bf I})\Rightarrow(\overline{{\bf R}},{\bf I})$ in ${\bf 2Mat}_{\Complex}$, which is indeed given by a collection of matrices ${\sf T}_{i',i}\in{\rm Mat}_{\overline{R}_{i',i}\times R_{i',i}}(\Complex)$ for all $(i',i)\in{\rm Sup}({\bf R})\cap{\rm Sup}(\overline{{\bf R}})$. We need to see that the naturality axiom (\ref{axioma_naturalitat_modificacio}) translates into condition (\ref{condicio_2-intertwiner}). Let $\Phi(g)$ be defined as in (\ref{def_phi_vs_S}) and similarly $\overline{\Phi}(g)$. Then, (\ref{axioma_2-intertwiner}) reads in this case
$$
\overline{\Phi}(g)\cdot(1_{({\bf P}(\rho'(g)),{\bf I})}\circ{\sf T})=({\sf T}\circ 1_{({\bf P}(\rho(g)),{\bf I})})\cdot\Phi(g)
$$
Then, an easy computation similar to that made in proving Proposition~\ref{primera_descripcio_1-intertwiners} shows that this is indeed equivalent to (\ref{condicio_2-intertwiner}). As regards (\ref{matriu_2-intertwiner_composicio_1}), it readily follows from (\ref{composicio_vertical}).
\end{proof}

This has a nice interpretation in terms of the vector bundles and projective actions $\{(p:E_{\Oo}\To{\Oo},\Theta_{\Oo})\}_{\Oo}$ and $\{(\overline{p}:\overline{E}_{\Oo}\To\Oo,\overline{\Theta}_{\Oo})\}_{\Oo}$ describing each 1-intertwiner. Thus, for each point $(i',i)\in{\rm Sup}({\bf R})\cap{\rm Sup}(\overline{{\bf R}})$, belonging to some intertwining orbit $\Oo$, the matrix ${\sf T}_{i',i}$ corresponds to a linear map between the respective fibers $(E_{\Oo})_{i',i}$ and $(\overline{E}_{\Oo})_{i',i}$ and all of them together define vector bundle morphisms $f({\sf T})_{\Oo}:E_{\Oo}\To\overline{E}_{\Oo}$, one for each intertwining orbit $\Oo$ in ${\rm Sup}({\bf R})\cap{\rm Sup}(\overline{{\bf R}})$. Equation~(\ref{matriu_2-intertwiner_composicio_1}) clearly corresponds to composing such morphisms. Furthermore, if we denote both actions by $\cdot$ for short, it is (cf. (\ref{definicio_accio_projectiva})) 
$$
f({\sf T})_{\Oo}(x\cdot g)=f({\sf T})_{\Oo}({\sf S}_{i',i}(g)^{-1}(x))=({\sf T}_{\rho'(g)^{-1}(i'),\rho(g)^{-1}(i)}\ {\sf S}_{i',i}(g)^{-1})(x)
$$
and
$$
f({\sf T})_{\Oo}(x)\cdot g=\overline{{\sf S}}_{i',i}(g)^{-1}(f({\sf T})_{\Oo}(x))=(\overline{{\sf S}}_{i',i}(g)^{-1}\ {\sf T}_{i',i})(x)
$$
for all $x\in(E_{\Oo})_{i',i}$, all $(i',i)\in\Oo$, all intertwining orbits $\Oo$ in ${\rm Sup}({\bf R})\cap{\rm Sup}(\overline{{\bf R}})$ and all $g\in\pi_0(\G)$. Hence, (\ref{condicio_2-intertwiner}) is nothing but the condition that each morphism $f({\sf T})_{\Oo}$ preserves the corresponding projective right actions $\Theta_{\Oo}$ and $\overline{\Theta}_{\Oo}$.

Putting together all previous results, we finally get the following more geometric description for the categories of morphisms between two non zero 2-matrix representations:

\begin{thm} \label{segona_descripcio_1-intertwiners}
For any non zero 2-matrix representations $(n,\rho,\beta,c)$ and $(n',\rho',\beta',c')$, with $n,n'\geq 1$, let ${\rm Orb}_I(n,\rho,\beta;n',\rho',\beta')$ be the associated set of intertwining $\pi_0(\G)$-orbits of $M(n',n)$ and $z_{\Oo}$ the induced 2-cocycles defined above (cf. Equation~(\ref{definicio_cocicles})). Then, there is an equivalence of categories
\begin{equation} \label{categories_morfismes}
{\bf Hom}((n,\rho,\beta,c),(n',\rho',\beta',c'))\simeq \prod_{\Oo\in{\rm Orb}_I(n,\rho,\beta;n',\rho',\beta')}\ {\bf PBund}_{\pi_0(\G),z_{\Oo}}(\Oo),
\end{equation}
where ${\bf PBund}_{\pi_0(\G),z_{\Oo}}(\Oo)$ denotes the category with objects the (homogeneous) vector bundles over $\Oo$ equipped with a projective (right) action of $\pi_0(\G)$ covering its action on $\Oo$ and of cocycle $z_{\Oo}$ and morphisms the vector bundle morphisms preserving these actions. In particular, ${\bf Hom}((n,\rho,\beta,c),(n',\rho',\beta',c'))\simeq{\bf 1}$ when ${\rm Orb}_I(n,\rho,\beta;n',\rho',\beta')=\emptyset$.
\end{thm}

Note that the categories of morphisms depend on the 2-cochains $c,c'$ coding the ``weak'' character of the involved representations uniquely via the induced 2-cocycles $z_{\Oo}$ giving the projective character of the actions.

\subsection{Composition functors}

For later use, let us make explicit also the composition functors in ${\repg}_{{\bf 2Mat}_{\Complex}}(\G)$. This means making explicit the composition law between 1-intertwiners and the horizontal composition law between 2-intertwiners. They respectively correspond to 
vertically composing pseudonatural transformations, according to
Equations~(\ref{composicio_vertical_transformacions_pseudo_1}) and (\ref{composicio_vertical_transformacions_pseudo_2}), and to horizontally composing modifications, according to Equation~(\ref{composicio_modificacions_2}). Since it seems more involved to describe such laws in the above geometric picture of vector bundles (over the intertwining orbits) and morphisms between these, we shall content ourselves with the description in terms of the data appearing in Propositions~\ref{primera_descripcio_1-intertwiners} and \ref{descripcio_2-intertwiners}.

\begin{prop} \label{composicio_1-intertwiners}
Let $\xi=({\bf R},s,{\sf S}):(n,\rho,\beta,c)\To(n',\rho',\beta',c')$ and $\xi'=({\bf R}',s',{\sf S}'):(n',\rho',\beta',c')\To(n'',\rho'',\beta'',c'')$ be two composable 1-intertwiners. Then, the rank matrix $\hat{{\bf R}}$, gauge $\hat{s}$ and weakening maps $\hat{{\sf S}}_{i'',i}$, for all $(i'',i)\in M(n'',n)$, describing the composite 1-intertwiner $\xi'\circ\xi$ are given by
\begin{equation} 
\hat{{\bf R}}={\bf R}'\ {\bf R}
\label{matriu_rangs_1-intertwiner_composicio}
\end{equation}
\begin{equation} \label{gauge_1-intertwiner_composicio}
\hat{s}_{i''}({\bf a})=s'_{i''}({\bf R}({\bf
a}))
\left(\bigoplus_{i'=1}^{n'}\ {\bf I}_{R'_{i''i'}}\otimes
s_{i'}({\bf a})\right)\ {\bf P}({\bf
R}'_{i''},{\bf R},{\bf a})\left(\bigoplus_{i=1}^n\
s'_{i''}({\bf R}({\bf e}_i))^{-1}\otimes {\bf
I}_{a_i}\right)
\end{equation}
\begin{align}
\hat{{\sf S}}_{i''i}(g)=s'_{i''}({\bf R}({\bf
e}_{i}))&\left(\bigoplus_{i'=1}^{n'}\ {\bf
I}_{R'_{i'',i'}}\otimes{\sf S}_{i',i}(g)\right)s'_{i''}({\bf R}({\bf
e}_{i}))^{-1} \nonumber
\\ &
s'_{i''}({\bf R}({\bf e}_{\rho(g)^{-1}(i)}))\left(\bigoplus_{i'=1}^{n'}\ {\sf S}'_{i'',\rho'(g)(i')}(g)\otimes {\bf I}_{R_{\rho'(g)(i'),i}}\right)s'_{i''}({\bf R}({\bf e}_{\rho(g)^{-1}(i)}))^{-1}
\label{weakening_map_1-intertwiner_composicio}
\end{align}
\end{prop}
\begin{proof}
Equations~(\ref{matriu_rangs_1-intertwiner_composicio}) and (\ref{gauge_1-intertwiner_composicio}) are an immediate consequence of (\ref{composicio_vertical_transformacions_pseudo_1}) and (\ref{composicio_1-morfismes})-(\ref{seccions_functor_composicio}). To prove (\ref{weakening_map_1-intertwiner_composicio}), let $\Phi(g):({\bf P}(\rho'(g)){\bf R},s)\Rightarrow({\bf RP}(\rho(g)),s)$ and $\Phi'(g):({\bf P}(\rho''(g)){\bf R}',s')\Rightarrow({\bf R}'{\bf P}(\rho'(g)),s')$ be the 2-isomorphisms in ${\bf 2Mat}_{\Complex}$ respectively associated to $\xi$ and $\xi'$, related to the weakening maps ${\sf S}_{i',i}$ and ${\sf S}'_{i'',i'}$ as in Equation~(\ref{def_phi_vs_S}). Then, by Equation~(\ref{composicio_vertical_transformacions_pseudo_2}), the 2-isomorphism $\hat{\Phi}(g):({\bf P}(\rho''(g))\hat{{\bf R}},\hat{s})\Rightarrow(\hat{{\bf R}}{\bf P}(\rho(g)),\hat{s})$ associated to the composite $\xi'\circ\xi$ is given by
$$
\hat{\Phi}(g)=(1_{({\bf R}',s')}\circ\Phi(g))\cdot(\Phi'(g)\circ 1_{({\bf
R},s)}),\quad g\in G
$$
By definition, we have $\hat{{\sf S}}_{i'',i}(g)=\hat{\Phi}(g)_{i'',\rho(g)^{-1}(i)}$, for all $g\in\pi_0(\G)$. Hence, Equation~(\ref{composicio_vertical}) gives that
$$
\hat{{\sf S}}_{i'',i}(g)=(1_{({\bf R}',s')}\circ\Phi(g))_{i'',\rho(g)^{-1}(i)}\ (\Phi'(g)\circ 1_{({\bf R},s)})_{i'',\rho(g)^{-1}(i)}
$$
and an easy computation using Equation~(\ref{composicio_horitzontal}) gives for the matrices in the right hand side the following expressions: 
\begin{align*}
(1_{({\bf R}',s')}\circ\Phi(g))_{i'',\rho(g)^{-1}(i)}&=s'_{i''}({\bf R}({\bf e}_{i}))\left(\bigoplus_{i'=1}^{n'}\ {\bf I}_{R'_{i'',i'}}\otimes{\sf S}_{i',i}(g)\right)s'_{i''}({\bf R}({\bf e}_{i}))^{-1}
\\ (\Phi'(g)\circ 1_{({\bf R},s)})_{i'',\rho(g)^{-1}(i)}&=s'_{i''}({\bf R}({\bf e}_{\rho(g)^{-1}(i)}))\left(\bigoplus_{i'=1}^{n'}\ {\sf S}'_{i'',\rho'(g)(i')}(g)\otimes {\bf I}_{R_{\rho'(g)(i'),i}}\right)s'_{i''}({\bf R}({\bf e}_{\rho(g)^{-1}(i)}))^{-1}
\end{align*}
Putting this into the previous equation gives (\ref{weakening_map_1-intertwiner_composicio}).
\end{proof}

Note that (\ref{gauge_1-intertwiner_composicio}) and (\ref{weakening_map_1-intertwiner_composicio}) in particular imply that the composite of two gauge trivial 1-intertwiners is non gauge trivial in general, because
of the permutation matrix ${\bf P}({\bf
R}^{\xi'}_{i''},{\bf R}^{\xi},{\bf a})$, and that the composite of two
strict 1-intertwiners (gauge trivial or not) is always strict as required.
Observe also that the previous equations together with the normalization conditions on the gauge and permutation matrices ${\bf P}({\bf R}^{\xi'}_k,{\bf R}^{\xi},{\bf a})$ (see
Section~\ref{2cat_2vect}) show that the triple $({\bf I}_n,{\bf I},{\bf I})$ acts as a unit for this composition. Hence, for any non zero 2-matrix representation $(n,\rho,\beta,c)$, its identity 1-arrow is given by the line bundle over the diagonal of $M(n,n)$ equipped with the trivial linear action of $\pi_0(\G)$ covering the $\pi_0(\G)$-action on the diagonal.  

\begin{prop} \label{composicio_horitzontal_2-intertwiners}
Let ${\sf T}:({\bf R},s,{\sf S})\Rightarrow(\overline{{\bf R}},\overline{s},\overline{{\sf S}}):(n,\rho,\beta,c)\To(n',\rho',\beta',c')$ and
${\sf T}':({\bf R}',s',{\sf S}')\Rightarrow(\overline{{\bf R}}',\overline{s}',\overline{{\sf S}}'):(n',\rho',\beta',c')\To(n'',\rho'',\beta'',c'')$ be any two horizontally composable 2-intertwiners. Then, the horizontal composite ${\sf T}'\circ{\sf T}$ is given by the $n''\times n$ matrix with components
\begin{equation} \label{matriu_2-intertwiner_composicio_2}
({\sf T}'\circ{\sf T})_{i'',i}=\overline{s}'_{i''}(\overline{{\bf
R}}({\bf e}_i))\ \left(\bigoplus_{i'=1}^{n'}\ {\sf
T}'_{i'',i'}\otimes{\sf T}_{i',i}\right)\
s'_{i''}({\bf R}({\bf e}_i))^{-1}
\end{equation}
\end{prop}
\begin{proof}
The proof is a computation similar to those made before and it is left to the reader (start with (\ref{composicio_modificacions_2}) and use (\ref{composicio_horitzontal})).
\end{proof}

\subsection{The set $\pi_0({\repg}_{{\bf 
2Mat}_{\Complex}}(\G))$}

\label{classes_equivalencia_repr}

An important (biequivalence) invariant of any 2-category $\Cgg$ is the set $\pi_0(\Cgg)$ of equivalence classes of objects in $\Cgg$. We wish now to compute this invariant for $\Cgg={\repg}_{{\bf 2Mat}_{\Complex}}(\G)$.

In general, $\pi_0(\Cgg)$ is a certain quotient of the set of isomorphism classes of objects in $\Cgg$. As the following result shows, however, both sets coincide in case $\Cgg={\repg}_{{\bf 
2Mat}_{\Complex}}(\G)$.

\begin{prop}
The non zero 2-matrix representations defined by the quadruples $(n,\rho,\beta,c)$ and $(n',\rho',\beta',c')$ are equivalent if and only if they are isomorphic. 
\end{prop}
\begin{proof}
The result is a consequence of the fact that ${\repg}_{{\bf
2Mat}_{\Complex}}(\G)$ is a pseudofunctor
2-category and of Lemmas \ref{lema_2-categoria_2-functors} and
\ref{morfismes_invertibles}. More explicitly,
suppose we have a 1-intertwiner $\xi=({\bf R},s,{\sf S}):(n,\rho,\beta,c)\To(n',\rho',\beta',c')$
which is invertible up to a 2-isomorphism. Let $\xi'=({\bf R}',s',{\sf
S}'):(n',\rho',\beta',c')\To(n,\rho,\beta,c)$ be a quasiinverse. In particular, $\xi'\circ\xi$ and
$\xi\circ\xi'$ are 2-isomorphic to ${\id}_{(n,\rho,\beta,c)}$ and 
${\id}_{(n',\rho',\beta',c')}$, respectively. Then, according to Lemma
\ref{lema_2-categoria_2-functors}, $({\bf
R}',s')\circ({\bf R},s)=({\bf R}'{\bf R},s'\ast s)$ and $({\bf
R},s)\circ({\bf 
R}',s')=({\bf R}{\bf R}',s\ast s')$ are 2-isomorphic in ${\bf
2Mat}_{\Complex}$ to
${\id}_n=({\bf I}_n,{\bf I})$ and ${\id}_{n'}=({\bf I}_{n'},{\bf I})$,
respectively. But Lemma~\ref{morfismes_invertibles}-(ii) says that this is
true if and only if ${\bf R}'{\bf 
  R}={\bf I}_n$ and ${\bf R}{\bf
  R}'={\bf I}_{n'}$, from which it follows that $n=n'$ and
that both {\bf R} and ${\bf R}'$ are permutation matrices. Using again Lemma~\ref{lema_2-categoria_2-functors} and Lemma~\ref{morfismes_invertibles}-(i), we conclude that $\xi$ is strictly invertible.
\end{proof}

We already have a description up to isomorphism of the objects of ${\repg}_{{\bf 2Mat}_{\Complex}}(\G)$, given by Theorem~\ref{representacions_modul_isomorfisme}. As pointed out repeteadly, however, this description is non faithful. To get the desired 1-1 correspondence with the set of isomorphism classes, we just need to take the following quotient on the set of quadruples $(n,\rho,\beta,c)$:

\begin{defn} \label{relacio_equivalencia_entre_quaternes}
Two quadruples $(n,\rho,\beta,c)$ and $(n',\rho',\beta',c')$ as in
Theorem~\ref{representacions_modul_isomorfisme} are called {\sl
  equivalent} if the following conditions hold: 
\begin{itemize}
\item[(i)]
$n=n'$, and
\item[(ii)]
there exists a permutation $\sigma\in S_n$ such that:
\begin{itemize}
\item[(a)]
$\rho'(g)=\sigma\rho(g)\sigma^{-1}$ for all $g\in\pi_0(\G)$.
\item[(b)]
$\beta'(u)=\sigma\cdot\beta(u)$ for all $u\in\pi_1(\G)$.
\item[(c)]
$[c']=[\sigma\cdot c]$ in
  $\widetilde{H}^2(\pi_0(\G),(\Complex^*)^n_{\rho'})$, where
  $\sigma\cdot c\in C^2(\pi_0(\G),(\Complex^*)^n_{\rho'})$ is the
  normalized 2-cochain defined by $(\sigma\cdot
  c)(g_1,g_2)=\sigma\cdot c(g_1,g_2)$ 
\end{itemize}
(the action of $S_n$ on $(\Complex^*)^n$ is that given by
(\ref{wreath_product})). 
\end{itemize} 
\end{defn}

\begin{thm} \label{representacions_modul_equivalencia}
There is a 1-1 correspondence between $\pi_0({\repg}_{{\bf 2Mat}_{\Complex}}(\G))$ and the set of equivalences classes of quadruples $(n,\rho,\beta,c)$ of Theorem~\ref{representacions_modul_isomorfisme}.
\end{thm} 
\begin{proof}
If $[n,\F]$ and $[n,\rho,\beta,c]$ respectively denote the equivalence classes of the $n$-dimensional 2-matrix representation $\F:\G\To{\sf Equiv}_{{\bf 2Mat}_{\Complex}}(n)$ and of the quadruple $(n,\rho,\beta,c)$, the correspondence is that given by $[n,\rho,\beta,c]\mapsto[n,\F(\rho,\beta,c)]$, with $\F(\rho,\beta,c)$ defined by (\ref{accio_representacio_sobre_objectes})-(\ref{isos_estructurals_representacio}). This is well defined. Thus, suppose $(n,\rho,\beta,c)$ and $(n,\rho',\beta,c')$ are equivalent via the permutation $\sigma\in S_n$, and let $x\in C^1(\pi_0(\G),(\Complex^*)^n_{\rho'})$ be a normalized 1-cochain such that $c'=\sigma\cdot c+\partial x$. Then, we clearly have $(i,\sigma^{-1}(i)\in I(\beta,\beta')$ for all $i=1,\ldots,n$. In particular, $I(\beta,\beta')\neq\emptyset$. Set ${\bf R}={\bf P}(\sigma)$ (hence, ${\rm Sup}({\bf R})=\{(i,\sigma^{-1}(i),\ i=1,\ldots,n\}$) and
$$
{\sf S}_{i,\sigma^{-1}(i)}(g)=x_i(g),\quad g\in\pi_0(\G),\ i=1,\ldots,n,
$$
The reader may easily check that {\bf R} is a $\left(\begin{array}{c} n\ \rho\ \beta \\ n\ \rho'\beta'\end{array}\right)$-admissible rank matrix and that the maps ${\sf S}_{i,\sigma^{-1}(i)}:\pi_0(\G)\To\Complex^*$ satisfy (\ref{condicio_normalitzacio}) and (\ref{condicio_morfisme_sobre_T}). Hence, by Proposition~\ref{primera_descripcio_1-intertwiners}, we have a 1-intertwiner $({\bf R},{\bf I},{\sf S}):(n,\F(\rho,\beta,c))\To(n,\F(\rho',\beta',c'))$, which moreover is a 1-isomorphism because of Lemmas~\ref{lema_2-categoria_2-functors} and \ref{morfismes_invertibles}-(i). Consequently, $[n,\F(\rho,\beta,c)]=[n,\F(\rho',\beta',c')]$. It clearly follows from Theorem~\ref{representacions_modul_isomorfisme} that the correspondence so defined is surjective. To prove that it is injective, it suffices to go back in the previous argument to see that $[n,\F(\rho,\beta,c)]=[n',\F(\rho',\beta',c')]$ implies $[n,\rho,\beta,c]=[n',\rho',\beta',c']$.
\end{proof}

\begin{ex} {\rm 
Let $\G$ be the dihedral group $D_{2m}$ ($m\geq 2$) thought of as a split 2-group as explained in Example~\ref{exemples_2-grups_escindits}-(2). By choosing as classifying 3-cocycle the trivial one, we get a 1-1 correspondence between equivalence classes of non zero 2-matrix representations of $D_{2m}$ and equivalence classes of quadruples $(n,\rho,\beta,z)$, with $n\geq 1$, $\rho:\Integer_2\To S_n$ a group morphism, $\beta:\Integer_m\To(\Complex^*)^n_{\rho}$ a map of $\Integer_2$-modules and $z\in Z^2(\Integer_2,(\Complex^*)^n_{\rho})$. In particular, the set of equivalence classes of 1-dimensional 2-matrix representations can be identified with $\Complex^*$ for $m=2$ or odd, and with $\Complex^*\cup\Complex^*$ for $m\neq 2$ even. Indeed, it is easily checked that for $m=2$ or odd, there is a unique map $\beta:\Integer_m\To\Complex^*$ of $\Integer_2$-modules (namely, the trivial one $\beta_1(\overline{k})=1$), while for $m\neq 2$ even we also have the map $\beta_2(\overline{k})=(-1)^k$. Furthermore, $H^2(\Integer_2,\Complex^*)\cong\Complex^*$. The result follows then from the fact that two pairs $(\beta,z)$ and $(\beta',z')$ are equivalent if and only if $\beta=\beta'$ and $[z]=[z']$ (cf. Definition~\ref{relacio_equivalencia_entre_quaternes}). 
}
\end{ex}

\subsection{2-category of characters of $\G$}

For any group $G$, a (linear) character of $G$, defined as a group morphism $\chi:G\To\Complex^*$, is exactly the same thing as a 1-dimensional matrix representation. This naturally suggests the following definition of character for 2-groups:

\begin{defn}
A {\sl character} of a 2-group $\G$ is a 1-dimensional 2-matrix representation of $\G$. 
\end{defn}
Thus, according to Example~\ref{exemples_representacions}-(1), any character of $\G$ is isomorphic to a representation $\I_{\beta,c}$ (or $\I_{\beta,z}$ if $\G$ is split) for some pair $(\beta,c)$  where $\beta:\pi_1(\G)\To\Complex^*$ is any morphism of $\pi_0(\G)$-modules such that $[\beta\circ\alpha]=0$ and $c\in C^2(\pi_0(\G),\Complex^*)$ a normalized 2-cochain such that $\partial c=\beta\circ\alpha$ (if $\G$ is split, $[\beta\circ\alpha]=0$ always hold and $c$ can be chosen to be a 2-cocycle $z\in Z^2(\pi_0(\G),\Complex^*)$). Furthermore, by Theorem~\ref{representacions_modul_equivalencia}, two such characters $\I_{\beta,c}$ and $\I_{\beta',c'}$ are isomorphic if and only if $\beta=\beta'$ and $[c]=[c']$ in $\widetilde{H}^2(\pi_0(\G),\Complex^*)$. 

In the group setting, characters define a full subcategory ${\bf Char}(G)$ of ${\bf Rep}_{{\bf Mat}_{\Complex}}(G)$, the category of matrix representations of $G$. Explicitly, it is easily seen that for any objects $\chi,\chi'$ of ${\bf Char}(G)$, it is
$$
{\rm Hom}(\chi,\chi')=\left\{ \begin{array}{ll} 0, & \mbox{if $\chi\neq\chi'$} \\ \Complex, & \mbox{if $\chi=\chi'$} \end{array} \right.
$$
with composition given by the product in $\Complex$.

Similarly, in the framework of 2-groups, characters define a full sub-2-category ${\sf Char}(\G)$ of $\repg_{{\bf 2Mat}_{\Complex}}(\G)$. Interestingly enough, the category of (projective) matrix representations of $\pi_0(\G)$ appears as category of morphisms between some pairs of objects in this 2-category. More explicitly, we have the following:

\begin{thm} \label{2-categoria_caracters}
For any $\G$ and pairs $(\beta,c)$ and $(\beta',c')$ as above, there is an equivalence of categories
\begin{equation} \label{homs_2-categoria_caracters}
{\bf Hom}(\I_{\beta,c},\I_{\beta',c'})\simeq\left\{ \begin{array}{ll}
{\bf 1}, & \mbox{if $\beta\neq\beta'$} \\ {\bf
  PRep}_{[c'-c]}(\pi_0(\G)), & \mbox{if $\beta=\beta'$} \end{array}
\right. 
\end{equation}
where ${\bf PRep}_{[c'-c]}(\pi_0(\G))$ denotes the category of
projective representations of $\pi_0(\G)$ with cohomology class
$[c'-c]\in H^2(\pi_0(\G),\Complex^*)$ (note that $c'-c$ is a 2-cocycle
if and only if $\beta=\beta'$). Furthermore, when
$(\beta,c)=(\beta',c')$, it is  
$$
{\bf End}(\I_{\beta,c})\simeq {\bf Rep}_{{\bf Mat}_{\Complex}}(\pi_0(\G))
$$
and this is an equivalence as monoidal categories (in particular, this
is true for the trivial representation $\I$). 
\end{thm}

\begin{proof}
The first statement is an immediate consequence of Theorem~\ref{segona_descripcio_1-intertwiners}. Thus, $M(1,1)$ is a singleton and the unique $\pi_0(\G)$-orbit is intertwining iff $\beta=\beta'$. Moreover, when $|\Oo|=1$, ${\bf PBund}_{\pi_0(\G),z_{\Oo}}(\Oo)$ clearly reduces to the category ${\bf PRep}_{z_{\Oo}}(\pi_0(\G))$, with $z_{\Oo}$ given by $z_{\Oo}=c'-c$ (see (\ref{definicio_cocicles})).

To prove the last assertion, observe that the equivalence ${\bf End}(\I_{\beta,c})\simeq {\bf Rep}_{{\bf Mat}_{\Complex}}(\pi_0(\G))$ follows from (\ref{homs_2-categoria_caracters}) and the fact that projective representations with trivial cohomology class are equivalent to usual matrix representations (see Kirillov \cite{aK76}, \S~2.3). Otherwise, it can also be deduced directly from Propositions~\ref{primera_descripcio_1-intertwiners} and \ref{descripcio_2-intertwiners}. Thus, for $n=n'=1$ and $c=c'$, the rank matrix of an object of ${\bf End}(\I_{\beta,c})$ reduces to a natural number $r$ and the (unique) weakening map to a group morphism ${\sf S}:\pi_0(\G)\To{\rm GL}(r)$, while a morphism between two such objects $(r,s,{\sf S})$, $(r',s',{\sf S}')$ reduces to an $r'\times r$ matrix {\sf T} such that ${\sf T}{\sf S}(g)={\sf S}'(g){\sf T}$ for all $g\in\pi_0(\G)$. It remains to see that the equivalence is as monoidal categories, the monoidal structure on ${\bf End}(\I_{\beta,c})$ being that defined by the corresponding composition functor (hence, by Propositions~\ref{composicio_1-intertwiners} and \ref{composicio_horitzontal_2-intertwiners}). Let us denote by $F$ this equivalence. In particular, it is $F(r,s,{\sf S})={\sf S}$ and $F({\sf T})={\sf T}$. The unit object in ${\bf End}(\I_{\beta,c})$ is ${\id}_{\I_{\beta,c}}$ and hence, it is mapped to the trivial 1-dimensional matrix representation, which is the unit object of ${\bf Rep}_{{\bf Mat}_{\Complex}}(\pi_0(\G))$. This allows us to take $F_0$ equal to the identity. By (\ref{matriu_rangs_1-intertwiner_composicio})-(\ref{weakening_map_1-intertwiner_composicio}), the tensor product $(r',s',{\sf S}')\otimes(r,s,{\sf S})=(r',s',{\sf S}')\circ(r,s,{\sf S})$ is the triple $(r'r,\hat{s},\hat{{\sf S}})$, with
$$
\hat{s}(a)=s'(ra)\ ({\bf I}_r\otimes s(a))\ (s'(r)^{-1}\otimes{\bf I}_a),\quad a\in\Natural
$$
(the permutation matrix in (\ref{gauge_1-intertwiner_composicio}) is here trivial) and
\begin{equation}
\hat{{\sf S}}(g)=s'(r)\ ({\bf I}_{r'}\otimes{\sf S}(g))\ ({\sf S}'(g)\otimes{\bf I}_r)\ s'(r)^{-1}=s'(r)\ ({\sf S}'(g)\otimes{\sf S}(g))\ s'(r)^{-1} \label{S_producte_tensorial}
\end{equation}
Therefore, we need to prove the existence of natural isomorphisms $F_2((r',s',{\sf S}'),(r,s,{\sf S})):{\sf S}'\otimes{\sf S}\To\hat{{\sf S}}$ satisfying (\ref{axioma_hexagon}) and (\ref{axiomes_unitat}) with $F_0$ an identity. Take
$$
F_2((r',s',{\sf S}'),(r,s,{\sf S}))=s'(r)\in{\rm GL}(r'r)
$$
By (\ref{S_producte_tensorial}), this indeed defines an isomorphism between both representations ${\sf S}'\otimes{\sf S}$ and $\hat{{\sf S}}$ and it readily follows from (\ref{matriu_2-intertwiner_composicio_2}) that these isomorphisms are natural in $((r',s',{\sf S}'),(r,s,{\sf S}))$. Furthermore, since both ${\bf End}(\I_{\beta,c})$ and ${\bf Rep}_{{\bf Mat}_{\Complex}}(\pi_0(\G))$ are strict as monoidal categories, conditions (\ref{axiomes_unitat}) amount to the triviality of $F_2((r',s',{\sf S}'),{\id}_{\I_{\beta,c}})$ and $F_2({\id}_{\I_{\beta,c}},(r',s',{\sf S}'))$. Now, the first morphism is trivial because of the normalization condition on the gauge (it is $r=1$ in this case) and the second one because the gauge $s'$ itself is trivial. Finally, taking again into account that both involved monoidal categories are strict, the reader may easily check that (\ref{axioma_hexagon}) is equivalent to
$$
s''(r'r)\ ({\bf I}_{r''}\otimes s'(r))\ (s''(r')^{-1}\otimes{\bf I}_r)\ (s''(r')\otimes {\bf I}_r)=s''(r'r)\ ({\bf I}_{r''}\otimes s'(r))
$$
which is trivially satisfied whatever they may be the gauges $s,s',s''$.  
\end{proof}

Note that, by taking $\beta$ equal to the constant map $\beta(u)=1$, all categories of $[z]$-projective representations of $\pi_0(\G)$, for any class $[z]\in H^2(\pi_0(\G),\Complex^*)$, are obtained as categories of morphisms between suitable 1-dimensional 2-matrix representations of $\G$.

\subsection{${\repg}_{{\bf 2Mat}_{\Complex}}(\G)$ vs the representation theory of $\pi_0(\G)$ on finite sets}

It has been shown in the previous paragraph that 2-matrix
representation theory of a 2-group $\G$ contains usual matrix
representation theory of its first homotopy group $\pi_0(\G)$. The
purpose of this paragraph is to show that the representation theory of
$\pi_0(\G)$ on finite sets is also implicit in the representation
theory of $\G$ on ${\bf 2Mat}_{\Complex}$. 

To make explicit this assertion, let us recall that, for any
2-category $\Cgg$, one defines its {\sl homotopy category}, which
we shall denote by ${\sf Ho}(\Cgg)$, as the category with the same
objects as $\Cgg$ and with sets of morphisms between any two objects
equal to the sets of isomorphism classes of 1-morphisms between these
objects (composition is given by $[g]\circ
[f]=[g\circ f]$ for $[f],[g]$ any composable arrows) \footnote{Unlike
the {\it underlying category} of $\Cgg$, i.e., the category having
the same objects as $\Cgg$ and 1-morphisms of $\Cgg$ as morphisms, the
homotopy category ${\sf
    Ho}(\Cgg)$ is a biequivalence invariant of $\Cgg$.}. For instance, for
locally discrete 2-categories it is 
${\sf Ho}(\Cc[0])=\Cc$. Furthermore, any 2-functor $\Ff:\Dgg\To\Cgg$
between 2-categories induces a functor ${\sf Ho}(\Ff):{\sf
  Ho}(\Dgg)\To{\sf Ho}(\Cgg)$ between the corresponding homotopy
categories, acting on morphisms by ${\sf Ho}(Ff)[f]=[\Ff(f)]$ (see
\cite{jB67}, \S~7). Then, if ${\bf FinSets}$ stands for the
category of finite sets, we have the following: 

\begin{thm}
For any 2-group $\G$, there is a 2-functor $\Ff:{\bf Rep}_{{\bf
FinSets}}(\pi_0(G))[0]\To{\repg}_{{\bf 2Mat}_{\Complex}}(\G)$ such
that the induced functor ${\sf Ho}(\Ff):{\bf Rep}_{{\bf
FinSets}}(\pi_0(G))\To {\sf Ho}({\repg}_{{\bf
2Mat}_{\Complex}}(\G))$ is faithful and essentially injective on
objects (i.e., nonisomorphic objects are mapped to nonisomorphic objects). 
\end{thm}
\begin{proof}
By choosing a total order in each object of {\bf FinSets}, objects in
${\bf Rep}_{{\bf FinSets}}(\pi_0(\G))[0]$ can be identified with pairs
$(n,\rho)$, for some $n\geq 1$ and some group morphism
$\rho:\pi_0(\G)\To S_n$, and 1-morphisms between two such pairs
$(n,\rho)$ and $(n',\rho')$ with $\pi_0(\G)$-equivariant maps
$f:\{1,\ldots,n\}\To\{1,\ldots,n'\}$, with $\pi_0(\G)$ acting on
$\{1,\ldots,n\}$ and $\{1,\ldots,n'\}$ via $\rho$ and $\rho'$,
respectively. Then, define $\Ff$ as follows: 
\begin{itemize}
\item on objects: $\Ff(n,\rho)=\R_{n,\rho}$, the 2-matrix
  representation defined in
  Example~\ref{exemples_representacions}-(2). 
\item on 1-morphisms: if $f:(n,\rho)\To(n',\rho')$, $\Ff(f)=({\bf
  R}(f),{\bf I},{\bf I})$, with ${\bf R}(f)$ the $n'\times n$ rank
  matrix given by 
$$
R(f)_{i',i}=\delta_{i',f(i)}
$$
for all $(i',i)\in M(n',n)$.
\end{itemize}
As the maps $\beta,\beta'$ and 2-cochains $c,c'$ associated to the
representations $\R_{n,\rho}$ and $\R_{n',\rho'}$ are all trivial,
(\ref{condicio_sobre_suport}) and
(\ref{condicio_normalitzacio})-(\ref{condicio_morfisme_sobre_T}) hold
trivially. Furthermore, it is easily checked that ${\bf R}(f)$ is
$(\rho,\rho')$-invariant due to the $G$-equivariance of $f$. Hence,
Proposition~\ref{primera_descripcio_1-intertwiners} ensures that the
triple $({\bf R}(f),{\bf I},{\bf I})$ indeed defines a 1-intertwiner
from $\R_{n,\rho}$ to $\R_{n',\rho'}$. Besides, if $f$ is an identity
map, ${\bf R}(f)$ is an identity matrix and $\Ff(f)$ is the
corresponding identity 1-intertwiner and, given another pair
$(n'',\rho'')$ and a $\pi_0(\G)$-equivariant map
$f':\{1,\ldots,n'\}\To\{1,\ldots,n''\}$, it follows from
Proposition~\ref{composicio_1-intertwiners} and the first
normalization condition on the permutation matrix which appears in
(\ref{gauge_1-intertwiner_composicio}) that 
$$
\Ff(f')\circ\Ff(f)=({\bf R}(f'),{\bf I},{\bf I})\circ({\bf R}(f),{\bf I},{\bf I})=({\bf R}(f'){\bf R}(f),{\bf I},{\bf I})=({\bf R}(f'\circ f),{\bf I},{\bf I})
$$
Hence, we also have $\Ff(f')\circ\Ff(f)=\Ff(f'\circ f)$. This proves
that $\Ff$ is a 2-functor. To show the essential injectivity of ${\sf
  Ho}(\Ff)$, suppose $\R_{n,\rho}$ and $\R_{n',\rho'}$ are isomorphic
in ${\sf Ho}({\repg}_{{\bf 2Mat}_{\Complex}}(\G))$. This means that
they are equivalent objects of ${\repg}_{{\bf 2Mat}_{\Complex}}(\G)$
and, by Theorem~\ref{representacions_modul_equivalencia}, that $n=n'$
and $\rho'=\sigma\rho\sigma^{-1}$ for some permutation $\sigma\in
S_n$. Hence, $(n,\rho)$ and $(n',\rho')$ are isomorphic in ${\bf
  Rep}_{{\bf FinSets}}(\pi_0(G))$. Finally, given morphisms
$f,\overline{f}:(n,\rho)\To(n',\rho')$, the equality ${\sf
Ho}(\Ff)(f)={\sf Ho}(\Ff)(\overline{f})$ means that both
1-intertwiners $({\bf R}(f),{\bf I},{\bf I}),({\bf
  R}(\overline{f}),{\bf I},{\bf I}):\R_{n,\rho}\To\R_{n',\rho'}$ are
isomorphic and hence, ${\bf R}(f)={\bf R}(\overline{f})$. But this
implies $f=\overline{f}$.   
\end{proof}

\begin{cor}
For any 2-group $\G$, ${\bf Rep}_{{\bf
FinSets}}(\pi_0(G))$ is equivalent to a subcategory of the homotopy
category of ${\repg}_{{\bf 2Mat}_{\Complex}}(\G)$.
\end{cor}

Note, however, that ${\sf Ho}(\Ff)$ is not full because all 1-intertwiners in
the image are line bundles, while the dimensions of the various vector
bundles describing a 1-intertwiner between $\R_{n,\rho}$ and
$\R_{n',\rho'}$ can be chosen arbitrarily. Hence, ${\bf Rep}_{{\bf
FinSets}}(\pi_0(G))$ is not equivalent to a full subcategory of ${\sf Ho}({\repg}_{{\bf 2Mat}_{\Complex}}(\G))$.

\vspace{0.5 truecm}
\noindent{{\bf Acknowledgements.}} This is a completely new version of
a previous work on the same subject. This new version was mainly
motivated by the many discussions with B. Toen during my stay in the
Laboratoire \'{E}mile Picard, Universit\'{e} Paul Sabatier
(Toulouse). My most sincere acknowledgement to him for his kindness. I
also ackownedge L. Crane, thanks to whom I became interested in the
representation theory of 2-groups, and with whom I discussed a lot on
the representations of such objects on Yetter's 2-category when we met
in the IST of Lisbon 
and in the Universit\'e de Montpelier. Finally, let me thank the DURSI
(Generalitat de Catalunya) and the Universitat Polit\`ecnica de 
Catalunya for their financial support and for allowing me to visit the
Laboratoire \'Emile Picard while I should have been teaching.

\bibliographystyle{plain}
\bibliography{2cat_repr_2grup_v2}

\begin{thebibliography}{10}

\bibitem{BC03}
J.~Baez and A.~Crans.
\newblock Higher-dimensional algebra vi: Lie 2-algebras.
\newblock {\em Preprint (arXiv: math.QA/0307263)}, 2003.

\bibitem{BL03}
J.~Baez and A.~Lauda.
\newblock Higher-dimensional algebra v: 2-groups.
\newblock {\em Preprint (arXiv: math.QA/0307200)}, 2003.

\bibitem{BM04}
J.~Barrett and M.~Mackaay.
\newblock Categorical representations of categorical groups.
\newblock {\em Preprint (arXiv: math.CT/0407463)}, 2004.

\bibitem{jB67}
J.~Benabou.
\newblock Introduction to bicategories.
\newblock In {\em Reports of the Midwest Category Seminar (LNM, volume 47)},
  pages 1--77. Springer, 1967.

\bibitem{fB94}
F.~Borceux.
\newblock {\em Handbook of Categorical Algebra 1}, volume~50 of {\em
  Encyclopedia of Mathematics and Its Applications}.
\newblock Cambridge University Press, 1994.

\bibitem{lB92}
L.~Breen.
\newblock Th\'{e}orie de schreier sup\'{e}rieure.
\newblock {\em Ann. Scient. \'{E}cole Norm. Sup.}, 25:465--514, 1992.

\bibitem{BRup}
R.~Brown and R.~Sivera.
\newblock {\em Non abelian algebraic topology}.
\newblock Draft of a planned book available on the web:
  www.bangor.ac.uk/~mas010/nonab-a-t.html.

\bibitem{BS76}
R.~Brown and C.B. Spencer.
\newblock $\mathcal{G}$-groupoids, crossed modules and the classifying space of
  a topological group.
\newblock {\em Proc. Kon. Akad. v. Wet.}, 79:296--302, 1976.

\bibitem{CC96}
P.~Carrasco and A.M. Cegarra.
\newblock Schreier theory for central extensions of categorical groups.
\newblock {\em Communications in Alg.}, 28:2585--2613, 1996.

\bibitem{CF94}
L.~Crane and I.~Frenkel.
\newblock Four dimensional topological quantum field theory, hopf categories
  and the canonical bases.
\newblock {\em J. Math. Phys.}, 35:5136--5154, 1994.

\bibitem{CS03}
L.~Crane and M.~Sheppheard.
\newblock 2-categorical representations of the poincar\'e group.
\newblock {\em arXiv: math.QA/0306440}, 2003.

\bibitem{CY03}
L.~Crane and D.~Yetter.
\newblock Measurable categories and 2-groups.
\newblock {\em preprint (arXiv: math.QA/0305176)}, 2003.

\bibitem{jE5}
J.~Elgueta.
\newblock A reconstruction theorem for finite 2-groups.
\newblock {\em (in preparation)}.

\bibitem{jE3}
J.~Elgueta.
\newblock A strict and totally coordinatized version of kapranov and
  voevodsky's 2-category {\bf 2vect}.
\newblock {\em Macquarie Mathematics Report, No. 860081}, 1986.

\bibitem{jE1}
J.~Elgueta.
\newblock Cohomology and deformation theory of monoidal 2-categories i.
\newblock {\em Adv. Math.}, 182:204--277, 2004.

\bibitem{HKK01}
K.H.~Hardie {\it et al.}
\newblock A homotopy bigroupoid of a topological space.
\newblock {\em Appl. Cat. Str.}, 9:311--327, 2001.

\bibitem{GPS95}
R.~Gordon {\it et al}.
\newblock {\em Coherence for tricategories}, volume 117 of {\em Mem. Am. Math.
  Soc.}
\newblock AMS, 1995.

\bibitem{FB02}
M.~Forrester-Barker.
\newblock Group objects and internal categories.
\newblock {\em Preprint (arXiv: math.CT/0212065)}, 2002.

\bibitem{GI01}
A.~Garzon and H.~Inassaridze.
\newblock Semidirect product of categorical groups, obstruction theory.
\newblock {\em Hom., Hom. and Applications}, 3:111--138, 2001.

\bibitem{JS93}
A.~Joyal and R.~Street.
\newblock Braided monoidal categories.
\newblock {\em Adv. Math.}, 102:20--78, 1993.

\bibitem{JS86}
A.~Joyal and R.~Street.
\newblock Braided monoidal categories.
\newblock {\em arXiv: math.CT/0406475}, 2004.

\bibitem{KV94}
M.~Kapranov and V.~Voevodsky.
\newblock 2-categories and zamolodchikov tetrahedra equations.
\newblock In {\em Proc. Sympos. Pure Math.}, volume 56(2), pages 177--260.
  American Mathematical Society, 1994.

\bibitem{aK76}
A.A. Kirillov.
\newblock {\em Elements of the theory of representations}.
\newblock Springer Verlag, 1976.

\bibitem{sM98}
S.~MacLane.
\newblock {\em Categories for the Working Mathematician}, volume~5 of {\em
  Graduate Texts in Mathematics}.
\newblock Springer, 1998.

\bibitem{PW04}
P.~Polesello and I.~Waschkies.
\newblock Higher monodromy.
\newblock {\em Arxiv.org: math.AT/0407507}, 2004.

\bibitem{SR72}
N.~Saavedra Rivano.
\newblock {\em Cat\'{e}gories tannakiennes}, volume 265 of {\em Lecture Notes
  in Mathematics}.
\newblock Springer-Verlag, 1972.

\bibitem{aR03}
A.~Rousseau.
\newblock Bicat\'{e}gories monoidales et extensions de gr-cat\'{e}gories.
\newblock {\em Homology, Homotopy and Appl.}, 5(1):437--547, 2003.

\bibitem{hxS75}
H.~X. Sinh.
\newblock {\em Gr-cat\'{e}gories}.
\newblock th\`{e}se de doctorat, Univ. Paris VII, 1975.

\bibitem{zT96}
Z.~Tamsamani.
\newblock Sur des notions de n-cat\'{e}gorie et n-groupoide non strictes via
  des ensembles multi-simpliciaux.
\newblock {\em K-Theory}, 16:51--99, 1999.

\bibitem{vT94}
V.~Turaev.
\newblock {\em Quantum invariants of knots and 3-manifolds}, volume~18 of {\em
  de Gruyter Studies in Mathematics}.
\newblock Walter de Gruyter, 1994.

\bibitem{dYnp}
D.~Yetter.
\newblock Categorical linear algebra: a setting for questions from physics and
  low-dimensional topology.
\newblock {\em Kansas State University, preprint}.

\bibitem{dY03}
D.~Yetter.
\newblock Measurable categories.
\newblock {\em arXiv: math.CT/0309185}, 2003.

\end{thebibliography}

\end{document}